\documentclass{amsart}
\usepackage{etex}
\usepackage[T1]{fontenc}
\usepackage[utf8]{inputenx}
\usepackage{lmodern}
\newcommand{\texorpdfstring}[2]{#1}
\usepackage{mathtools}
\usepackage{amssymb}
\usepackage{booktabs}
\usepackage{euscript}
\usepackage{tikz}
\usepackage{tikz-cd}
\usepackage{standalone}
\usepackage{braket}
\usepackage{varioref}
\vrefwarning % Don't give errors for references over pages.
\usepackage{xcolor}
\newtheorem{theorem}{Theorem}
\newtheorem{corollary}{Corollary}
\newtheorem{lemma}{Lemma}
\newtheorem{proposition}{Proposition}

\newcommand*{\defeq}{\mathrel{\vcenter{\baselineskip0.5ex \lineskiplimit0pt
                     \hbox{\scriptsize.}\hbox{\scriptsize.}}}%
                     =}

\usepackage{mathrsfs}
\newcommand{\pr}[1]{\ensuremath{\left(#1\right)}}

\newcommand{\hook}{\ensuremath{\mathbin{
\hbox{\vrule height1.4pt width4pt depth-1pt
\vrule height4pt width0.4pt depth-1pt}}}}

\newcommand{\pd}[2]{\ensuremath{\frac{\partial{#1}}{\partial{#2}}}}
\newcommand{\R}[1]{\ensuremath{\mathbb{R}^{#1}}}
\newcommand{\C}[1]{\ensuremath{\mathbb{C}^{#1}}}
\newcommand{\CP}[1]{\ensuremath{\mathbb{CP}^{#1}}}
\newcommand{\RP}[1]{\ensuremath{\mathbb{RP}^{#1}}}
\newcommand{\Gro}[2]{\widetilde{\operatorname{Gr}}^{#1}_{#2}}
\newcommand{\GrC}[2]{\operatorname{Gr}_{\C{},#2}^{#1}}
\newcommand{\SO}[1]{\operatorname{SO}\!\left({#1}\right)}
\newcommand{\GL}[1]{\operatorname{GL}\!\left({#1}\right)}

\newcommand{\gl}[1]{\mathfrak{gl}\!\left({#1}\right)}
\newcommand{\SL}[1]{\operatorname{SL}\!\left({#1}\right)}
\newcommand{\PSL}[1]{\mathbb{P}\SL{#1}}
\renewcommand{\sl}[1]{\mathfrak{sl}\!\left({#1}\right)}
\newcommand{\Lin}[2]{ {#1}^* \otimes {#2}}
\newcommand{\Cstrucs}[1]{\mathscr{J}_{#1}}
\newcommand{\Lm}[2]{\ensuremath{\Lambda^{#1}\!\left ( {#2} \right )}}
\newcommand{\nForms}[2]{\ensuremath{\Omega^{#1}\!\left ( {#2} \right )}}
\DeclareMathOperator{\Ad}{Ad}

\DeclareMathOperator{\Det}{Det}
\DeclareMathOperator{\tr}{tr}

\renewcommand{\Im}{\operatorname{Im}}
\DeclareMathOperator{\Aut}{Aut}
\definecolor{remark-color}{gray}{0.96}
\definecolor{remark-border-color}{gray}{0.98}
\theoremstyle{definition}

\newcommand{\ray}[1]{\ensuremath{\left[{#1}\right]}}
\newcommand{\Line}[1]{\ensuremath{\left<{#1}\right>}}
\newcommand{\OO}[1]{\ensuremath{\mathcal{O}\left({#1}\right)}}
\newcommand{\SatoSpace}{\ensuremath{\CP{}\left(V_{\C{}}\right)}} 
\newcommand{\Sato}{\ensuremath{S_V}}
\newcommand{\Gcircle}{\ensuremath{G_{\textit{circle}}}}

\newcommand{\M}{\ensuremath{M_V}}

\usepackage{colortbl}
\arrayrulecolor{gray!30}
\makeatletter
    \def\rulecolor#1#{\CT@arc{#1}}
    \def\CT@arc#1#2{%
      \ifdim\baselineskip=\z@\noalign\fi
      {\gdef\CT@arc@{\color#1{#2}}}}
    \let\CT@arc@\relax
\rulecolor{gray!30}
\makeatother

\pgfdeclarefunctionalshading{sphere}{\pgfpoint{-25bp}{-25bp}}{\pgfpoint{25bp}{25bp}}{}{
25 div exch
25 div exch
2 copy 
dup mul exch
dup mul add
1.0 sub
0.3 dup mul
-0.5 dup mul add
1.0 sub
mul abs sqrt
exch 0.3 mul add
exch -0.5 mul add
dup abs add 2.0 div 
0.6 mul 0.4 add
dup
0.4
}

\makeatletter
\renewcommand*\env@matrix[1][\arraystretch]{%
  \edef\arraystretch{#1}%
  \hskip -\arraycolsep
  \let\@ifnextchar\new@ifnextchar
  \array{*\c@MaxMatrixCols c}}
\makeatother
\newcommand{\tallmatrix}{1.5}

\usepackage{varwidth}
\makeatletter
\def\@tocline#1#2#3#4#5#6#7{\relax
  \ifnum #1>\c@tocdepth % then omit
  \else
    \par \addpenalty\@secpenalty\addvspace{#2}%
    \begingroup \hyphenpenalty\@M
    \@ifempty{#4}{%
      \@tempdima\csname r@tocindent\number#1\endcsname\relax
    }{%
      \@tempdima#4\relax
    }%
    \parindent\z@ \leftskip#3\relax \advance\leftskip\@tempdima\relax
    #5\leavevmode\hskip-\@tempdima #6\nobreak\relax
    ,~#7\par
    \endgroup
  \fi}
\makeatother

\begin{document}
\title[]{The Blaschke conjecture and great circle fibrations of spheres}
\author{Benjamin McKay}
\address{University College Cork \\
Cork, Ireland}
\email{b.mckay@ucc.ie}
\date{\today}
\thanks{Thanks to Daniel Allcock, Werner Ballmann and Karsten Grove.}
\begin{abstract}
We construct an explicit diffeomorphism taking any ``nondegenerate'' fibration of a sphere by great circles into the Hopf fibration. 
The diffeomorphism is a local differential invariant, algebraic in derivatives. 
\end{abstract}
\maketitle
\begin{center}
%\begin{varwidth}{\textwidth}
\tableofcontents
%\end{varwidth}
\end{center}
\section{Introduction}
\begin{quote}
``Notice that the classification of fibrations  of spheres by great circles
is an interesting but almost untouched subject\ldots''
Arthur L. Besse \cite{Besse:1978} pg. 135.
\end{quote}
In studying the Blaschke conjecture \cite{Besse:1978} and in the theory of elliptic partial differential equations \cite{McKay:2001}, 
one encounters fibrations of the standard round sphere
by great circles.\footnote{Gluck, Warner and Yang \cite{GluckWarnerYang:1983} is a nice introduction to great sphere fibrations of spheres.}
The best known example is the Hopf fibration
\[
\begin{tikzcd}
S^1 \arrow{r} & S^{2n+1} \arrow{d} \\
           & \CP{n}
\end{tikzcd}
\]
given by the circle action
\[
e^{i \theta} \left ( z_0 , \dots , z_n \right )
=
\left ( e^{i \theta} z_0, \dots, e^{i \theta} z_n \right )
\]
on the unit sphere inside \(\C{n+1}\). 

Yang's conjecture \cite{Yang:1990} states that for any smooth great circle fibration of a sphere, there is a diffeomorphism of the sphere carrying it to the Hopf fibration.
Sato proved Yang's conjecture for PL-homeomorphism rather than diffeomorphism \cite{Sato:1984}.
Yang \cite{Yang:1990,Yang:1993} stated that Sato's proof is flawed, and proved his conjecture for spheres of dimension at least seven, missing the 3-sphere and 5-sphere, although he claimed that his proof works on the 5-sphere modulo the Poincar{\'e} conjecture.
Ballmann and Grove (personal communication from Grove) stated that Yang's proof is flawed.
Gluck \& Warner \cite{GluckWarner:1983} proved Yang's conjecture for the 3-sphere.
In an earlier version of this paper, I attempted to prove Yang's conjecture in any dimension.
Yang's proof of his conjecture is quite difficult, involving a mixture of differential geometry and surgery, and employing the \(h\) and \(s\) cobordism theorems and the signature theorem; the diffeomorphism is not explicit. 
For great circle fibrations which are, in a suitable sense, nondegenerate, I will give an explicit diffeomorphism, which can be  written in algebraic functions of the first and second derivatives of the functions cutting out the great circle fibration, in any local coordinates.
Our diffeomorphism is ``linear'' on each fiber, i.e. a linear fractional transformation of each circle.
The nondegeneracy is satisfied by the Hopf fibration, and by every great circle fibration near the Hopf fibration, but I do not know if there are degenerate great circle fibrations.

A \emph{Blaschke manifold} is a Riemannian manifold whose injectivity radius equals its diameter.
Every Blaschke manifold has the cohomology of a unique compact rank one symmetric space, its \emph{model}.
The Blaschke conjecture claims that every Blaschke manifold is isometric to its model, up to a constant rescaling.
The conjecture above has the corollary that every Blaschke manifold modelled on a complex projective space is homeomorphic to its model.
Yang \cite{Yang:1990} claims diffeomorphic, but I cannot see how to prove such a result, even after a careful reading of Yang's paper, nor can Karsten Grove, and even modulo Yang's conjecture.
An abridged version of the paper you have before you has appeared in the \emph{American Journal of Mathematics}, with an attempt at a proof of Yang's conjecture.
In 2016, I rewrote parts of this paper to clarify the proofs, correct errors discovered by Karsten Grove and Werner Ballmann, and to improve the citations. 
Below I give credit to Reznikov for discovery of invariants of great circle fibrations, and I remove the claim that every Blaschke manifold modelled on a complex projective space is diffeomorphic, or even homeomorphic, to its model.

\section{Overview}

Section~\ref{sec:ElemTop} presents a review of some results on the topology of circle fibrations of spheres, which are intended for the reader's interest, but are not employed subsequently.

The long section~\ref{sec:MovingFrame} develops  a description of a great circle fibration in terms of local data, following Cartan's method of the moving frame. This associates to each great circle fibration of a sphere a principal bundle over the sphere, with a collection of differential forms on it representing
the local data of the great circle fibration.
The idea is to succesively determine subbundles of this bundle, by algebraic
equations among those differential forms.
This process is also a part of Cartan's method. 
In order to make it work, one needs to find algebraic relations among the coefficients of the various differential forms, and then show that the functions
expressing these relations (called the \emph{torsion functions}), which are differential invariants of the great circle fibration, satisfy regularity hypotheses. 
These regularity hypotheses are strong enough that the subset on which the torsion satisfies some algebraic condition is a submanifold.
Recall that the original bundle is a principal bundle---we will
see that these coefficients vary according to an action of the structure group of that bundle. 
This will force the subbundle on which the torsion satisfies an appropriately chosen algebraic equation to be a principal subbundle.

The first torsion to show up is essentially an endomorphism of the normal bundle of each great circle. 
(It is an endomorphism twisted by a real line bundle, but this is irrelevant.) 
Subsection~\ref{subsec:Follow} shows that this endomorphism satisfies an ordinary differential equation as we move along the great circles (using the Levi--Civita connection to differentiate). 
This ordinary differential equation is explicitly integrated, and we find that consistent initial data for it which will remain finite as we tranverse the great circle consists in endomorphisms with no real eigenvalues. 

This leads to a digression: section~\ref{subsec:LinTrans} shows that an endomorphism of a vector space which has no real eigenvalues determines invariantly a complex structure which commutes with it. 
Applied to the endomorphism of the normal bundle, we obtain an almost complex structure on the base manifold of the great circle fibration. 
In subsection~\vref{subsec:Back}, I change bases for the differential forms, splitting forms into complex linear and conjugate linear parts, and find that the Cayley transform of our endomorphism is a complex linear matrix with eigenvalues in the unit disk. 
The structure group of our principal bundle acts on this matrix, moving its spectrum around. 
Section~\vref{subsec:LFT}  finds that one can normalize it to have vanishing trace.

In section~\vref{subsec:Reduce}, an elementary step enables reduction of the structure group of our original principal bundle to a group \(\Gamma_0\). 
This group is the same group which appears as structure group for the Hopf
fibration. 
But this is exactly the isotropy group of a point of the base manifold in the Hopf fibration, signalling an end to the method of the moving frame, since there can in general be no further reduction, as the Hopf fibration admits no invariant subbundle contained in this one.

Section~\vref{sec:Analogy} shows that the structure equations satisfied by our differential forms now appear similar to those of complex projective space.
Section~\vref{sec:Homog} proves that the Hopf fibration is precisely the great circle fibration with maximal symmetry group.

Section~\vref{sec:Embedding} turns to another description of a great circle fibration: every great circle in a sphere spans a 2-plane in the ambient vector space containing the sphere, so a family of great circles is a family of 2-planes.
More precisely, the base of a great circle fibration is an embedded compact connected submanifold of the Grassmannian of 2-planes. 
Section~\vref{sec:Charac} characterizes the submanifolds of the Grassmannian
which represent great circle fibrations as being precisely those which 
are compact, connected, and \emph{elliptic}.
Ellipticity is a local condition on an immersed submanifold of the Grassmannian; ellipticity is a first order partial differential inequality. 
I believe that this inequality satisfies the \(h\)-principle of Gromov; I explain in section~\vref{sec:H} why one is apparently unable to use Gromov's techniques to prove the \(h\)-principle. 
The space of all great circle fibrations, viewed as the space of all compact,
connected, elliptic submanifolds of the Grassmannian, is straightforward
to parameterize locally, and shows itself as an infinite dimensional
manifold. 
Its topology is unknown; for example it is not known to be connected.

In section~\vref{sec:Osculating}, the geometry of the principal bundle
we have constructed is used to determine an osculating complex structure at
each point of the base manifold of the great circle fibration. 
This is a complex structure whose associated Hopf fibration has base manifold
sitting inside the Grassmannian just touching the base manifold
of our given great circle fibration, and is in some sense the best
approximating Hopf fibration.

All of the theory developed in this article is based on the conviction
that great circle fibrations provide a very natural mechanism for
deforming complex geometry. 
This article is a contribution to the microlocal theory of such deformations. 
So a great circle fibration of a sphere should be thought of as a kind of 
nonlinear complex structure on the vector space containing the sphere. 
The base of the great circle fibration is a kind of deformed complex projective space. 
The next step is to define complex hyperplanes in that space. 
We do this by looking at hyperspheres in our sphere, and asking for the family of great circles from our fibration which lie entirely inside the hypersphere.
We prove that this ``hyperplane'' is a smooth submanifold in the base manifold, of real codimension 2.
It is not an almost complex submanifold in general.

Hajime Sato \cite{Sato:1984} attempted to prove part of the topological Blaschke conjecture using a map, which was perhaps not well defined \cite{Yang:1990}.
Therefore the next step is to define Sato's map. 
This map embeds the base of a great circle fibration into a complex projective space of much higher dimension.
A diffeomorphism from the base manifold to a complex projective space is constructed out of the Sato map, essentially using linear projections, assuming a nondegeneracy condition.
We carry over this map into an isomorphism of the given great circle fibration with the Hopf fibration, achieving our main result.

The analogy between great circle fibrations of spheres and complex
structures of vector spaces is significantly strengthened in section~\vref{sec:NonlinearJ}.
We describe the notion of twisted complex structure, or nonlinear \(J\). 
To each great circle fibration, we assign a nonlinear \(J\), using local invariants. 
To each nonlinear \(J\), we associate a great circle fibration. 
However, the two concepts are not equivalent: rather the set of nonlinear \(J\) form an infinite dimensional fiber bundle with contractible fibers over the
infinite dimensional manifold of great circle fibrations. 
There is a canonical section of this bundle. 
This picture tells us that these nonlinear \(J\) are really superfluous, and that the real object of our theory, generalizing the concept of complex structure, is the great circle fibration.

\section{Elementary topology}{\label{sec:ElemTop}

A sphere of even dimension does not admit a circle fibration. 
A circle fibration determines a unit vector field up to sign, and therefore determines a unit vector field on a 2-1 cover. 
Because spheres of dimension two or more are simply connected, this determines two unit vector fields, and therefore shows the vanishing
of the Euler characteristic.
The Euler characteristic of a sphere is zero just when the dimension of the sphere is odd.
Therefore we restrict our attention to the odd dimensional spheres.

Fix a smooth foliation of  \(S^{2n+1}\) by great circles.
The foliation is a principal circle bundle, say \(S^1 \to S^{2n+1} \to X,\) since the fibers can be consistently oriented, and then we can apply rotations by various angles (measured in the usual Riemannian geometry  on the sphere \(S^{2n+1}\)) to implement a circle action. 
(We  will see this from another point of view below.)
The compactness of \(S^{2n+1}\) forces \(X\) to be compact.

Using the orientations of the circles and of the sphere, we have a quotient orientation on \(X\).
Chasing through the relevant exact sequences \cite{DNF:1985} p. 230, \cite{McCleary:1985} p. 134, the homotopy groups, Whitehead products and cohomology ring of \(X\) (with coefficients in any  ring; see McCleary \cite{McCleary:1985}, page 134) are those of \(\CP{n}\). 
(These authors carry out their calculations with the hypotheses that the base manifold \(X\) is \(\CP{n}\) and that the fibration is the Hopf fibration. However, they do not use these hypotheses. 
All they require is that the total space be a homotopy sphere, and that the fiber be a circle.) 
Turning to  characteristic classes \cite{McCleary:1985} p. 199, follow the Leray spectral sequence to see that the Chern class of the bundle \(S^{2n+1} \to X\) (i.e. the transgression of the \(S^1\) generating class), which we write as \(-[H]\), generates the cohomology.  
Keeping track of signs (using the Hopf fibration as our guide), find that \([H]^n\) is Poincar{\'e} dual to the fundamental class \([X].\)

\section{The moving frame}\label{sec:MovingFrame}

\subsection{Structure equations of a flat projective structure}

For my purposes, the sphere \(S^{2n+1}\) and the real projective space \(\RP{2n+1}\) are equally reasonable spaces to work on, since a great circle fibration of a sphere is the same thing as a fibration of the [antipodal quotient] projective space by projective lines: choose the sphere.
The group of symmetries of a fiber bundle is always infinite dimensional, but in our case we want the concept of great circle to be preserved. 
\begin{theorem}
Every invertible (not necessarily continuous) map of real projective space of dimension \(2\) or more to itself taking projective lines to projective lines is the action of a projective linear transformation.
\end{theorem}
\begin{proof}
The result is well known for projective planes \cite{Hartshorne:1967} p. 48 theorem 3.13.
Assume by induction that \(n \ge 3\) and that we have proven the result for all dimensions less than \(n\). 
Take an invertible map \(f \colon \RP{n} \to \RP{n}\) taking lines to lines.
Recall that a set of points of projective space is in general position if any \(k\) of them lie in a unique projective subspace of dimension \(k-1\) for an integer \(k\) with \(1 \le k \le n\).
Note that the projective general linear group acts transitively on \((n+2)\)-tuples of points in general position.
Take \(n+2\) vectors \(v_0, v_1, \dots, v_n, v_{n+1}\) in \(\R{n+1}\) whose images in \(\RP{n}\) are in general position.
For example, we can assume that \(v_0=e_0, \dots, v_n=e_n\) are the standard basis vectors and that \(v_{n+1}=\sum e_i\).
Composing \(f\) with a projective linear transformation, we can assume that \(f\) fixes the images of these points in projective space.
So \(f\) leaves invariant the projective subspaces of all dimensions spanned by any subset of these points.
By induction on dimension, \(f\) is projective linear on each of those projective subspaces.
In particular, any \(n\) of our points lie in a unique projective subspace of dimension \(n-1\), on which we find that \(f\) is projective linear for a linear transformation with the various \(v_i\) as eigenvectors:
\[
f\left[ \sum_{i \ne j} a_i e_i \right]
=
\left[ \sum_{i \ne j} a_i \lambda_i e_i \right]
\]
for some numbers \(\lambda_i\).
Allowing \(a_k = 0\) for some \(k\), we find that the eigenvalue \(\lambda_i\) is independent of \(j\).
By the same token
\[
f\left[a e_0 + \sum_j e_j\right]
=
\left[a \lambda_0 e_0 + \mu \sum_j e_j\right],
\]
for some eigenvalue \(\mu \ne 0\).
Setting \(a=-1\), 
\[
f\left[\sum_{j>0} e_j\right]
=
\left[\left(\mu-\lambda_0\right)e_0 + \mu e_1 + \dots + \mu e_n\right].
\]
Linearity also gives
\[
f\left[\sum_{j>0} e_j\right]
=
\left[\lambda_1 e_1 + \dots + \lambda_1 e_n\right].
\]
Therefore all eigenvalues are equal, i.e. \(f\) fixes each \(\left[e_i\right]\) in projective space, and \(f\) fixes all points on all of the projective subspaces spanned by any set of these \(\left[e_i\right]\), in particular all points of all coordinate hyperplanes are fixed.
Any point of projective space is the intersection point of a pair of lines through points of coordinate hyperplanes.
\end{proof}

\begin{corollary}
Any continuous bijection of a sphere taking great circles to great circles is the action of a linear transformation.
\end{corollary}
\begin{proof}
Any such map of the sphere takes antipodal points to antipodal points, since these are intersection points of great circles.
So it quotients to a map of real projective space taking lines to lines, a projective linear automorphism.
After composing with a linear transformation, we arrange that this projective transformation is the identity, so our continuous bijection induces the trivial projective automorphism, so each point is preserved or reflected to its antipode.
By continuity, our bijection is either the identity or minus the identity.
\end{proof}

Let \(V=\R{2n+2}\) have basis \(e_0,\dots,e_{2n+1}\).
Let \(G = \SL{V}\) act on the sphere \(S^{2n+1} = \left(V\backslash{0}\right)/\R{+}\).
Let \(\mathfrak{g}\) be the Lie algebra of \(G\).
Write \(\ray{v}=\R{+}v\) for the ray through a  vector \(v \in V\backslash 0\).
Let \(G_0\) be the stabilizer of \(\ray{e_0}\), i.e. the group of matrices of the form
\[
\begin{pmatrix}
g^0_0 & g^0_{\nu} \\
0 & g^{\mu}_{\nu}
\end{pmatrix}
\]
with real entries satisfying \(g^0_0 > 0\),
\(g^{0}_{0} \det g^{\mu}_{\nu} = 1\), and Greek indices
\(\mu, \nu = 1,\dots,2n+1\). For reference,
our index conventions in this paper are:
\begin{align*}
\mu,\nu,\sigma&=1,\dots,2n+1 \\
i,j,k&=2,\dots,2n+1 \\
p,q,r&=1,\dots,n \\
P,Q,R&=0,\dots,n.
\end{align*}
Let \(\mathfrak{g}_0\)
be the Lie algebra of \(G_0\). 
The sphere \(S^{2n+1}\) is the quotient
\(G/G_0\) via the right action of \(G_0\).
More concretely, the map \(G \to S^{2n+1}\)
is the map \(g \mapsto \ray{g e_0}.\)

We will follow {\'E}lie Cartan's method of the moving
frame \cite{Clelland:2016}.
The left invariant Maurer--Cartan 1-form 
\[
\omega = g^{-1} \, dg \in \nForms{1}{G} \otimes \mathfrak{g}
\]
satisfies
\(
d \omega = - \omega \wedge \omega.
\)
Our subgroup \(G_0\)
acts on the right on \(G\), thus not preserving
the Maurer--Cartan 1-form, but instead
if \(R_{g_0}\) is the right action on \(G\)
of an element \(g_0 \in G_0\) then
\[
R_{g_0}^* \omega = \Ad^{-1}_{g_0} \omega.
\]
Let us divide \(\omega\) into 1-forms
according to
\[
\omega =
\begin{pmatrix}
\omega^0_0 & \omega^0_{\nu} \\
\omega^{\mu}_0 & \omega^{\mu}_{\nu}
\end{pmatrix}.
\]
The 1-forms \(\omega^{\mu}_0\), which we will
write as \(\omega^{\mu}\), are semibasic
for the fiber bundle map \(G \to S^{2n+1}\),
which we see because they are linearly
independent on \(G\) but vanish on \(G_0\)
and therefore on the left translates of \(G_0\).
Their exterior derivatives are
\(
d \omega^{\mu} = - \gamma^{\mu}_{\nu} \wedge \omega^{\nu}
\)
where we define
\[
\gamma^{\mu}_{\nu} \defeq \omega^{\mu}_{\nu} - \delta^{\mu}_{\nu} \omega^0_0.
\]
Their exterior derivatives are
\[
d \gamma^{\mu}_{\nu}
=
- \gamma^{\mu}_{\sigma} \wedge \gamma^{\sigma}_{\nu}
+ \left ( \delta^{\mu}_{\nu} \omega_{\sigma} 
+ \omega_{\nu} \delta^{\mu}_{\sigma} \right ) \wedge \omega^{\sigma}
\]
where we write \(\omega_{\mu}\) for \(\omega^0_{\mu}\).
Our structure equations can now be rewritten
as
\begin{align*}
d \omega^{\mu} &= - \gamma^{\mu}_{\nu} \wedge \omega^{\nu} \\
d \gamma^{\mu}_{\nu} &= - \gamma^{\mu}_{\sigma} \wedge
\gamma^{\sigma}_{\nu}
+ \left ( \delta^{\mu}_{\nu} \omega_{\sigma} 
+ \omega_{\nu} \delta^{\mu}_{\sigma} \right ) \wedge \omega^{\sigma} \\
d \omega_{\mu} &= \gamma^{\nu}_{\mu} \wedge \omega_{\nu}.
\end{align*}
These are the structure equations of a flat projective
structure \cite{Kobayashi/Nagano:1964}.

\subsection{Canonically defined vector bundles on a manifold with flat projective structure}{\label{subsec:VBS}
It is not essential to work out the theory of invariantly defined vector bundles on the sphere determined by the projective structure, but it makes clearer the interpretation of the invariantly defined vector bundles which we will produce from a great circle fibration in subsection~\vref{subsec:Later}.

We can now take any representation of the group
\(G_0,\) say \(\rho \colon G_0 \to \GL{W}\), and
use it to define a vector bundle \(\tilde{W} \to S^{2n+1}\)
by 
\[
\tilde{W} = G \times^{G_0} W = \left(G \times W\right)/G_0
\]
where the quotient is taken by the \(G_0\) action
\[
(g,w)g_0 = \left(gg_0,\rho\left(g_0\right)^{-1}w\right)
\]
for
\[
g_0 \in G_0, g \in G, w \in W.
\]
A section of \(\tilde{W}\to S^{2n+1}\)
is precisely a \(G_0\)-equivariant map
\( f \colon G \to W. \)
If we pick a basis \(w_{\alpha}\)
of \(W\), \(f\) has components \(f^{\alpha}\). 
Write \(\rho: \mathfrak{g}_0 \to \gl{W}\)
for the Lie algebra morphism induced
by our morphism \(\rho \colon G_0 \to \GL{W}\) of Lie groups.
The differential of \(f^{\alpha}\) is
\[
df^{\alpha} + \rho
\begin{pmatrix}
\omega^0_0 & \omega^0_{\nu} \\
0 & \omega^{\mu}_{\nu}
\end{pmatrix}
^{\alpha}_{\beta}
f^{\beta} = f^{\alpha}_{\mu} \omega^{\mu}
\]
for some functions \(f^{\alpha}_{\mu}\) on \(G\),
or in other words
\[
df + \rho(g_0^{-1}dg_0)f = \nabla f \omega
\]
where
\[
\nabla f \colon G \to e_0^{\perp} \otimes W
\]
is the covariant derivative of the section \(f\).
We will say that \(\rho\) \emph{solders}
the bundle \(\tilde{W}\). Since
all of the Lie groups I will
employ in this article are connected,
the Lie algebra representation
will suffice for our purposes
to identify the group representation, 
and we will
usually only indicate the
Lie algebra representation,
saying that the bundle is
soldered by the expression
\[
\rho
\begin{pmatrix}
\omega^0_0 & \omega^0_{\nu} \\
0 & \omega^{\mu}_{\nu}
\end{pmatrix}
^{\alpha}_{\beta}.
\]

For example, we define \(\OO{-1}\) to be
the bundle  
\(\OO{-1}=\tilde{W}\)
where \(W=\Line{e_0}\subset V\) is the span of \(e_0\)
in \(V\). It is not difficult
to see that if we were to quotient our
sphere down to the underlying real projective
space, then this bundle \(\OO{-1}\)
would become the algebraic geometer's
usual universal line bundle.
Since \(G_0\) preserves an orientation
in this line \(\Line{e_0}\), the bundle \(\OO{-1}\)
has oriented fibers. 
Write any section \(\sigma\) of \(\OO{-1}\)
(perhaps only defined on an open subset
of the sphere) as
\[
\sigma\left(g\ray{e_0}\right)
= f(g) g e_0.
\]
Then \(f \colon G \to \R{}\) and
\(
df + \omega^0_0 f = f_{\mu} \omega^{\mu}
\)
for some functions \(f_{\mu} \colon G \to \R{}\).
So \(\OO{-1}\) is soldered by \(\omega^0_0\).

Define 
\(
\OO{1} \defeq \OO{-1}^*
\)
and similarly
\[
\OO{p} \defeq \OO{1}^{\otimes{p}} = \left(\Line{e_0}^*\right)^{\otimes p}.
\]
Then \(\OO{p}\) is soldered by \(-p\omega^0_0\).

We will also want to consider the bundle
\(\tilde{V}\) soldered by 
\[
\begin{pmatrix}
\omega^0_0 & \omega^0_{\nu} \\
0 & \omega^{\mu}_{\nu}
\end{pmatrix}
\]
(i.e. \(\rho\) is the identity).
This is a trivial bundle, since 
any element \(v \in V\) gives rise
to a global section \(f_v(g)=g^{-1}v\)
of \(\tilde{V} \to S^{2n+1}\). However,
there is no trivialization invariant
under \(\SL{V}\). Therefore we will prefer
to consider \(\tilde{V}\) as
a vector bundle.
The sections of this bundle
correspond to functions \(f \colon G \to V\)
so that in terms of our usual basis of \(V\)
\[
d
\begin{pmatrix}
f^0 \\
f^{\mu}
\end{pmatrix}
+
\begin{pmatrix}
\omega^0_0 & \omega^0_{\nu} \\
0 & \omega^{\mu}_{\nu}
\end{pmatrix}
\begin{pmatrix}
f^0 \\
f^{\nu}
\end{pmatrix}
=
\begin{pmatrix}
f^0_{\nu} \\
f^{\mu}_{\nu}
\end{pmatrix}
\omega^{\nu}.
\]
Take a fixed vector \(v \in V\) and let \(f = g^{-1}v\), so that
\(
df = - \omega f
\)
or
\[
d
\begin{pmatrix}
f^0 \\
f^{\mu}
\end{pmatrix}
+
\begin{pmatrix}
\omega^0_0 & \omega^0_{\nu} \\
\omega^{\mu}_0 & \omega^{\mu}_{\nu}
\end{pmatrix}
\begin{pmatrix}
f^0 \\
f^{\nu}
\end{pmatrix}
=
0.
\]
We see the covariant derivatives when
we write it as
\[
d
\begin{pmatrix}
f^0 \\
f^{\mu}
\end{pmatrix}
+
\begin{pmatrix}
\omega^0_0 & \omega^0_{\nu} \\
0 & \omega^{\mu}_{\nu}
\end{pmatrix}
\begin{pmatrix}
f^0 \\
f^{\nu}
\end{pmatrix}
=
-
\begin{pmatrix}
0 \\
f^{0} \delta^{\mu}_{\nu}
\end{pmatrix}
\omega^{\nu}.
\]

If a vector bundle \(\tilde{W}\) is soldered by
\(\rho^{\alpha}_{\beta}\)
then its dual \(\widetilde{W^*}=\tilde{W}^*\)
is soldered by \(-\rho^{\beta}_{\alpha}\),
i.e. the negative transpose.

When we add two representations,
say \(U\) and \(W\), with bases
\(u_{\alpha}\) and \(w_{M}\)
we obtain a representation \(U \oplus W\)
with basis \(z_I\) where
\(I\) runs over first the \(\alpha\)
indices and then the \(M\) indices.
The matrix elements
of the Lie algebra representation
(or the Lie group) on the sum are
\[
\rho^I_J = \rho^{\alpha}_{\beta} \delta^I_{\alpha} \delta^{\beta}_J
+ \delta^I_{M} \delta^{N}_{J} \rho^M_N.
\]
Similarly under tensor product,
the new index \(I\) runs over pairs \((\alpha,M)\)
and gives matrix elements
of the Lie algebra representation
(not the same as the group representation)
\[
\rho^I_J = \rho^{\alpha}_{\beta} \delta^{M}_{N}
+ \delta^{\alpha}_{\beta} \rho^{M}_{N}
\]
where
\[
I=\left(\alpha,M\right), \ J = \left(\beta,N\right).
\]
For a quotient representation,
\(W/U\) where \(U\) is an invariant subspace
of \(W\), we have indices \(\alpha,\beta\) for \(U\),
\(M,N\) for \(W/U\), and \(I,J\) for \(W\), and 
the relation
\[
\left ( \rho^I_J \right )
=
\begin{pmatrix}
\rho^{\alpha}_{\beta} & \rho^{\alpha}_{N} \\
0 & \rho^{M}_{N}
\end{pmatrix}
\]
so that we can dig out the matrix
elements of the quotient from
those of the original representation.

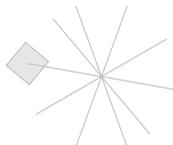
\begin{center}
\begin{tikzpicture}[rotate=50]
\fill[gray!20,draw=gray!50] 
	({cos(2*60)-.2},{sin(2*60)-.2}) 
	-- 
	({cos(2*60)+.2},{sin(2*60)-.2}) 
	-- 
	({cos(2*60)+.2},{sin(2*60)+.2}) 
	-- 
	({cos(2*60)-.2},{sin(2*60)+.2}) -- cycle; 
\foreach \i in {80,120,...,240}
{
	\draw[gray!50] ({cos(\i+180)},{sin(\i+180)}) -- ({cos(\i)},{sin(\i)});
}
\end{tikzpicture}
\end{center}

\begin{lemma}\label{lm:HotCereal}
The tangent bundle of the sphere is soldered
by
\(
\gamma^{\mu}_{\nu}
\)
and is canonically and
\(\SL{V}\)-equivariantly isomorphic to
\[
\OO{1}\otimes\left(\tilde{V}/\OO{-1}\right).
\]
\end{lemma}
\begin{proof}
Lets write \(\Line{e_0}\) for the 
line in \(V\) through \(e_0\).
The tangent space to the sphere
\(S^{2n+1}\) at \(\ray{e_0}\)
is given by
\[
0 \to T_{e_0} \ray{e_0} = \Line{e_0} \to T_{e_0} V = V \to 
T_{\ray{e_0}} S^{2n+1} =V/\Line{e_0} \to 0.
\]
But if we change the choice of the
point \(e_0\), by a positive multiple,
then we have to rescale \(V\) by
this multiple, and rescale \(\Line{e_0}\)
by the same multiple, while we don't
rescale the sphere at all. Therefore
this description of the tangent
space to the sphere is certainly
not \(G_0\)-equivariant, and it will only
become scale invariant if we rewrite
it as
\[
0 \to \Line{e_0}^* \otimes \Line{e_0} 
\to \Line{e_0}^* \otimes V \to T_{\ray{e_0}} S^{2n+1} \to 0 
\]
with the first map just an inclusion, and
the second map defined by
\[
\xi \otimes v \mapsto \xi\left(e_0\right)v/\Line{e_0}
\in V/\Line{e_0} = T_{\ray{e_0}}S^{2n+1}.
\]
We have to check that this map is
\(G_0\)-equivariant, invariant under positive rescaling,
because this cancels out from each
factor. Under linear transformations
\(g_0 \in G_0\), which take \(e_0\) to \(g^0_0 e_0\)
with \(g^0_0 > 0\),
\[
\begin{tikzcd}[ampersand replacement=\&]
\xi \otimes v \arrow{r} \arrow{d} \& \xi(e_0) \otimes v / \Line{e_0} \arrow{d} \\ 
\left(g_0^* \xi\right) \otimes \left(g_0 v\right) \arrow{r} \&
\xi\left(e_0\right) g_0 v_0/ \Line{e_0}
\end{tikzcd}
\]
we see that the map is \(G_0\)-equivariant,
so that
\[
T_{\ray{e_0}} S^{2n+1} = 
\left(\Line{e_0}^* \otimes V\right)
/\left(\Line{e_0}^* \otimes \Line{e_0}\right)
=
\Line{e_0}^* \otimes \left(V/\Line{e_0}\right)
\]
and therefore
\[
TS^{2n+1} = 
\widetilde{\Line{e_0}^*} \otimes \left(\tilde{V}/\tilde{\Line{e_0}}\right)
=
\OO{1}\otimes\left(\tilde{V}/\OO{-1}\right)
=
\OO{1}\otimes\tilde{V}/\OO{0}.
\]

Finally, we will consider 
the soldering of the tangent
bundle. The representation \(\omega^{\mu}_{\nu}\)
solders \(\tilde{V}/\OO{-1}\) and 
\(-\omega^0_0\) solders \(\OO{1}.\)
Therefore the tensor product
\[
\OO{1}\otimes\left(\tilde{V}/\OO{-1}\right)
\]
is soldered by
\[
\left( - \omega^0_0 \right ) \delta^{\mu}_{\nu}
+ \delta^0_0 \omega^{\mu}_{\nu} 
=
\omega^{\mu}_{\nu} - \delta^{\mu}_{\nu} \omega^0_0.
\]
\end{proof}

\begin{corollary}
\[
T^*S^{2n+1} = \OO{0}^{\perp} \subset \OO{-1} \otimes \tilde{V}^*
\]
is soldered by 
\(
- \left ( \gamma^{\mu}_{\nu} \right )^t.
\)
while 
\[
\Det TS^{2n+1} = \OO{-(2n+2)}
\]
is soldered by
\(
-(2n+2) \omega^0_0.
\)
\end{corollary}

\subsection{Structure equations of a geodesic
foliation in a flat projective structure}

 Let the Roman indices \(i,j,k,l\) run from \(2\)
to \(2n+1\). Inside \(G\) we have a subgroup
\(\Gcircle\)
consisting of matrices of the form
\[
\begin{pmatrix}
g^0_0 & g^0_1 & g^0_j \\
g^1_0 & g^1_1 & g^1_j \\
0 & 0 & g^i_j \\
\end{pmatrix}
\]
with \(\det = 1\) and 
\[
\det\begin{pmatrix}
g^0_0 & g^0_1 \\
g^1_0 & g^1_1
\end{pmatrix}
> 0. 
\]
This is the subgroup of all elements 
of \(G\)
which preserve the oriented plane spanned
by \(e_0\) and \(e_1.\) Under the 
map \(G \to S^{2n+1}\), this subgroup
projects to the oriented great circle which
is the image of the \(e_0, e_1\) plane.
This subgroup satisfies
\(
\omega^i = \gamma^i_1 = 0
\)
as do all of its left translates,
by left invariance of the Maurer--Cartan
1-form.
Thus the geodesics (the great circles) 
are the curves in the sphere which are the projections
of the integral manifolds of
\(
\omega^i = \gamma^i_1 = 0
\) 
in \(G\). The manifold \(G/\Gcircle\)
is the manifold of all oriented great circles
on the sphere \(S^{2n+1}\).

Similarly, any pointed great circle on the sphere is the projection of a left translate of the subgroup \(G_1\) of matrices of the form
\[
\begin{pmatrix}
g^{0}_{0} & g^0_1 & g^0_j \\
0 & g^1_1 & g^1_j \\
0 & 0 & g^i_j
\end{pmatrix}
\]
where \(g^0_0,g^1_1>0\) and 
\(
g^0_0 g^1_1 \det\left(g^i_j\right) = 1
\)
(the subgroup of \(G\) preserving not only the oriented \(e_0,e_1\) plane, but also fixing the point \(e_0\), up to positive factor).
The left translates of \(G_1\) are precisely the leaves of the foliation of \(G\)
given by the equations
\(
\omega^1 = \omega^i = \gamma^i_1 = 0.
\)
Hence the quotient space \(G/G_1\) is the space of pointed great circles.

Starting with a foliation \(F\) by curves, look
inside \(G\) and consider the subbundle 
\(B_1\) whose fiber above any 
point \(x \in S^{2n+1}\) consists of 
the linear maps \(g \colon V \to \R{2n+2}\)
which identify our point \(x\) of the sphere, i.e. 
a ray in \(V\), with a given ray in \(\R{2n+2}\),
say the ray through \(e_0\), and
which identify the tangent line to the
\(F\) curve through \(x\) with a given 2-plane
in \(\R{2n+2}\), say the span of \(e_0,e_1\).
The map \(B_1 \to S^{2n+1}\) is
the pullback
\[
\begin{tikzcd}[ampersand replacement=\&]
B_1 \arrow{r} \arrow{d} \& G \arrow{d} \\
S^{2n+1} \arrow{r} \& G/G_1
\end{tikzcd}
\]
above the map \(S^{2n+1} \to G/G_1\)
(taking any point \(p \in S^{2n+1}\) 
to the great circle which is tangent at \(p\)
to the leaf of \(F\) through \(p\))
and so \(B_1 \to S^{2n+1}\)
is a smooth right principal \(G_1\) bundle.
The fibers of the bundle \(B_1\) are 
cut out by the equations 
\(\omega^1=\omega^i=\gamma^i_1=0.\)
But the \(\omega^1,\omega^i\) are 
linearly independent 1-forms on the bundle \(B_1\) 
(since we have made no restriction on motions in
the base manifold, the sphere \(S^{2n+1}\)).
Therefore on this bundle \(B_1\)
\(
\gamma^i_1 = t^i_1 \omega^1 + t^i_j \omega^j
\)
for some functions \(t^i_1\) and \(t^i_j\).
The equation
\(
\omega^i = 0
\)
cuts out a foliation \(F_1\) of \(B_1\),
by the Frobenius theorem. Indeed 
the leaves of \(F_1\) are precisely
the preimages in \(B_1\) of the leaves 
of \(F\) down on the sphere.
To have \(F\) constitute a foliation by
geodesics, we will need each
leaf of \(F_1\) to sit inside a right translate
of a subgroup 
of \(G\) satisfying \(\omega^i=\gamma^i_1=0\).
Therefore we need \(t^i_1=0\),
and henceforth we will assume this;
i.e. on \(B_1\):
\begin{equation}\label{eqn:ColdCereal}
\gamma^i_1 = t^i_j \omega^j.
\end{equation}
The \emph{Reznikov invariant} is \(t=\pr{t^i_j}\) \cite{Reznikov:1985} p. 89.

Differentiating this last equation
gives an expression for the derivatives of the Reznikov invariant 
\[
\nabla t^i_j = d t^i_j 
- t^i_j \gamma^1_1 + \gamma^i_k t^k_j
- t^i_k \gamma^k_j - \delta^i_j \omega_1
\]
and we can calculate that
\begin{equation}\label{eqn:NablaA}
\nabla t^i_j = - t^i_k t^k_j \omega^1 + t^i_{jk} \omega^k
\end{equation}
where \(t^i_{jk} = t^i_{kj}\).

Note that this ``covariant derivative'' is determined not with a connection on the tangent bundle, but using the flat projective connection on the sphere.

\begin{lemma}\label{lm:Porridge}
Suppose that \(S^{2n+1} \to X^{2n}\) is a great circle fibration.
The tangent bundle to \(X^{2n}\) is soldered
by
\[
\gamma^i_j = \omega^i_j - \delta^i_j \omega^0_0.
\]
\end{lemma}
\begin{proof}
The proof is essentially
the same as that of
lemma~\vref{lm:HotCereal},
except that we use the representation
\[
\left ( \Line{e_0}^* \otimes V \right ) 
/
\left ( \Line{e_0}^* \otimes \Line{e_0,e_1} \right )
\]
where \(\Line{e_0,e_1}\) is the span of \(e_0,e_1\).
\end{proof}

\subsection{Structure equations of the Hopf fibration}
\label{subsec:Hopf}

To fix the Hopf fibration as well as the
flat projective structure,
transformations must take
complex lines to complex lines, since
the Hopf fibration on the sphere
is the quotient (by rescaling by positive numbers) 
of the fibration
of \(\C{n+1} \backslash 0\) into complex lines through
the origin. 
The circle fibers can be oriented
by using the natural orientation on
complex lines,
and we will take them to be thus
oriented.
\begin{lemma}
Let \(V\) be a
complex vector space
of dimension at least two.
Every invertible real linear map of 
\(V\)
which takes complex lines
to complex lines, preserving 
the natural orientation
of complex lines, is 
complex linear. 
\end{lemma}

The result is not true for \(V\) of one complex dimension.
\begin{proof}
Take \(e_1, e_2 \in V\) any two
vectors which are linearly independent
over the complex numbers. Then
\(ge_1,ge_2\) must still be linearly
independent over the complex numbers,
because \(e_1,e_2\) belong to distinct
complex lines, so \(g\) must take
these to distinct complex lines.
Take \(h\) any complex linear transformation
taking \(g e_1 \mapsto e_1, ge_2 \mapsto e_2\).
To show that \(g\) is complex linear on
the complex 2-plane spanned by \(e_1,e_2\)
it suffices to show that \(hg\) is.
So without loss of generality, we can
assume \(g e_1 = e_1\) and \(g e_2 = e_2\).
Consequently there must be real constants
\(a_j,b_j\) so that
\[
g \sqrt{-1} e_j = a_j e_j + b_j \sqrt{-1} e_j.
\]
The map \(g\) also must preserve the
complex line spanned by \(e_1+e_2\),
which forces \(a_1=a_2,b_1=b_2\).
Preserving the complex linear spanned
by \(e_1 + \sqrt{-1} e_2\) forces \(a_1=0,b_1=1\).
This makes \(g\) complex linear on the 
2-plane spanned by \(e_1,e_2\). Since 
\(e_1\) and \(e_2\) are arbitrary this
proves the lemma.
\end{proof}

So the group of symmetries
of the Hopf fibration as a geodesic
foliation of a flat projective structure
is the group
\[
\Gamma = \SL{2n+2,\R{}} \cap \GL{n+1,\C{}}.
\]
The subgroup preserving the point \(e_0\)
up to positive rescaling (i.e. fixing
the north pole of the sphere) is
the group \(\Gamma_0\) of complex
matrices of the form
\[
\begin{pmatrix}
g^0_0 & g^0_q \\
0 & g^p_q  
\end{pmatrix}
\]
where the indices \(p,q\) here run from \(1\) to \(n\),
with \(g^0_0 > 0\), the \(g^0_q\) and \(g^p_q\)
are complex numbers, and
\[
\left|g^0_0 \, \det\left(g^p_q\right)\right|^2=1
\]
(the real linear determinant must be 1).
The Hopf fibration is represented by the fiber bundle
\[
\begin{tikzcd}[ampersand replacement=\&]
\Gamma_0 \arrow{r} \& \Gamma \arrow{d} \\ 
                \& S^{2n+1}
\end{tikzcd}
\]
given by the obvious right action of \(\Gamma_0\)
on \(\Gamma\).
The group \(\Gamma\) is a subgroup of the
group \(G=\SL{2n+2,\R{}}\) that we encountered 
previously. 
Let 
\[
J_0 =
\begin{pmatrix}
0 & -1 &         &    &  \\
1 & 0  &         &    & \\
  &    & \ddots  &    & \\
  &    &         & 0  & -1 \\
  &    &         & 1  & 0
\end{pmatrix}
\]
be the usual complex structure on \(\R{2n+2}=\C{n+1}\),
and let
\[
K_0 =
\begin{pmatrix}
0 & 1 &         &    &  \\
1 & 0  &         &    & \\
  &    & \ddots  &    & \\
  &    &         & 0  & 1 \\
  &    &         & 1  & 0
\end{pmatrix}
\]
be the usual complex conjugation.
Taking the Maurer--Cartan
1-form \(\omega\) from \(G\), we can split
it into 
\[
\omega = \Omega^{1,0} + \Omega^{0,1}K_0
\]
by
\begin{align*}
\Omega^{1,0} &= \frac{1}{2} \left ( \omega - J_0 \omega J_0 \right ) \\
\Omega^{0,1} &= \frac{1}{2} \left ( \omega + J_0 \omega J_0 \right ) K_0.
\end{align*}
We can now write these in complex 
components, since each of \(\Omega^{1,0}\)
and \(\Omega^{0,1}\) is a matrix built
out of \(2 \times 2\) blocks like
\[
\begin{pmatrix}
a & -b \\
b & a
\end{pmatrix}
\]
which can be identified with the complex 
number \(a+b\sqrt{-1}\).
Write out
\begin{align*}
\Omega^{1,0} &= 
\begin{pmatrix}[\tallmatrix]
\Omega^0_0 & \Omega^0_{q} \\
\Omega^{p}_0 & \Omega^{p}_{q}
\end{pmatrix} \\
\Omega^{0,1} &= 
\begin{pmatrix}[\tallmatrix]
\Omega^0_{\bar{0}} & \Omega^0_{\bar{q}} \\
\Omega^{p}_{\bar{0}} & \Omega^{q}_{\bar{q}}
\end{pmatrix}.
\end{align*}
To be more explicit,
\begin{align*}
\Omega^P_Q &= \frac{1}{2} 
\left ( 
  \omega^{2P}_{2Q} +
  \omega^{2P+1}_{2Q+1} 
\right )
+
\frac{\sqrt{-1}}{2}
\left (
  \omega^{2P+1}_{2Q}
  -
  \omega^{2P}_{2Q+1}
\right ) \\
\Omega^P_{\bar{Q}} &= \frac{1}{2} 
\left ( 
  \omega^{2P+1}_{2Q} +
  \omega^{2P}_{2Q+1} 
\right )
+
\frac{\sqrt{-1}}{2}
\left (
  \omega^{2P+1}_{2Q+1}
  -
  \omega^{2P}_{2Q}
\right )
\end{align*}
for \(P,Q=0,\dots,n\).
Since the real trace of \(\omega\) vanishes, so does
\(\Omega^P_P + \Omega^{\bar{P}}_{\bar{P}}\).
We will write 
\(\Omega^{\bar{P}}_{\bar{Q}}\) for
the conjugate of \(\Omega^{P}_{Q}\)
and
\(\Omega^{\bar{P}}_Q\) for
the conjugate of \(\Omega^{P}_{\bar{Q}}\).
The structure equations of \(G\) can now
be written in this notation as
\begin{align*}
d \Omega^P_Q &= - \Omega^P_R \wedge \Omega^R_Q 
- \Omega^P_{\bar{R}} \wedge \Omega^{\bar{R}}_Q \\
d \Omega^P_{\bar{Q}} &= - \Omega^P_R \wedge \Omega^R_{\bar{Q}} 
- \Omega^P_{\bar{R}} \wedge \Omega^{\bar{R}}_{\bar{Q}}.
\end{align*}

The Hopf
fibration satisfies the equations
\(\Omega^P_{\bar{Q}}=0\).
These imply
\begin{align*}
\omega_1 &= - \omega^1 \\
\omega_{2P} &= \gamma^1_{2P+1} \\
\omega_{2P+1} &= - \gamma^1_{2P} \\
\gamma^1_1 &= 0 \\
\gamma^{2P}_1 &= - \omega^{2P+1} \\
\gamma^{2P+1}_1 &= \omega^{2P} \\
\gamma^{2P+1}_{2Q+1} &= \gamma^{2P}_{2Q} \\
\gamma^{2P}_{2Q+1} &= - \gamma^{2P+1}_{2Q} \\
\end{align*}
So for the Hopf fibration, the
invariant \(t\) is
\[
\left ( t^i_j \right ) =
\begin{pmatrix}
0 & -1 &         &    &  \\
1 & 0  &         &    & \\
  &    & \ddots  &    & \\
  &    &         & 0  & -1 \\
  &    &         & 1  & 0
\end{pmatrix}.
\]

\subsection{Structure equations of the standard round metric on the sphere}
It is helpful to compare the structure equations of a geodesic foliation
on the sphere to the equations we obtain with the standard metric in place.
We obtain the structure equations of \(\SO{2n+2}\) from those of \(\GL{2n+2,\R{}}\) by imposing the relations 
\begin{equation} \label{eqn:SO}
\begin{split}
\omega_{\mu} &= - \omega^{\mu} \\
\gamma^{\mu}_{\nu} &= - \gamma^{\nu}_{\mu}.
\end{split}
\end{equation}
In other words, \(\SO{2n+2}\) is the connected subgroup of \(\GL{2n+2,\R{}}\) of
largest dimension on which these equations are satisfied. 
\begin{lemma} 
The integral manifolds of the equations~\vref{eqn:SO} are precisely the left translates of \(\SO{2n+2}\) inside \(\GL{2n+2,\R{}}\).
\end{lemma}
\begin{proof} 
The Lie algebra of \(\SO{2n+2}\) satisfies these equations, so all of \(\SO{2n+2}\), as do all left translates of \(\SO{2n+2}\) by left invariance of the equations.
The left translates form a foliation of \(\GL{2n+2,\R{}}\) and by the Frobenius theorem, they are all of the integral manifolds. 
\end{proof}

Consider the equations
\begin{equation} \label{eqn:CircleSubgp}
\begin{split}
\omega_1 &= - \omega^1 \\
\omega_i &= 
\omega^i = 0 \text{ for } i > 1 \\
\gamma &= 0. 
\end{split}
\end{equation}
These are also the equations of a subgroup---in this case
the circle subgroup of \(\SO{2n+2}\) which turns the \(e_0,e_1\)
plane and fixes the perpendicular directions.
\begin{lemma} The integral curves of equation~\vref{eqn:CircleSubgp} are precisely the left translates of this circle subgroup of \(\SO{2n+2}\).
In particular these integral curves are compact and
\[
\int \omega^1 = 2 \pi
\]
when integrating over any of the integral curves.
\end{lemma}

\subsection{Following the invariants around the circles}\label{subsec:Follow}

Consider the behaviour of our invariant \(t\) on a circle subgroup left translate.
Imposing the relations satisfied by a circle subgroup translate given in equation~\vref{eqn:CircleSubgp}, we find that, if \(t= \left ( t^i_j \right )\)
is treated as a matrix then
\[
dt = - \left ( I + t^2 \right ) \omega^1.
\]
Such a circle subgroup translate lives inside our bundle \(B_1\) above each of our great circles, by construction of \(B_1\), since above each circle on the sphere, \(B_1\) contains a left translate of the subgroup \(G_1\) which
contains a circle subgroup. Writing \(\omega^1=d \theta\)
we have the ordinary differential
equation
\[
\frac{dt}{d \theta} = - \left ( I + t^2 \right ).
\]

Consider first how solutions of this ordinary differential equation behave if \(t\) is just a complex number.
The solutions are \(t=-\tan(\theta+c)\), except for the two singular solutions \(t=\pm \sqrt{-1}\).
So all solutions have period \(\pi\).
Indeed on the Riemann sphere, this equation is a rotation, and has the two exceptional points \(\pm\sqrt{-1}\) as rotational fixed points.
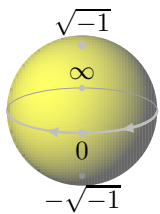
\begin{center}
\begin{tikzpicture}
  \shade[white,shading=sphere] (0,0) circle [radius=1cm];
  \draw[gray!50,thick] (1,0) arc (0:-180:1 and .3);
  \draw[gray!50,thick,-latex] (1,0) arc (0:-120:1 and .3);
  \draw[gray!50,thick,-latex] (1,0) arc (0:-60:1 and .3);
  \draw[gray!75,very thin] (1,0) arc (0:180:1 and .3);
  \fill[gray!50] (0,-.3) circle (1.5pt) node[below,black] {\(0\)}; 
  \fill[gray!50] (0,.3) circle (1.2pt) node[above,black] {\(\infty\)}; 
  \fill[gray!50] (0,.87) circle (1.5pt) node[above,black] {\(\sqrt{-1}\)}; 
  \fill[gray!85] (0,-.87) circle (1.2pt) node[below,black] {\(-\sqrt{-1}\)}; 
\end{tikzpicture}
\end{center}
To have our complex number \(t\) remain finite, it can never be real-valued, since
it would then go to infinity in time at most \(\pi/2\) in one direction or the other.

For a matrix \(t\) the solution is
\[
t(\theta) = \left( t(0) - \tan(\theta)I \right)
\left(I+\tan(\theta) t(0) \right)^{-1} 
\]
and this tells us that \(t\) has no real spectrum. 
Since \(t\) is a real matrix, the eigenvalues and eigenspaces of \(t\) come in conjugate pairs. Since the equation is invariant
under change of linear coordinates,
the eigenspaces will remain invariant
under the flow,
while the eigenvalues will change.
For any initial conditions \(t(0)\)
with no real eigenvalues, we find
that \(t\) will remain defined
as a function of \(\theta\) for 
any positive or negative \(\theta\).
For generic initial conditions,
for example for \(t\) diagonalizable
at \(\theta=0\),
we find that \(t\) is \(\pi\) periodic.
Therefore \(T\) is also \(\pi\) periodic
for any initial conditions with
no real eigenvalues.

Notice that we have not so far invoked the hypothesis that the sphere has
odd dimension. In fact, this is an immediate consequence of the Reznikov invariant having no real eigenvalues, since \(t\) is a square matrix whose size is one less than the dimension of the sphere. 
If \(t\) had odd size, then its characteristic polynomial would have odd degree, and so would have a real root---hence a real eigenvector.

\subsection{Linear transformations without
real eigenvalues}\label{subsec:LinTrans}
Take any linear transformation \(T \colon V \to V\)
of a real vector space \(V\),
with no real eigenvalues. 
We will also write \(V\) as \(V_{\R{}}\)
to emphasis that we are studying
its real points, and write
\(V_{\C{}}\) for \(V \otimes_{\R{}} \C{}\).
We naturally have \(V_{\C{}}\) split
into generalized eigenspaces of \(T\),
say 
\[
E_{\lambda} T =
\left \{ v \in V_{\C{}} \, | \, 
(T-\lambda I)^k v = 0,
\text{for some } k > 0 \right \}.
\]
Since \(T\) is real, 
\(
\overline{E_{\lambda} T}
=
E_{\bar{\lambda}} T.
\)
Pick the eigenspaces \(E_{\lambda}T\) of \(T\) where the eigenvalues \(\lambda\) have positive imaginary parts, and define a subspace
\[
V^{1,0}_T \defeq \bigoplus_{\Im{\lambda} > 0} E_{\lambda} T \subset V_{\C{}}.
\]
Let \(V^{0,1}\) be the conjugate of \(V^{1,0}\),
i.e. the sum of the eigenspaces
whose eigenvalues have negative 
imaginary part.
Note that \(V^{1,0}\) and \(V^{0,1}\) are complex subspaces of \(V_{\C{}}\) and are complementary:
\[
V^{1,0} \cap V^{0,1} = 0,
V_{\C{}} = V^{1,0} \oplus V^{0,1}.
\]
Define a linear transformation \(J \colon V_{\C{}} \to V_{\C{}}\) by letting \(J\) act via \(\sqrt{-1}\) on \(V^{1,0}\) and by \(-\sqrt{-1}\) 
on \(V^{0,1}\). 
This \(J\) is real, i.e. acts as a real linear transformation on
\(V_{\R{}}\), since \(J\) is conjugation invariant, with \(J^2=-I\), since this equation holds on \(V^{1,0}\) and on \(V^{0,1}\) by construction.
Write \(J\) as \(J_T\).

\begin{lemma} Let \(H_V\)
be the set of linear transformations
of a finite dimensional 
real vector space \(V\) which have
no real eigenvalues (think of it as
a ``generalized upper half plane''), 
and \(\Cstrucs{V}\)
be the subset of real linear transformations
\(J\) which satisfy \(J^2 = -I\) (the
complex structures). Then
the map
\[
T \in H_V \mapsto J_T \in \Cstrucs{V}
\]
is a smooth fiber bundle, with the
inclusion \(\Cstrucs{V} \subset H_V\)
as a section. This fiber
bundle is \(\GL{V}\)-equivariant:
\(
J_{gTg^{-1}} = gJ_T g^{-1}
\)
for any \(g \in \GL{V}\).
Its fiber above any point \(J\) consists
precisely of the \(J\)-complex linear
maps \(T \colon V \to V\), i.e. \(TJ=JT\),
all of whose eigenvalues on \((J=\sqrt{-1})\subset V_{\mathbb{C}}\)
have positive imaginary part.
If \(c \ne 0\) is a real number, then
\(
J_{T + cI} = J_T
\)
and
\(
J_{cT} = \operatorname{sign}(c) J_T.
\)
\end{lemma}
\begin{proof}
Notice that if \(V\) has odd dimension,
then (looking at the characteristic
polynomial, which is of odd degree)
we find \(H_V\) is empty.
Similarly, taking determinant, we find
that \(\Cstrucs{V}\) is empty.
So we can assume \(V\) has even dimension.

The equivariance of \(T \mapsto J_T\)
under \(\GL{V}\) is elementary.

First, we construct the space \(Z \subset  H_V \times \Cstrucs{V}\)
which consists of pairs \((T,J)\) satisfying
\[
TJ=JT, J^2 = -I
\]
and with \(T\) having no real eigenvalues.
We have maps
\[
\begin{tikzcd}[ampersand replacement=\&]
\& Z \arrow{dl} \arrow{dr} \&  \\
H_V \& \& \Cstrucs{V}
\end{tikzcd}
\]
given by taking a pair \((T,J)\) and
either forgetting \(J\) or forgetting \(T\).
Now given any \(T \in H_V\),
we know how to construct a \(J=J_T\)
commuting with it so that \(J^2=-I\),
i.e. the map
\(
T \mapsto J_T
\)
has graph lying in \(Z\).

On the other hand, if we take
any \(T\) with no real eigenvalues,
i.e. \(T \in H_V\),
then we can could also have constructed
a \(J\) commuting with it in other ways;
in fact if we just put any collection
of eigenspaces of \(T\) together
and call their sum \(V^{1,0}\) and call the sum of their 
conjugates \(V^{0,1}\), so that every eigenspace
is in one or the other, then we can
define \(J\) to be \(\sqrt{-1}\) on \(V^{1,0}\)
and \(-\sqrt{-1}\) on \(V^{0,1}\). 
Because \(T\) has no real eigenvalues, this procedure
unambiguously determines \(J\), up to
the choices of which eigenspaces of \(T\) go
to \(V^{1,0}\) and which go to
\(V^{0,1}\).

Conversely, if we pick a pair \((T,J)\)
in \(Z\), the fact that \(T\) and \(J\) 
commute ensures that each one leaves
the eigenspaces of the other invariant.
Looking at the minimal polynomial of \(J\),
\(x^2+1\), we see that \(J\) has exactly
two eigenspaces, and we call them
\(V^{1,0}\) (the \(\sqrt{-1}\)-eigenspace)
and \(V^{0,1}\) (the \(-\sqrt{-1}\)-eigenspace). The map \(T\) takes
\(V^{1,0}\) to itself, and takes \(V^{0,1}\)
to itself, because \(T\) commutes with \(J\),
and \(T\) is complex linear on each. 
Decomposing \(V^{1,0}\) into complex eigenspaces of \(T\), and 
decomposing \(V^{0,1}\) into the conjugates of those eigenspaces,
we find \(V_{\C{}} = V^{1,0} \oplus V^{0,1}\) decomposed into conjugate 
eigenspaces of \(T\).
Because \(T\)  is real, its eigenspaces are
conjugate, with half of them in \(V^{1,0}\)
and half in \(V^{0,1}\). Therefore
\(J\) comes about from \(T\) by the construction
outlined in the last paragraph.

Thus the map \(Z \to H_V\)
has finitely many points in each 
stalk---each point corresponds to a 
choice of which eigenspaces of \(T\) 
go into \(V^{1,0}\).
We want to show that the points which
make up the graph of \(T \mapsto J_T\)
form a smooth subvariety of \(Z\). We have only
to show that they are transverse to the
fibers of the map 
\[
H_V \times \Cstrucs{V} \to H_V
\]
given by \((T,J) \mapsto T\).

The algebraic equations cutting out \(Z\) are \(TJ=JT\) and \(J^2=-I\).
Differentiating these in motions up the fiber, we find that the equations of a vertical tangent vector are \(T \dot{J} = \dot{J} T\) and \(\dot{J} J + J \dot{J} = 0\).
The first equation tells us that \(\dot{J}\) preserves the eigenspaces of \(T\), while the second tells us that \(\dot{J}\) swaps the eigenspaces of \(J\). 
But the eigenspaces of \(T\) are entirely contained inside those of \(J\).
Therefore \(\dot{J}=0\), and there are no vertical tangent vectors. 
This shows that \(Z\) is smooth and transverse to the fibers of 
\(H_V \times \Cstrucs{V} \to \Cstrucs{V}.\) 
So the map \(T \mapsto J_T\) is smooth, being just a single branch of \(Z\), and  \(J_T=T\) for any \(T\) with \(T^2=-I\). 
Hence this map \(T \mapsto J_T\) is a smooth surjection \(H_V \to \Cstrucs{V}\).

To find the rank of the map \(T \mapsto J_T\) differentiate the equations \(J^2=-I\) and \(JT=TJ\) to find 
\[
\dot{J}J + J \dot{J} = 0 \text{ and } 
\dot{J} T + J \dot{T} = \dot{T} J + T \dot{J}.
\]
The kernel of the derivative of \(T \mapsto J_T\) is the set of \(\dot{T}\) satisfying \(J \dot{T} = \dot{T}J\), i.e. the \(J\) complex linear maps. 
Since the real general linear group acts transitively on complex structures, we find that the space of such \(\dot{T}\) has dimension independent of \(J\). 
Indeed if \(\dim_{\R{}} V = 2n\), then the derivative of the map \(T \mapsto J_T\) 
has kernel of dimension \(2n^2\), and so fibers of dimension also
\(2n^2\). 
The dimension of the base is also \(2n^2\), so the map \(T \mapsto J_T\) is a smooth submersion.

Next we want to show that this map is locally trivial, so that it will be a fiber bundle.
First, we note once again that \(\Cstrucs{V}\) is a homogeneous space under the
action of \(\GL{V}\). 
So if we pick a particular \(J_0 \in \Cstrucs{V}\) then we have the fiber bundle
\[
\begin{tikzcd}[ampersand replacement=\&]
\GL{V,J_0} \arrow{r} \& \GL{V} \arrow{d} \\
                  \& \Cstrucs{V}.
\end{tikzcd}
\]
Every homogeneous space is the base spaces of a fiber bundle, with total space given by the transitively acting group.
Thus this bundle is locally trivial, i.e. to every sufficiently small open set
\(U \subset \Cstrucs{V}\) we can associate a map \(U \to \GL{V}\), say \(g(J)\), so that if \(g=g(J)\) then \(gJg^{-1} = J_0\).

We map the fiber \(H_{V,J_0}\) above \(J_0\) to the fiber \(H_{V,J}\) above \(J\) by
\[
T_0 \mapsto T=T(J,T_0) = gT_0g^{-1}.
\]
By \(\GL{V}\)-equivariance, this is a local trivialization of \(H_V \to \Cstrucs{V}\), so that this map is a fiber bundle.
\end{proof}

The reader who wishes to understand this lemma might work out the whole story in matrices for \(V=\R{2}\).
This complex structure is \(J\) due to Reznikov \cite{Reznikov:1985} p. 89.

\subsection{Back to the great circle fibration}\label{subsec:Back}

Inside \(\SO{2n+2}\) we do not have enough room to put the Reznikov invariant \(t\) into a normal form. 
If we return to \(G=\SL{2n+2,\R{}}\) we have enough room to arrange by moving
up and down the fibers of \(B_1\) that at least the associated linear map \(J_t\) is normalized:
\[
J_t = J_0 =
\begin{pmatrix}
0 & -1 &         &    &  \\
1 & 0  &         &    & \\
  &    & \ddots  &    & \\
  &    &         & 0  & -1 \\
  &    &         & 1  & 0
\end{pmatrix}.
\]
There is a subbundle \(B_2 \subset B_1\) on which these equations hold.
This bundle \(B_2\) is a principal \(G_2\) subbundle, where \(G_2\) is the group of matrices of the form
\[
\begin{pmatrix}
g^0_0 & g^0_1 & g^0_j \\
0 & g^1_1 & g^1_j \\
0 & 0 & g^i_j 
\end{pmatrix}
\]
with \(g^0_0,g^1_1 > 0\) and \(g^0_0 g^1_1 \det \left ( g^i_j \right ) = 1\)
and where \(\left ( g^i_j \right )\) commutes with \(J\).

This matrix \(t^i_j\) is now \(J_0\) complex linear,
so that we can say that the 1-forms
\[
\omega^{2p}_1 + \sqrt{-1} \omega^{2p+1}_1
\]
are complex multiples of the 1-forms
\[
\omega^{2p}_0 + \sqrt{-1} \omega^{2p+1}_0.
\]
Written in terms of the complex
\(\Omega\) notation, 
\[
\omega^{2p}_1 + \sqrt{-1} \omega^{2p+1}_1
= 
\Omega^p_{\bar{0}} + \sqrt{-1} \Omega^p_0
\]
and
\[
\omega^{2p}_0 + \sqrt{-1} \omega^{2p+1}_0
= 
\Omega^p_{0} + \sqrt{-1} \Omega^p_{\bar{0}}
\]
so that
the equation \(\gamma^i_1=t^i_j \omega^j\)
can be written as
\begin{equation} \label{eqn:PPPPPPPP}
\Omega^p_{\bar{0}} + \sqrt{-1} \Omega^p_0
=
t^p_q \left (
        \Omega^q_0 + 
        \sqrt{-1} \Omega^q_{\bar{0}}
        \right )
\end{equation}
with \(p,q=1,\dots,n\), and \(t^p_q\)
a complex matrix whose eigenvalues
lie in the upper half plane.

Lets write linear fractional transformations
using the notation
\[
\begin{bmatrix}
a & b \\
c & d
\end{bmatrix}
z = \frac{az+b}{cz+d}.
\]
Define the Cayley map
\[
C(z) = 
\begin{bmatrix}
\sqrt{-1} & 1 \\
1 & \sqrt{-1}
\end{bmatrix}z.
\]
Solving equation~\vref{eqn:PPPPPPPP}
for \(\Omega^p_{\bar{0}}\) in terms
of \(\Omega^q_0\), we find
\(
\Omega^p_{\bar{0}} = s^p_q \Omega^q_0
\)
expressed in terms of the Cayley mapped
matrix \(s=C(t)\).
Henceforth we employ \(s\) instead of \(t\) and call \(s\) the \emph{Reznikov invariant}.
The condition that \(t\) have all of its eigenvalues in the upper half plane is equivalent to its image \(s=C(t)\) under the Cayley map
having all of its eigenvalues inside the unit disk. 
In particular, for the Hopf fibration, on the subbundle \(\Gamma \subset B_2\) we find \(s=0\).

Consider the ordinary differential equation 
\[
\frac{dt}{d \theta} = - \left( 1 + t^2 \right).
\]
If we set \(s=C(t)\), then
\[
\frac{ds}{d \theta} = - 2 \sqrt{-1} s
\]
so that 
\(
s(\theta) = s(0) e^{-2 \sqrt{-1} \theta}
\)
evolves by rotation.
We will see that the Reznikov invariant is a section of a complex vector bundle over the base manifold of the circle fibration (i.e. the manifold parameterizing the circles). 
For the Hopf fibration, this base manifold is \(\CP{n}\).

The restriction of structure group imposed by the equation \(J_t=J_0\)
requires \(\Omega^p_{\bar{q}}\) to be semibasic as well, for \(p,q=1,\dots,n\).
Recall that the 1-forms \(\omega^1,\omega^i\) form a basis for the semibasic 1-forms.
The 1-forms \(\Omega^p_0\) are semibasic, but (even together with their conjugates) they do not span the semibasic 1-forms---only
the \(\omega^i\) (\(i=2,\dots,2n+1\)) are multiples of them, while \(\omega^1\) is not. 

Taking the exterior derivative of both sides of the equation
\(\Omega^p_{\bar{0}} = s^p_q \Omega^q_0\),
\begin{align*}
0=&\left (
ds^p_q 
+ \Omega^p_r s^r_q - s^p_r \Omega^r_q
- \Omega^0_{\bar{0}} \delta^p_q
+ \left ( \Omega^0_0 - \Omega^{\bar{0}}_{\bar{0}} \right) s^p_q
+ s^p_r s^r_q \Omega^{\bar{0}}_{0}
\right )
\wedge \Omega^q_0
\\&+
\left (
\Omega^p_{\bar{q}} - s^p_r \Omega^r_{\bar{t}} s^{\bar{t}}_{\bar{q}}
\right )
\wedge \Omega^{\bar{q}}_{\bar{0}}.
\end{align*}
By Cartan's lemma, this implies that
there are constants 
\[
s^p_{qr}=s^p_{rq}, \ 
s^p_{q \bar{r}} = s^p_{\bar{r}q}, \
s^p_{\bar{q}\bar{r}} = s^p_{\bar{r}\bar{q}}
\]
so that
\begin{align}
ds^p_q 
+ \Omega^p_r s^r_q - s^p_r \Omega^r_q
- \Omega^0_{\bar{0}} \delta^p_q
+ \left ( \Omega^0_0 - \Omega^{\bar{0}}_{\bar{0}} \right) s^p_q
+ s^p_r s^r_q \Omega^{\bar{0}}_{0}
&=
s^p_{qr} \Omega^r_0 + s^p_{q\bar{r}} \Omega^{\bar{r}}_{\bar{0}} 
\label{eqn:ds} \\
\Omega^p_{\bar{q}} - s^p_r \Omega^r_{\bar{t}} s^{\bar{t}}_{\bar{q}} 
&=
s^p_{\bar{q}r} \Omega^r_0 + s^p_{\bar{q}\bar{r}} \Omega^{\bar{r}}_{\bar{0}}.
\label{eqn:Popping}
\end{align}

We will need to see that it is possible to solve for \(\Omega^p_{\bar{q}}.\)
Consider the operation on matrices \(M \mapsto M - sM\bar{s}\).
We need to show that it is invertible, as long as the matrix \(s\) has all of its eigenvalues in the unit disk. 
The kernel of this operation consists in matrices
\(M\) satisfying \(M=sM\bar{s}\). 
This implies that \(sM\bar{s} = s^2M\bar{s}^2\) etc. so that \(M=s^k M \bar{s}^k\)
for all positive integers \(k\). 
But sufficiently high powers of \(s\) are strictly contracting, since all eigenvalues of \(s\) are in the unit disk, so \(M=0\).
The equation~\vref{eqn:Popping} can be solved for \(\Omega^p_{\bar{q}}\)
as a complex linear combination of
the 1-forms \(\Omega^r_0\) and \(\Omega^{\bar{r}}_{\bar{0}}\).

The Reznikov invariant \(s\) transforms under
motions through the fibers of the bundle \(B_2\)
(via the action of the structure group \(G_2\))
in a very complicated action, which consists
of conjugations via the \(\Omega^p_q\) 
``variables,'' and, as we will see, linear fractional
transformations via the \(\Omega^0_0-\Omega^{\bar{0}}_{\bar{0}}\)
and \(\Omega^0_{\bar{0}},\Omega^{\bar{0}}_{0}\)
``variables.''

\subsection{Linear fractional transformations acting on matrices}\label{subsec:LFT}
Recall that the linear fractional transformations
\[
\begin{bmatrix}
a & b \\
c & d 
\end{bmatrix}
z
=
\frac{az+b}{cz+d}
\]
constitute an action of \(\PSL{2,\C{}}=\SL{2,\C{}}/\pm{1}\)
on the Riemann sphere, 
\[
\begin{bmatrix}
a & b \\
c & d 
\end{bmatrix}
\in \PSL{2,\C{}},
\ z \in \C{} \cup \infty.
\]
The infinitesimal action is
\[
\begin{pmatrix}
a & b \\
c & -a
\end{pmatrix}
z
=
\left ( b + 2az - cz^2 \right ) \pd{}{z}
\]
for
\[
\begin{pmatrix}
a & b \\
c & -a
\end{pmatrix} \in \sl{2,\C{}}.
\]

The subgroup of \(\PSL{2,\C{}}\) which
preserves the unit disk \(D\), call it \(\Aut{D}\),
is the quotient modulo \(\pm{1}\) of 
the group of matrices 
of the form
\[
\begin{pmatrix}
a & b \\
\bar{b} & \bar{a}
\end{pmatrix}
\in \SL{2,\C{}}
\]
subject to \(|a|^2 - |b|^2=1\).
The Lie algebra of \(\Aut{D}\) consists of the matrices of the form
\[
\begin{pmatrix}[\tallmatrix]
\sqrt{-1}a & b \\
\bar{b} & - \sqrt{-1}a
\end{pmatrix}
\]
where \(a \in \R{}, b \in \C{}\).
The Cayley map identifies the upper half
plane with the unit disk, identifying \(\Aut{D}\)
with \(\SL{2,\R{}}\). There is a (unique
up to conjugation)
connected 2-dimensional subgroup of \(\SL{2,\R{}}\);
it is nilpotent and consists of the
matrices of the form
\[
\begin{pmatrix}
a & b \\
0 & 1/a
\end{pmatrix}
\in \SL{2,\R{}}
\]
with \(a,b \in \R{}\) and \(a>0.\)
The quotient of this group by \(\pm{1}\)
is the unique connected 2-dimensional
subgroup of \(\PSL{2,\R{}}\).
Such a matrix acts on the element \(\sqrt{-1}\)
in the upper half plane by
\[
\begin{bmatrix}
a & b \\
0 & 1/a
\end{bmatrix}
\sqrt{-1}
=
a^2 \sqrt{-1} + ab,
\]
a transitive action.

Under the Cayley map, this 2-dimensional
subgroup gets mapped to a subgroup
\(N\) of \(\Aut{D}\). One can readily 
conjugate
with the Cayley mapping to calculate that 
the elements of \(N\) are 
precisely the elements of \(\PSL{2,\C{}}\) 
of the form
\begin{equation}\label{eqn:POP}
\begin{bmatrix}[\tallmatrix]
a+\frac{1}{a} + \sqrt{-1}b &
b+\sqrt{-1} \left ( a - \frac{1}{a} \right ) \\
b-\sqrt{-1} \left ( a - \frac{1}{a} \right ) &
a+\frac{1}{a} -\sqrt{-1}b
\end{bmatrix}
\end{equation}
where \(a,b \in \R{}\) and \(a>0.\)
This group \(N\) acts transitively on the
unit disk \(D\), because its conjugate
under the Cayley map acts transitively
on the upper half plane.
In particular, the Lie algebra of
\(N\), call it \(\mathfrak{n}\), consists
of matrices of the form
\[
\begin{pmatrix}[\tallmatrix]
\frac{\sqrt{-1}}{2}\left(Q+\bar{Q}\right) & Q \\
\bar{Q} & -\frac{\sqrt{-1}}{2}\left(Q+\bar{Q}\right)
\end{pmatrix}
\]
where \(Q\) can be any complex number.
The infinitesimal action on the unit disk of 
such a Lie algebra
element, say \(M\), is
\[
Mz
=
\left(
  Q
  -\sqrt{-1}\left(Q+\bar{Q}\right)z
  -\bar{Q}z^2
\right)\pd{}{z}.
\]

A matrix can be plugged in to a linear fractional transformation as long as its spectrum lies in the domain where the linear fractional transformation is finite. 
In particular, for \(s\) a complex matrix whose spectrum lies inside the unit disk, all elements of \(\Aut{D}\) can act on \(s\).
A 1-parameter family of motions of a matrix \(s\) by elements of \(N\) is the same as an ordinary differential equation like
\[
\frac{ds}{dt} =
 Q
 -\sqrt{-1}\left(Q+\bar{Q}\right)s
 -\bar{Q}s^2
\]
where the \(Q(t)\) is any smooth complex-valued function of a real variable \(t\).

We will henceforth orient the group \(N\) using the orientation of the Lie algebra
given by the usual orientation of the \(Q\) complex plane. 
Write the Maurer--Cartan 1-form on \(N\) as
\[
\begin{pmatrix}[\tallmatrix]
-\frac{\sqrt{-1}}{2}\left(\psi+\bar{\psi}\right) & \psi \\
\bar{\psi} & \frac{\sqrt{-1}}{2}\left(\psi+\bar{\psi}\right)
\end{pmatrix}
\]

Reconsidering our equation~\vref{eqn:ds},
we now see that this is precisely
the sort of motion that \(s\) undergoes
when we move up the fibers of the
bundle \(B_2\), at least as long as
the motion is only in the \(\omega_1,\omega^1_1,\omega^0_0\)
directions. Indeed equation~\vref{eqn:ds}
tells us that
\[
ds^p_q = 
- \Omega^0_{\bar{0}} \delta^p_q
+ \left ( \Omega^0_0 - \Omega^{\bar{0}}_{\bar{0}} \right) s^p_q
+  \Omega^{\bar{0}}_{0} s^p_r s^r_q
\pmod{\Omega^r_0,\Omega^r_s}.
\]
Converting this into \(\omega\) notation:
\begin{align*}
\Omega^0_{\bar{0}} &= \frac{1}{2} \omega_1 
+ \frac{\sqrt{-1}}{2}\gamma^1_1 \\
\Omega^0_0 - \Omega^{\bar{0}}_{\bar{0}}
&= -\sqrt{-1}\omega_1
\end{align*}
modulo semibasic terms. 
This implies that
\[
ds^p_q = - \left(\frac{1}{2}\omega_1
+\frac{\sqrt{-1}}{2}\gamma^1_1\right)\delta^p_q
-\sqrt{-1}\omega_1 s^p_q
+ s^p_r s^r_q \left(\frac{1}{2}\omega_1-\frac{\sqrt{-1}}{2}\gamma^1_1\right)
\pmod{\omega^1,\Omega^r_0,\Omega^r_s}.
\]
Therefore we can set
\(
\psi=-\Omega^0_{\bar{0}}
\)
and find that any motion through the leaves of the
foliation \(\omega^1=\Omega^r_0=\Omega^r_s=0\) effects
an action of \(N\) on \(s^p_q\). 
Any element of \(N\) arises in this manner, since \(N\)
is a connected Lie group, as is obvious from equation~\vref{eqn:POP}.

\begin{lemma}
Suppose that \(s\) is a square matrix with complex entries, and that the
spectrum of \(s\) is contained in the unit disk. 
Then there is a unique linear fractional transformation belonging to the group \(N\), say
\[
g=
\begin{bmatrix}[\tallmatrix]
a+\frac{1}{a} + \sqrt{-1}b &
b+\sqrt{-1} \left ( a - \frac{1}{a} \right ) \\
b-\sqrt{-1} \left ( a - \frac{1}{a} \right ) &
a+\frac{1}{a} -\sqrt{-1}b
\end{bmatrix}
\]
so that \(g(s)\) has trace zero. 
This transformation \(g\) depends analytically on \(s\).
\end{lemma}
\begin{proof}
Assume \(s\) is an \(n \times n\) matrix.
First, suppose that \(s\) has trace zero, with 
eigenvalues \(\lambda_1,\lambda_2,\dots,\lambda_n.\)
Then infinitesimal motions under the group \(N\)
affect the trace by
\[
\frac{1}{n} ds^p_p = \psi 
- \bar{\psi} \frac{1}{n} s^p_q s^q_p.
\]
The first term has larger complex coefficient
than the second,
since having spectrum in the unit disk
forces
\[
\frac{1}{n} s^p_q s^q_p
=
\frac{1}{n} \lambda^2_p
=
\text{average squared eigenvalue}
\] 
inside the unit disk. Therefore the differential
\(ds^p_p\) as a linear map from the tangent space
of \(N\) to \(\C{}\)
is orientation preserving, and has full rank,
at every zero of \(s^p_p.\)
In particular, all zeros of \(s^p_p\) are nondegenerate
and positively oriented, and the set
of elements \(g \in N\) at which \(\tr(g(s))=0\)
is discrete, for any matrix \(s\)
with spectrum in the unit disk.

Since \(N\) acts
transitively and via isometries of
the hyperbolic metric on the unit disk,
we can arrange that the center
of mass of the spectrum, in the hyperbolic
metric, sits wherever we like.
The problem of arranging a vanishing
trace is that we have to get the 
center of mass in the Euclidean metric
to vanish, and it is not obvious that this
is possible.

We will analyze the behaviour of
elements of \(N\) ``near infinity,'' i.e.
far away from the identity element.
The points of the spectrum
retain their hyperbolic distances
under maps from \(N\), because the
elements of \(N\) are hyperbolic
isometries, but the points of 
the spectrum become
very close in Euclidean norm if
any one of them approaches 
the boundary of the
unit disk, since  near the boundary of the disk,
the hyperbolic balls
of fixed radius are contained
in Euclidean balls of very small
radius.
Hence we can easily control the
average of the eigenvalues of \(s\), 
to get the average of the
eigenvalues close to
any number \(e^{i \theta}\) on the
boundary of the unit disk, using 
linear fractional transformations
from \(N\).

Using the coordinates \(a,b\) for the 
group \(N\), we can see that elements
of \(N\) with \(a\) close to zero (but
positive) are close to the constant
map
\(z \mapsto -\sqrt{-1}\)
given by the matrix
\[
\begin{bmatrix}[\tallmatrix]
1 & -\sqrt{-1} \\
\sqrt{-1} & 1
\end{bmatrix}.
\]
This matrix is not an invertible
matrix, and obviously the map it
generates is not either.
Nonetheless, applying elements of \(N\)
near \(a=0\) to our matrix \(s,\) we can arrange
that the entire spectrum of \(s\) 
lies close to \(-\sqrt{-1}\).
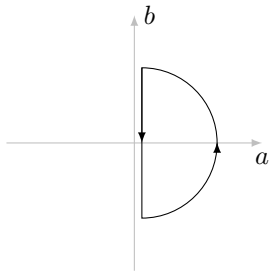
\begin{center}
\begin{tikzpicture}
\draw[gray!50,-latex] (-1.7,0) -- (1.7,0) node[black,below] {\(a\)};
\draw[gray!50,-latex] (0,-1.7) -- (0,1.7) node[black,right] {\(b\)};
\draw (0.1,1) arc (90:-90:1) -- cycle;
\draw[-latex] (1.1,-0.01) -- (1.1,0.01);
%\draw[-latex] (0.1,-1) arc (-90:0:1);
\draw[-latex] (0.1,1) -- (0.1,0);
\end{tikzpicture}
\end{center}
Consider the half circle
\[
a = r \cos \theta, \ b = r \sin \theta
\]
for values of \(\theta\) from a little
above \(-\pi/2\) to a little less than
\(\pi/2\). We find that for \(r\) a large
positive number, the origin of
the unit disk is taken by this linear
fractional transformation to
\[
\frac{r^2 \sin \theta \cos \theta - 
\sqrt{-1}\left(1 - r^2 \cos^2 \theta\right)}
{
\left(1+r^2 \cos \theta\right)
-\sqrt{-1}r^2 \sin \theta \cos \theta
}
\]
which, for large \(r\), is close to \(\sqrt{-1}\)
for all values of \(\theta.\)
Now try the path
\(
a = 1/\sqrt{r}, b = \sqrt{r}t
\), 
\(-\sqrt{r} \le t \le \sqrt{r}\).
The origin is taken to the point
\[
\frac{t + \left(\frac{1}{r} - 1\right)\sqrt{-1}}
{1+\frac{1}{r}-t\sqrt{-1}}
\]
which, for large \(r\), is very close to the point
\[
\frac{t - \sqrt{-1}}{1-t\sqrt{-1}}
\]
which traverses the unit circle
as \(t\) runs from \(-\infty\) to \(\infty.\)

Consider any 1-parameter family \(g_{\theta}\) of
linear fractional transformations from \(N\) which stays close to the contour
\[
a = r \cos \theta, b = r \sin \theta, \quad
-\pi/2+\varepsilon < \theta < \pi/2-\varepsilon
\]
(some small \(\varepsilon\)) and then approaches a vertical line in the \(a,b\) plane, close to, but just to the right of  \(a=0\), say at \(a=1/\sqrt{r},\) completing a loop. 
Applying \(g_{\theta}\) to a matrix \(s\), we find that (if the radius \(r\) of the half circle in our contour is set large enough) the
average of the eigenvalues of \(g_{\theta}(s)\) travels quite
near the boundary of the unit disk, through a 
single rotation. 
In fact, since the hyperbolic distance between eigenvalues is preserved by these linear fractional transformations, the Euclidean distance shrinks as
we approach the boundary of the unit disk, and all of the eigenvalues are very close to the average.
This ensures that there must be a zero of the average of the
eigenvalues of \(g(s)\) for some value of \(a,b\) inside the semicircle
of radius \(r\), since otherwise the winding number of the average
eigenvalue would be unchanged as we shrunk the contour down to a point.

The winding number
of the average eigenvalue
around the loop
is \(1,\) and the zeros of \(\tr(g(s))\)
are positively oriented, so there can only be one
such zero, and it must
depend analytically on \(s\) by
the implicit function theorem.
\end{proof}

For example, consider the transformations
\(g\) which take \(s=0\) to a traceless matrix.
They must look like
\begin{align*}
g(0) &= 
\begin{bmatrix}[\tallmatrix]
a+\frac{1}{a} + \sqrt{-1}b &
-b+\sqrt{-1} \left ( a - \frac{1}{a} \right ) \\
-b-\sqrt{-1} \left ( a - \frac{1}{a} \right ) &
a+\frac{1}{a} -\sqrt{-1}b
\end{bmatrix}0 \\
&=\left(-b+\sqrt{-1}\left(a-\frac{1}{a}\right)\right)
\left(a+\frac{1}{a}-\sqrt{-1}b\right)^{-1}I
\end{align*}
which can not have vanishing trace unless \(a=1\) and \(b=0\),
in which case \(g\) represents the identity
transformation.

\subsection{Reducing the structure group}\label{subsec:Reduce}

We have seen that we can arrange \(s^p_p=0\), and
that this occurs on a subbundle of \(B_2\),
call it \(B_3\).
Plugging into equation~\vref{eqn:ds}
and taking trace, we find that there
are some functions 
\(s^0_q,s^0_{\bar{q}},s^p_{\bar{q}r}\)
and
\(s^p_{\bar{q}\bar{r}}\)
satisfying
\[
s^p_{\bar{q}\bar{r}}-s^p_u s^u_{\bar{t}\bar{r}} s^{\bar{t}}_{\bar{q}}
=
s^p_{\bar{r}\bar{q}}-s^p_u s^u_{\bar{t}\bar{q}} s^{\bar{t}}_{\bar{r}}
\]
so that
\begin{equation}\label{eqn:Putty}
\begin{pmatrix}
\Omega^0_{\bar{0}} \\
\Omega^p_{\bar{0}} \\
\Omega^p_{\bar{q}} 
\end{pmatrix}
= 
\begin{pmatrix}
s^0_r & s^0_{\bar{r}} \\
s^p_r & 0 \\
s^p_{\bar{q} r} & s^p_{\bar{q} \bar{r}}
\end{pmatrix}
\begin{pmatrix}
\Omega^r_{0} \\
\Omega^{\bar{r}}_{\bar{0}} 
\end{pmatrix}
\end{equation}
with \(s^p_p=0.\) 

\begin{corollary}
Inside
our bundle \(B_2\) there is a smaller
bundle \(B_3\) of points at which
the trace \(s^p_p\) vanishes.
The bundle \(B_3\) is a principal
\(G_3\) bundle where \(G_3\) is
the group of matrices of the
form
\[
\begin{pmatrix}
g^0_0 & 0 & g^0_j \\
0 & g^0_0 & g^1_j \\
0 & 0 & g^i_j
\end{pmatrix}
\]
where \(g^0_0 > 0,\) 
\(\left(g^0_0\right)^2 \det\left(g^i_j\right) = 1\) 
and \(\left(g^i_j\right)\)
commutes with \(J\).
\end{corollary}
\begin{proof} This is the only
closed Lie subgroup of \(G_2\) with the required Lie
algebra.
\end{proof}

Differentiating the equation for \(\Omega^0_{\bar{0}}\)
in equations~\vref{eqn:Putty},
we find that modulo semibasic terms:
\[
ds^0_{\bar{q}} = - \Omega^0_{\bar{q}} 
+ s^0_{\bar{p}} \Omega^{\bar{p}}_{\bar{q}}
+ s^0_{\bar{q}} \Omega^{\bar{0}}_{\bar{0}}
\pmod{\Omega^r_{0},\Omega^{\bar{r}}_{\bar{0}}}
\]
so that as we move up the fibers of
\(B_3\), we can arrange that \(s^0_{\bar{q}}\)
vanish on a subbundle; call it \(B \subset B_3\).

\begin{corollary}
Inside
our bundle \(B_3 \to S^{2n+1}\) there is a smaller
bundle \(B \to S^{2n+1}\) of points at which
the functions \(s^0_{\bar{p}}\) vanish.
The bundle \(B\) is a principal
\(\Gamma_1\) bundle where \(\Gamma_1\) is
the group of complex matrices of the
form
\[
\begin{pmatrix}
g^0_0 & g^0_q \\
0 & g^p_q
\end{pmatrix}
\]
where \(g^0_0>0\) and 
\(\left(g^0_0\right)^2 \left|\det\left(g^p_q\right)\right|^2=1\).
(This group is identified with a subgroup
of \(G_3\) in the manner outlined
in subsection~\vref{subsec:Hopf}.)
\end{corollary}
\begin{proof} This is the
Lie subgroup of \(G_3\) satisfying
\(\Omega^0_{\bar{q}}=0\).
\end{proof}

Henceforth, we will forget all of
the other \(B_j\) bundles, and
work exclusively with the bundle
\(B\).

Our equations~\vref{eqn:Putty}
simplify to
\begin{equation}\label{eqn:PuttyTwo}
\begin{pmatrix}
\Omega^0_{\bar{0}} \\
\Omega^p_{\bar{0}} \\
\Omega^p_{\bar{q}} 
\end{pmatrix}
= 
\begin{pmatrix}
s^0_r & 0 \\
s^p_r & 0 \\
s^p_{\bar{q} r} & s^p_{\bar{q} \bar{r}}
\end{pmatrix}
\begin{pmatrix}
\Omega^r_{0} \\
\Omega^{\bar{r}}_{\bar{0}} 
\end{pmatrix}
\end{equation}
with \(s^p_p=0.\)
But \(\Omega^0_{\bar{q}}\)
is semibasic on \(B\),
and in fact from the same
equations, we see that
\[
\Omega^0_{\bar{p}} = s^0_{\bar{p}q}\Omega^q_0 +
s^0_{\bar{p}\bar{q}}\Omega^{\bar{q}}_{\bar{0}}
\]
for some functions \(s^0_{\bar{p}q},s^0_{\bar{p}\bar{q}}\)
on the bundle \(B,\) with symmetries in the
lower indices. The invariants 
\(s^0_q\) are related to these by, setting:
\[
\nabla s^0_q = ds^0_q
+ \left(2 \Omega^0_0 - \Omega^{\bar{0}}_{\bar{0}}
\right)s^0_q
-s^0_p \Omega^p_q + \Omega^0_p s^p_q
\]
and calculating:
\[
\nabla s^0_q = s^0_{qr} 
\Omega^r_0 + s^0_{q \bar{r}} \Omega^{\bar{r}}_{\bar{0}}
\]
with \(s^0_{qr}=s^0_{rq}\).
If we define
\[
\nabla s^p_q \defeq
ds^p_q + \Omega^p_r s^r_q - s^p_r \Omega^r_q
+ \left( \Omega^0_0 - \Omega^{\bar{0}}_{\bar{0}} \right) s^p_q
\]
we find
\[
\nabla s^p_q = \left( s^p_{qr} - \delta^p_q s^0_r \right ) \Omega^r_0 
+
\left(-s^p_t s^t_q s^{\bar{0}}_{\bar{r}}
+
s^p_{\bar{r}q}
- s^p_t s^t_{\bar{u}q}s^{\bar{u}}_{\bar{r}}
\right) \Omega^{\bar{r}}_{\bar{0}}.
\]
Differentiating the last of our structure equations
we obtain the covariant derivatives
\begin{align*}
\nabla s^p_{\bar{q}t} &=
ds^p_{\bar{q}t} 
-s^p_t \Omega^{\bar{0}}_{\bar{q}}
+ \Omega^p_r s^r_{\bar{q}t}
-s^p_{\bar{r}t} \Omega^{\bar{r}}_{\bar{q}}
+s^p_{\bar{q}t} \Omega^0_0
-s^p_{\bar{q}r} \Omega^r_t \\
\nabla s^p_{\bar{q}\bar{t}}
&=
ds^p_{\bar{q}\bar{t}}
+ \Omega^p_r s^r_{\bar{q}\bar{t}}
-s^p_{\bar{r}\bar{t}} \Omega^{\bar{r}}_{\bar{q}}
+s^p_{\bar{q}\bar{t}} \Omega^{\bar{0}}_{\bar{0}}
-s^p_{\bar{q}\bar{r}} \Omega^{\bar{r}}_{\bar{t}}
\end{align*}
and we find these satisfy
\begin{align*}
\nabla s^p_{\bar{q} t} +
\tau^p_{\bar{q} t}
&=
s^p_{\bar{q}tu} \Omega^u_0 + s^p_{\bar{q}t\bar{u}} \Omega^{\bar{u}}_{\bar{0}}
\\ 
\nabla s^p_{\bar{q} \bar{t}} +
\tau^p_{\bar{q}\bar{t}} 
&=
s^p_{\bar{q}\bar{t}u} \Omega^u_0 + s^p_{\bar{q}\bar{t}\bar{u}} \Omega^{\bar{u}}_{\bar{0}}
\end{align*}
where the functions \(s^{\cdot}_{\cdot \cdot \cdot}\) on the right
hand sides are symmetric in all lower indices,
and the \(\tau\) 1-forms are given by
\begin{align*}
\tau^p_{\bar{q} t}
&=
\left (
s^0_{\bar{q}t} \delta^p_u 
-
s^p_{\bar{q}\bar{r}} s^{\bar{r}}_{vu} s^v_t 
\right )
\Omega^u_0 
+
\left (
s^p_{\bar{q}\bar{r}}s^{\bar{r}}_{\bar{u}}s^0_t
+s^p_{\bar{q}\bar{r}}s^{\bar{r}}_{v\bar{u}} s^v_t 
\right )
\Omega^{\bar{u}}_{\bar{0}}
\\
\tau^p_{\bar{q}\bar{t}}
&=
\left(
s^0_{\bar{q}\bar{t}} \delta^p_u
-s^p_{\bar{q}r}s^r_u s^{\bar{0}}_{\bar{t}}
-s^p_{\bar{q}r}s^r_{\bar{v}u} s^{\bar{r}}_{\bar{t}}
\right )
\Omega^u_0
-
s^p_{\bar{q}r} s^r_{\bar{v}\bar{u}}
s^{\bar{v}}_{\bar{t}} \Omega^{\bar{u}}_{\bar{0}}
\end{align*}

Applying the Cartan--K{\"a}hler theorem \cite{Bryant:2014}, these structure equations are involutive with general solution depending on \(2n\) functions of \(2n\) variables.
We will see this from another point of view below. 
Since the equations are involutive, there is no need to proceed further along the path of the method of the moving frame; no further local invariants will appear, except for covariant derivatives of the invariants we have already found.

\section{Analogy with complex projective structures}
{\label{sec:Analogy}

Let us once again (for the last time) consider another notation: define
\begin{align*}
\Omega^p &\defeq \Omega^p_0 \\
\Gamma^p_q &\defeq \Omega^p_q - \delta^p_q \Omega^0_0 \\
\Omega_p &\defeq \Omega^0_p.
\end{align*}

Calculating exterior derivatives using the equations we have derived so far gives:
\[
d \Omega^p 
= 
- \Gamma^p_q \wedge \Omega^q 
- \left(s^p_r s^{\bar{0}}_{\bar{t}}
+s^p_{\bar{q}r}s^{\bar{q}}_{\bar{t}}\right)
\Omega^r\wedge\Omega^{\bar{t}}
-s^p_{\bar{q}\bar{r}}s^{\bar{q}}_{\bar{t}}\Omega^{\bar{r}}
\wedge\Omega^{\bar{t}}
\]
and
\begin{align*} 
d \Gamma^p_q =& - \Gamma^p_r \wedge \Gamma^r_q
+\left(\delta^p_q\Omega_t+\Omega_q\delta^p_t\right)\wedge\Omega^t \\
&-\left(s^p_rs^{\bar{0}}_{qt} + s^p_{\bar{u}r}s^{\bar{u}}_{qt}\right)
\Omega^r\wedge\Omega^t\\ 
&-\left(s^p_rs^{\bar{0}}_{q\bar{t}} 
+s^p_{\bar{u}r}s^{\bar{u}}_{q\bar{t}} 
-s^p_{\bar{u}\bar{t}}s^{\bar{u}}_{qr}
-\delta^p_q s^0_r s^{\bar{0}}_{\bar{t}} 
-\delta^p_q s^0_{\bar{u}r}s^{\bar{u}}_{\bar{t}} 
\right)\Omega^r\wedge\Omega^{\bar{t}} \\
&-\left(
s^p_{\bar{u}\bar{r}}s^{\bar{u}}_{q\bar{t}}
-\delta^p_q s^0_{\bar{u}\bar{r}} s^{\bar{u}}_{\bar{t}}
\right)\Omega^{\bar{r}}\wedge\Omega^{\bar{t}}
\end{align*}
and
\begin{align*}
d \Omega_p =& 
\Gamma^q_p \wedge \Omega_q
-\left(s^0_r s^{\bar{0}}_{pq}+s^0_{\bar{t}r}s^{\bar{t}}_{pq}\right)
\Omega^r\wedge\Omega^q
\\&
-\left(s^0_r s^{\bar{0}}_{p\bar{q}} 
+ s^0_{\bar{t}r}s^{\bar{t}}_{p\bar{q}}
- s^0_{\bar{t}\bar{q}}s^{\bar{t}}_{pr}\right)
\Omega^r\wedge\Omega^{\bar{q}}
\\&
-s^0_{\bar{t}\bar{r}}s^{\bar{t}}_{p\bar{q}}
\Omega^{\bar{r}}\wedge\Omega^{\bar{q}}.
\end{align*}
Modulo the various \(s\) functions, i.e.
the torsion, these are the equations
of a flat holomorphic projective structure
on the base manifold \(S^{2d+1} \to X\)
of our great circle fibration.
However, when we include the \(s\) functions,
we find that they are not the structure
equations of a projective structure
at all---in fact the base manifold
is only equipped with an almost complex
structure, which we will soon see.

\subsection{Invariantly defined vector bundles
on the base of a great circle fibration}\label{subsec:Later}

\begin{lemma}
Suppose that \(S^{2n+1} \to X^{2n}\) is our
great circle fibration.
The bundle \(B \to X\) is a principal
\(\Gamma_0\) bundle where (as in 
subsection~\vref{subsec:Hopf})
\(\Gamma_0\) is the group
of matrices of the form 
\[
\begin{pmatrix}
g^0_0 & g^0_q \\
0 & g^p_q
\end{pmatrix}
\]
where all entries are complex numbers
and \(|g^0_0 \det\left(g^p_q\right)|^2=1\). 
The representation 
\[
\begin{pmatrix}
\Gamma^p_q & 0 \\
0 & \Gamma^{\bar{p}}_{\bar{q}} 
\end{pmatrix}
\]
solders
the tangent bundle of \(X^{2n}\).
There is an invariantly defined
almost complex structure on \(X\)
whose holomorphic tangent space
\(T^{1,0} X\) is soldered by
\(\Gamma^p_q.\)
The 1-form \(\Omega^0_0\) solders
the principal 
bundle \(S^{2n+1} \to X\).
\end{lemma}
\begin{proof}
We need
to show that the fibers of
\(B \to X\) are connected.
But this factors as
\(B \to S^{2n+1} \to X\)
so that the fibers of the
first map are copies of \(\Gamma_1\),
which is connected, and the fibers
of the second map are circles,
hence connected.

The 1-forms \(\Omega^p\) are semibasic
for the projection to \(X\).
On the fibers of \(B \to X\)
the structure equations reduce
to the structure equations of
\(\Gamma_0\) and its left translates.
As before, this shows by connectedness
of \(\Gamma_0\) and of the fibers
of \(B \to X\) that the
fibers are in fact left translates
of \(\Gamma_0\). Therefore \(B \to X\)
is a principal right \(\Gamma_0\) bundle.

The soldering forms of the tangent bundle
of \(X\) are unchanged from lemma~\vref{lm:Porridge},
but rewritten in complex notation.
The almost complex structure
on \(X\) is immediately visible
from the structure equations,
since the \(\Gamma^p_q\) 1-forms
solder in a complex representation.

I have not explained what it
means to solder a principal
bundle out of another one,
but it should be obvious. The
soldering of \(\Omega^0_0\)
is just the quotient by \(\Gamma_1\)
of the soldering of \(\Gamma_0\),
which is just \(S^{2n+1}\to X\).
\end{proof}

The base manifold \(X\) is analogous
to a complex projective space,
and the circle fibration \(S^{2n+1} \to X\)
analogous to the Hopf fibration.
Therefore henceforth we will refer to the
complex line bundle soldered
(on the bundle \(B \to X\))
by \(\Omega^0_0\)  as \(\OO{-1}\),
and similarly define
the line bundles 
\[
\OO{p} \defeq \OO{-1}^{\otimes (-p)}.
\]
Be careful to note that these are
bundles on \(X\) and are 
complex line bundles. We will never
again refer to the similarly
named real line bundles on \(S^{2n+1}\),
which were introduced simply
to encourage an analogy between
the soldering of these \(\OO{p} \to X\)
bundles with the similarly named
bundles on projective
spaces.

From here on, we will use the notation \(\tilde{W}\)
for any complex representation \(W\) of \(\Gamma_0\)
to mean the vector bundle \(\tilde{W} \to X\) given as
\[
\tilde{W} 
= 
B \times^{\Gamma_0} W 
= 
\left(B \times W\right)/\Gamma_0.
\]
constructed in the same manner
as in subsection~\vref{subsec:VBS}.
For example, we will introduce
the complex vector bundle \(\tilde{V}\)
out of the complex representation
of \(\Gamma_0\) given by the
identity representation (recall that elements
of \(\Gamma_0\) are complex matrices).

\begin{lemma}
\[
T^{1,0} X = \OO{1} \otimes \left ( \tilde{V} / \OO{-1} \right ). 
\]
\end{lemma}
\begin{proof}
This is immediate from working
out the soldering in 1-forms.
\end{proof}

The equations satisfied on the bundle \(B\)
by a section of \(\tilde{V}\) are
\[
d
\begin{pmatrix}
f^0 \\
f^{\bar{0}} \\
f^{p} \\
f^{\bar{p}}
\end{pmatrix}
=
\begin{pmatrix}
\Omega^0_0 & \Omega^0_{\bar{0}} & \Omega^0_q & \Omega^0_{\bar{q}} \\
\Omega^{\bar{0}}_0 & \Omega^{\bar{0}}_{\bar{0}} 
        & \Omega^{\bar{0}}_q & \Omega^{\bar{0}}_{\bar{q}} \\
\Omega^p_0 & \Omega^p_{\bar{0}} & \Omega^p_q & \Omega^p_{\bar{q}} \\
\Omega^{\bar{p}}_0 & \Omega^{\bar{p}}_{\bar{0}} 
        & \Omega^{\bar{p}}_q & \Omega^{\bar{p}}_{\bar{q}} 
\end{pmatrix}
\begin{pmatrix}
f^0 \\
f^{\bar{0}} \\
f^q \\
f^{\bar{q}}
\end{pmatrix}.
\]
Plugging in our equations on \(\Omega^0_{\bar{0}}\)
and \(\Omega^0_{\bar{q}}\) we might first get
the impression that \(f^0\) is a holomorphic
section of some line bundle. But the \(\Omega^0_{q}\)
factor is not semibasic, so this is not 
actually soldering any line bundle.

%Although
%\(\tilde{V}\) is topologically trivial,
%it is not obvious that it is trivial
%as a complex vector bundle. Any complex line
%bundle which is topologically trivial
%as a real rank 2 bundle is also 
%trivial as a complex line bundle,
%since the space of complex structures
%on a real 2-plane of fixed orientation
%is contractible. But this is no
%longer the case in higher dimensions.
%The universal example is 
%\(\Det\tilde{V}\to\Cstrucs{V}\).
%I still don't understand that example; 
%although it has nonzero
%curvature, that is not incompatible
%with flatness, because the curvature
%might not be of fixed sign, and
%more importantly the base space
%\(\Cstrucs{V}\) is not compact.

\section{Homogeneous great circle fibrations}\label{sec:Homog}

A symmetry of a great circle fibration
which preserves the flat projective
structure on the sphere 
acts as a symmetry of the right 
principal bundle \(B \to X\).
Suppose that the symmetry group
acts transitively on \(B\). Then
the invariantly defined functions
\(s^0_q,s^p_q\) etc. are all
constant on \(B\). Plugging this hypothesis
into our equations for covariant derivatives,
we find that all of the invariants
vanish, and applying the Frobenius theorem,
along with simple connectivity of \(X\),
we see directly that the
great circle fibration is a Hopf fibration.

\begin{theorem}
The symmetry group of a great circle
fibration acts transitively on the
adapted frame bundle precisely
if the fibration is a Hopf fibration.
\end{theorem}

\begin{theorem}
The 
symmetry group of a great circle fibration
injects into the bundle
\(B\), so that it is always of dimension
at most that of \(B\), i.e.
\[
\dim \Aut \le 2(n+1)^2-1.
\] 
Equality occurs only for the Hopf fibration.
\end{theorem}
\begin{proof}
If 
\(\phi \colon B \to B\) is a symmetry, then
by definition it preserves great circles, so \(\phi=g \in \SL{V}\).
If we arrange that the identity element \(I \in \SL{V}\)
belongs to \(B\), then \(gI = g \in B\), so 
in fact the diffeomorphism \(\phi\) is
an element of \(B\). Therefore the symmetry
group is actually a subgroup of \(\SL{V}\)
sitting inside \(B\). The symmetry group is closed, since the
condition of being a symmetry is a
closed condition. 
Therefore it is 
a Lie subgroup of \(\SL{V}\) lying inside \(B\).

If the symmetry group has the same dimension as
\(B\), then it is an open subset of \(B\),
say \(U\), and a closed subgroup of \(\SL{V}\).
But \(B \subset \SL{V}\) is also closed,
so \(U\) is both open and closed in \(B\),
and therefore a union of path components
of \(B\). But \(B\) is connected, since it
is a bundle
\[
\begin{tikzcd}
\Gamma_1 \arrow{r} & B \arrow{d} \\
                & S^{2n+1} 
\end{tikzcd}
\]
with connected fibers and base.
Therefore the symmetry group is precisely \(B\).
This forces all of the invariant functions \(s\)
on \(B\) to be constants, and then the
structure equations force them to vanish.
\end{proof}

\section{Embedding into the Grassmannian of oriented 2-planes}\label{sec:Embedding}

Recall that \(G=\SL{V}\) and \(\Gcircle\)
is the subgroup preserving the oriented great
circle on \(S^{2n+1}\) passing from
\(e_0\) through to \(e_1\).
The space \(G/\Gcircle\)
is naturally identified with the
space of all oriented great circles
on the sphere, or with the
space of oriented 2-planes in \(V\):
\[
G/\Gcircle=\Gro{2}{V}.
\]
We find that \(\Gamma_0 \subset \Gcircle\)
so that a fiber bundle
\[
G/\Gamma_0 \to G/\Gcircle=\Gro{2}{V}
\]
is defined, with fibers being
homogeneous spaces of \(\Gcircle\)-equivariantly 
diffeomorphic to \(\Gcircle/\Gamma_0\)
which is the space of all complex
structures on \(V\) which have the oriented
2-plane \(\Line{e_0,e_1}\) as a complex line.

We can map
\[
X = B/\Gamma_0 \to G/\Gamma_0 \to G/\Gcircle=\Gro{2}{V}.
\]

\begin{lemma}
This map \(X \to \Gro{2}{V}\) is an embedding.
\end{lemma}
\begin{proof}
If two points \(x,y \in X\) get mapped
to the same place, then they correspond
to the same great circle, since \(X\) 
parameterizes the great circles
of our fibration, and \(\Gro{2}{V}\)
is parameterizing all great circles.
So the map is injective.

We have to differentiate to see
that the map is an immersion. Pull back a local coframing 
by 1-forms on \(\Gro{2}{V}\)
to \(X\) and show that it contains a
coframing for \(X\). But then it is 
sufficient to pull the coframing
back to \(B\) and show that on \(B\)
every semibasic 1-form for the bundle
map \(B \to X\) can be expressed
as a linear combination (with 
real functions as coefficients)
of the 1-forms from the coframing from \(\Gro{2}{V}\).

The semibasic 1-forms for the
map \(G \to \Gro{2}{V}\) are precisely
those complimentary to the Lie
algebra of \(\Gcircle\),
i.e. they are the 1-forms
\(
\Omega^p_0, \Omega^p_{\bar{0}}
\)
and their complex conjugates.
Therefore any coframing on \(\Gro{2}{V}\)
pulls back to \(G\) to be a
combination of these 1-forms,
and conversely they are combinations
of the 1-forms from the coframing.
But pulling back to \(B\), we find
only the relations 
\(
\Omega^p_{\bar{0}} = s^p_q \Omega^q_0.
\)
The remaining semibasic 1-forms
\(\Omega^p_0\) span the semibasic
1-forms on \(X\).
Therefore the map is an immersion.
An injective immersion of a compact
manifold is an embedding.
\end{proof}

A consequence of the proof is that
these \(s^p_q\) invariants can be expressed
in terms of the first derivative of
the embedding \(X \to \Gro{2}{V}\).
The embedding determines
the great circle fibration. This embedding
is a more useful way
to examine great circle fibrations
than looking directly at the sphere,
because a perturbation of great circles
on the sphere can never be local,
while it can be local on \(X\),
as a small motion of \(X\) inside
the Grassmannian.

\section{Characterizing the submanifolds of the Grassmannian 
which represent great circle fibrations}\label{sec:Charac}

Consider an immersed submanifold \(\iota \colon X \hookrightarrow \Gro{2}{V}\)
in the Grassmannian of oriented 2-planes
of a vector space \(V\). Assume that \(\dim V = 2n+2\)
and \(\dim X = 2n\).
Our next problem is to characterize
when \(X\) represents a great circle
fibration. The tangent spaces
to \(\Gro{2}{V}\) are canonically identified
with 
\[
T_P \Gro{2}{V} \cong \Lin{P}{V/P}.
\]
So \(T_P X\) is identified with a linear subspace
\[
T_P X \cong U_P \subset \Lin{P}{V/P}.
\]
We have a map 
\(
\alpha_P \colon P \to \Lin{U_P}{V/P}
\)
defined by
\begin{equation}\label{eqn:Puke}
\alpha_P(p)(u)\defeq u(p) \text{ for } p \in P \text{ and } u \in U_P.
\end{equation}
Note that for each \(p \in P\)
the map
\[
\alpha_P(p) \in \Lin{U_P}{V/P}
\]
is a linear transformation
between vector spaces of the same
dimension.
For each \(p \in P\) we can define
the polynomial
\[
\xi_P(p) \defeq \det \alpha(p) 
\in \Lin{\Det\left(U_P\right)}{\Det\left(V/P\right)}.
\]
which is a polynomial not valued in real numbers,
but in the one dimensional vector 
space \( \Lin{\Det\left(U_P\right)}{\Det\left(V/P\right)}. \)
We define the \emph{characteristic
variety} \(\Xi_P\) to be the 
projective variety in the 
projective line \(\mathbb{CP}(P\otimes_{\R{}}\C{})\)
cut out by \(\xi_P\).

An immersed submanifold \(\iota \colon X^{2n} \hookrightarrow \Gro{2}{V}\)
(where \(\dim V = 2n+2\)) is \emph{elliptic} at a point \(P \in X\)
if the characteristic variety \(\Xi_P\) has no real points.

We have \(\RP{}(P) \subset \CP{}(P)\)
a real circle on a real sphere, cutting it
into two halves. The orientation of \(P\)
chooses one of these halves: \(P\) being
oriented orients \(\RP{}(P)\)
and \(\CP{}(P)\) is oriented by its complex structure.
This coorients \(\RP{}(P)\) and picks a half of \(\CP{}(P) - \RP{}(P)\).
The characteristic variety does intersect \(\RP{}(P)\), so consists in a finite set of points lying on that half, and their conjugates on the other half. 
The linear fractional transformations we used are just real reparameterizations
of \(P\). 
These enable us to move these points of the characteristic variety by isometries of the Poincar{\'e} metric. 
We choose to move the points to obtain vanishing trace of the Reznikov \(s^p_q\) invariant. 
This picks a specific choice of linear map from the family parameterized
by \(P\).

The Grassmannian \(\Gro{2}{V}\) is a homogeneous \(G=\SL{V}\) space, and we can write \(\Gro{2}{V} = G/\Gcircle.\)
Consider the pullback bundle
\[
\begin{tikzcd}
\iota^*G \arrow{d} \arrow{r} & G \arrow{d} \\
X \arrow{r} & \Gro{2}{V}
\end{tikzcd}
\]
which is a principal \(\Gcircle\)
bundle. On \(G\) we have our old
Maurer--Cartan 1-forms
\[
\begin{pmatrix}
\omega^0_0 & \omega^0_1 & \omega^0_j \\
\omega^1_0 & \omega^1_1 & \omega^1_j \\
\omega^i_0 & \omega^i_1 & \omega^i_j
\end{pmatrix}.
\]
The group \(\Gcircle\) is
precisely the connected subgroup
of \(G\) satisfying
\(
\omega^i_0 = \omega^i_1 = 0.
\)
Therefore the 1-forms \(\omega^i_0,\omega^i_1\)
span the semibasic 1-forms for the map \(G \to \Gro{2}{V}\).
On the pullback bundle over \(X\), these 1-forms
are no longer independent, since there 
are fewer degrees of freedom along \(X\),
in fact half as many, since \(\dim X = \frac{1}{2} \dim \Gro{2}{V}.\)

We need to consider how to express
the \(G\)-invariant identification
\[
T_P \Gro{2}{V} \cong \Lin{P}{V/P}
\]
in terms of these 1-forms.
Recall how the identification
is constructed: take any family
of 2-planes \(P(t) \in \Gro{2}{V}\),
and any family of linear maps
\(\phi(t) \colon V \to W\) for some
fixed vector space \(W\) of dimension
\(\dim W = \dim V - \dim P(t)\),
with the maps \(\phi(t)\) chosen
so that \(\ker \phi(t) = P(t)\).
Write \(\bar{\phi}(t)\) for the
induced map 
\[
\bar{\phi}(t) \colon v+P \in V/P(t) \to \phi(t)(v) \in W
\]
which is defined because \(\ker \phi(t) = P(t)\).
Then identify 
\[
P'(t) \sim \bar{\phi}(t)^{-1} \left.\phi'(t)\right|_{P} \colon P \to V/P. 
\]
This is well defined because if \(\psi(t)\)
is any other choice of maps replacing
\(\phi(t)\), with the same kernel, 
\[
\begin{tikzcd}
0 \arrow{r} & P(t) \arrow{r} & V \arrow{r}{\psi(t)} & U \arrow{r} & 0 
\end{tikzcd}
\]
then differentiating the equation
\[
\psi(t) = \bar{\psi}(t)\bar{\phi}(t)^{-1} \phi(t)
\]
shows that
\[
\bar{\phi}(t)^{-1} \left.\phi'(t)\right|_{P}
=
\bar{\psi}(t)^{-1} \left.\psi'(t)\right|_{P}
\]
so that the resulting map in \(\Lin{P(t)}{V/P(t)}\)
is independent of the choice of map \(\phi(t)\).

So far this only defines a map
\[
T_P \Gro{2}{V} \to \Lin{P}{V/P}.
\]
We need to pick some local coordinates
on \(\Gro{2}{V}\), and we will take
them as follows: for any 2-plane \(P_0 \in \Gro{2}{V},\)
we take coordinates \(x^1,x^{2},y^1,\dots,y^{2n}\)
on \(V\)
so that \(P_0\) is cut out by \(y=0.\)
Then 2-planes near \(P_0\) are cut out
by equations like \(y=px\) where \(p\)
is a \(2n \times 2\) matrix. These \(p\)
are our local coordinates on \(\Gro{2}{V},\)
and one can easily compute the
transformations of coordinates 
if we change the  choice of \(P_0\)
and the choice of coordinates \(x,y\).
In these coordinates, we can take
\(\phi_P(x,y)=px-y\) as our map, with
kernel \(P\). Then we find that for
any family \(P(t)\) with \(P(0)=P_0\)
we have \(p(0)=0\) and 
\[
\bar{\phi}(0)^{-1}\phi'(0)x=-p'(0)x
\] 
or
\[
\bar{\phi}(0)^{-1}\phi'(0)=-p'(0).
\]
Therefore this map
\[
T_P \Gro{2}{V} \to \Lin{P}{V/P}
\]
is an isomorphism. Equivariance
under linear transformations
of \(V\) is obvious.

Returning to the 1-forms \(\omega^i_0\) and
\(\omega^i_1\), recall that the
fiber of the bundle \(G \to \Gro{2}{V}\)
over a point \(P \in \Gro{2}{V}\)
consists precisely of the elements of
\(G\) taking the plane \(\Line{e_0,e_1}\) to the
plane \(P\), preserving
orientation. Therefore we can define
a map
\[
\phi_g \colon V \to V/\Line{e_0,e_1}
\]
with kernel \(P\) by
\[
\phi_g(v) \defeq g^{-1}v + \Line{e_0,e_1}.
\]
We find 
\[
\bar{\phi}^{-1} \, d \phi = -dg \, g^{-1} + P
\]
or
\begin{align*}
\bar{g}^{-1}\left( \bar{\phi}^{-1} \left. d\phi \right|_{P} \right ) g
&= - \left. \omega \right|_{\Line{e_0,e_1}} + \Line{e_0,e_1} \\
&= - 
\begin{pmatrix} 
\omega^i_0 & \omega^i_1 
\end{pmatrix}.
\end{align*}
Perhaps a little more concretely,
if \(Q \in TG\) is a tangent vector,
and \(a \in \Line{e_0,e_1}\), we have
\begin{equation}\label{eqn:Barf}
\bar{g}^{-1}\left( \bar{\phi}^{-1} \left. d \phi(Q) \right|_P \right ) ga
= - \left( \omega^i_0(Q) a^0 + \omega^i_1(Q) a^1 \right).
\end{equation}

\begin{lemma}
If \(X \subset \Gro{2}{V}\) is
the base manifold of a great
circle fibration, then \(X\)
is elliptic at every point.
\end{lemma}
\begin{proof}
Pick a point \(P \in X\) and
a point \(g \in B\) which is
taken by \(B \to X\) to \(P.\)
As expressed in terms of \(\omega^i_0,\omega^i_1\)
above, we found earlier that a great circle fibration
satisfies \(\omega^i_1 = t^i_j \omega^j_0,\)
with \(t^i_j\) a real \(2n \times 2n\) matrix
with no real eigenvalues.
For each vector \(u \in T_P X\)
we can write its components in the
coframe \(\omega^i_0\) as
\(u^i\). (The reader who
is keeping track of what spaces
we are working in will be puzzled
to read that the \(\omega^i_0\) 
are a coframe on \(X\). What we 
mean of course is that since the \(\omega^i_0\)
are semibasic, at each point \(g \in B\)
we can form the coframe \(\underline{\omega}^i_0\)
on \(T_P X\) from which these \(\omega^i_0\)
are pulled back. These \(\underline{\omega}^i_0 \in \Lm{1}{T_P X}\)
change as we move up the fibers of \(B \to X\).) Each point 
\(a^0e_0 +a^1e_1 \in \Line{e_0,e_1}\)
is carried by \(g\) to an element of \(P\)
and we identify these.
We find that in terms of equation~\vref{eqn:Barf},
the map \(\alpha_P\) defined in equation~\vref{eqn:Puke}
is given, up to factors of \(g\) and \(g^{-1}\)
(which won't affect the vanishing of the
relevant determinant) by
\[
\alpha\left(a^0e_0+a^1 e_1\right)(u)
=a^0 u^i + t^i_j a^1 u^i = \left(a^0 \delta^i_j + a^1 t^i_j \right )u^j.
\]
Therefore the polynomial \(\xi_P\) is
\[
\xi_P\left(a^0,a^1\right)
=
\det
\left ( a^0 \delta^i_j + a^1 t^i_j \right ).
\]
If \(a^1=0\) then a zero of \(\xi_P\) can
only occur at \(a^0=0\). Therefore, taking
\(a^1 \ne 0\), we find that 
\[
\xi_P = \left(a^1\right)^{2n} \det \left ( t^i_j +\frac{a^0}{a^1} \delta^i_j 
\right ).
\]
A real zero of \(\xi_P\) therefore corresponds precisely
to a real eigenvalue of \(t^i_j\). Consequently,
if \(X\) arises from a great circle fibration,
then \(X\) is elliptic.
\end{proof}

We will now take an elliptic immersed submanifold
\(X \subset \Gro{2}{V}\) and
apply the method of the moving
frame to calculate its structure equations.

\begin{lemma}\label{lemma:Popsicle}
Suppose that \(\iota \colon X \to \Gro{2}{V}\)
is an immersion of an elliptic
submanifold. Then there is an
invariantly defined 
principal right \(\Gamma_0\)
subbundle \(B \to X\) and a \(\Gamma_0\)-equivariant map
\[
\begin{tikzcd}
B \arrow{r} \arrow{d} & \SL{V} \arrow{d} \\
X \arrow{r} & \Gro{2}{V}
\end{tikzcd}
\]
which satisfies the structure equations
of a great circle fibration of the 
sphere \(S^{2n+1}\), and \(X\) is locally a great circle fibration.
\end{lemma}
\begin{proof}
Consider the bundle
\[
\begin{tikzcd}
\iota^* G \arrow{d} \arrow{r} & G \arrow{d} \\
X \arrow{r} & \Gro{2}{V}.
\end{tikzcd}
\]
As before we will write \(\phi_g(v) = g^{-1}v + \Line{e_0,e_1}\)
giving a map
\[
\phi_g \colon V \to V/\Line{e_0,e_1}
\]
with kernel \(P=g \Line{e_0,e_1}.\) 
Again we have the equation
\[
-\bar{g}^{-1} \left ( \bar{\phi}^{-1} d \left.\phi\right|_{P} \right )
ga = \omega^i_0 a^0 + \omega^i_1 a^1
\]
for \(a \in \Line{e_0,e_1}.\) We know that the
\(V/\Line{e_0,e_1}\) valued 1-form \(\omega^i_0 a^0 + \omega^i_1 a^1\)
is a coframing on \(T_P X\) for each \(a \ne 0\), 
precisely expressing ellipticity.
So \(\omega^i_1\) is a coframing, as is \(\omega^i_0\).
Consequently there is an invertible linear
transformation \(t^i_j\) so that
\(
\omega^i_1 = t^i_j \omega^i_0.
\)
We must have
\[
\omega^i_0 a^0 + \omega^i_1 a^1 
= 
\left ( a^0 \delta^i_j + a^1 t^i_j \right ) \omega^i_0
\]
also a coframing, as long as \(a \ne 0\).
This says precisely that \(t^i_j\) has
no real eigenvalues.

After this, we have obtained
exactly the same structure equations
as in equation~\vref{eqn:ColdCereal}.
Then we repeat the entire development
of those structure equations, identically.

Given the resulting structure
equations, we wish to construct a
local great circle fibration out
of a portion of \(X\).
Consider the circle bundle 
\(\Sigma \to X\) consisting
of pairs \((P,\ray{v})\) where 
\(P \in X\) and \(v \in P \subset V\)
with \(v \ne 0\)
and as usual \(\ray{v}\) means
\(v\) up to positive rescaling.
Consider the map
\(
\Phi \colon \Sigma \to S^{2n+1}
\)
given by
\(
\Phi\left(x,\ray{v}\right) = \ray{v}.
\)
We wish to show that \(\Phi\) is a local
diffeomorphism. This \(\Sigma\) is
the principal circle bundle associated
to the complex line bundle \(\OO{-1}_X\).
The isotropy group of a point of \(\Sigma\)
inside the structure group of \(B \to X\)
is precisely the group \(\Gamma_1.\)
We have maps
\[
\begin{tikzcd}
B \arrow{d} \arrow{r} & G \arrow{d} \\
\Sigma \arrow{d} \arrow{r} & S^{2n+1} \\
X.
\end{tikzcd}
\]
Working out the semibasic 1-forms
for the maps \(G \to S^{2n+1}\)
and \(B \to \Sigma\), we find
they are identical. Therefore the
map is a local diffeomorphism.
The fibers of \(B \to \Sigma\)
are contained in the left translates
of the group \(\Gcircle\) by
structure equations, so these fibers
sit in great circles on the sphere.
\end{proof}

Taking any matrix \(t^i_j\) with no
real eigenvalues, build the subspace
of linear maps of the form
\[
a \in \R{2} \mapsto \left (a^0 \delta^i_j + a^1 t^i_j \right ) u^j
\]
for \(u \in \R{2n}\) and you have a
linear subspace inside \(\Lin{\R{2}}{\R{2n}}\),
which you think of as a subspace of some 
tangent space
to a Grassmannian. Then you can take
any submanifold of the Grassmannian
with that tangent space at that point,
and you have (at least near this point
of the submanifold) an elliptic
submanifold. So there are lots of 
elliptic submanifolds, locally.  

\begin{proposition}
A submanifold of the Grassmannian
is the base of a great circle fibration
precisely if it is elliptic, compact and connected.
\end{proposition}
\begin{proof}
We have seen that the
base of a great circle fibration
is elliptic, compact and connected.
Let \(X \subset \Gro{2}{V}\) be
elliptic, compact and connected.
The map \(\Phi \colon \Sigma \to S^{2n+1}\)
from lemma~\vref{lemma:Popsicle}
is a local diffeomorphism, taking
fibers of \(\Sigma \to X\) to great 
circles. The space
\(\Sigma\) is the total space of 
circle bundle over \(X\), so compact.
Because the sphere \(S^{2n+1}\) is
simply connected, this forces \(\Phi\)
to be a diffeomorphism. The map
\(S^{2n+1} \to X\) is therefore defined,
and satisfies our structure equations,
so is a global great circle fibration.
\end{proof}

\begin{corollary}\label{cor:Pumpkin} Given a great circle
fibration \(S^{2n+1} \to X\), every
\(C^1\) small motion of \(X\) inside
\(\Gro{2}{V}\) is the base of a 
great circle fibration. So there
are lots of great circle fibrations,
and they admit lots of deformations.
The normal bundle \(\nu X\) of \(X\) inside \(\Gro{2}{V}\)
has fibers
\[
\nu_P X = \Lm{0,1}{P}\otimes_{\C{}} V/P
\]
for \(P \in X\) (so \(P \subset V\) a 2-plane)
where the relevant complex structure
on \(P\) and \(V/P\) is \(J_P\).
So \(\nu X = T^{0,1}X\).
The great circle fibrations
near a given great circle fibration
\(S^{2n+1} \to X\) are ``parameterized''
by sections of \(\nu X\) close to the
zero section (for example, by using
a Riemannian metric on \(\Gro{2}{V}\)).
Hence the general great circle
fibration depends on \(2n\) real functions
of \(2n\) real variables.
\end{corollary}

\begin{corollary}
The space of \(C^k\)
great circle fibrations is an infinite dimensional
manifold, for \(k>1\).
\end{corollary}

\section{A taste of elliptic partial differential equations}
{
Given a real surface \(C\) and an immersion \(\phi \colon C \to V\) to a vector space \(V\), construct the map \(\phi_1 \colon C \to \Gro{2}{V}\) giving
the tangent plane.
Pick a great circle fibration \(X \subset \Gro{2}{V}\).
If the image of \(\phi_1\) is always in \(X \subset \Gro{2}{V}\), call \(\phi\) an immersed \(X\)-curve.
The \(X\)-curves are the solutions of an involutive elliptic first order system of 2 partial differential equations for \(2n\) functions of 2 variables. 
Another approach: consider pseudocomplex manifolds, as in
section~\vref{sec:NonlinearJ}.

\section{\texorpdfstring{Inapplicability of methods to prove the $h$-principle}{Inapplicability of methods to prove the h-principle}}{\label{sec:H}
The relation of ellipticity for a submanifold of the Grassmannian is open but not ample in the sense of Gromov's theory of convex integration. 
This makes it unlikely that convex integration can be applied here. 

Because the relation is open, it is microflexible.
I believe that some deformations of elliptic submanifolds nowhere
parallel to a given hinge (see section~\vref{section:hinge}) are \emph{not}
microcompressible (see page 81 of Gromov \cite{Gromov:1986}); if so this makes the sheaf of elliptic submanifolds with given hinge not flexible, and then
the same is true for the sheaf of all elliptic submanifolds.
So probably we can not apply the method of sheaves.

The ellipticity relation is not the complement of a hypersurface in the 1-jet bundle, so we can not apply the method of elimination of singularities. 
It is not defined by a differential operator, so we can not apply the method of 
inversion of differential operators.

\section{The osculating complex structure}\label{sec:Osculating}
Recall that \(J_0\) is our fixed complex structure on \(V\). 
The space of all complex structures on \(V\) with the standard orientation is
\[
\Cstrucs{V}=G/\Gamma=\SL{V}/\pr{\SL{V}\cap\GL{V,J_0}}.
\]
Now consider an elliptic submanifold \(X \subset \Gro{2}{V}.\)
Since \(B \subset G\) we can map
\[
B \to G \to G/\Gamma=\Cstrucs{V}.
\]
In fact, since \(\Gamma_0 \subset \Gamma,\)
we can map 
\[
X = B/\Gamma_0 \to G/\Gamma_0 \to G/\Gamma = \Cstrucs{V}.
\]
We can obviously do better: the space
\(
G/\Gamma_0
\)
is naturally identified with the space
of all pairs \((J,P)\) where \(J\) is a
complex structure on \(V\), and \(P\) is a 
real 2-plane which is a \(J\)-complex line,
i.e. \(JP=P\).
Then we can map
\[
X = B/\Gamma_0 \to G/\Gamma_0
\]
to this space. 
When \(X\) is the base manifold of
a great circle fibration, 
this map
is an embedding, since \(B \subset G\)
is embedded. More generally, it
is an immersion. We will write
for \(P \in X\) the corresponding
point in \(\Cstrucs{V}\) (i.e.
complex structure on \(V\)) as
\(J_P\). So \(P \subset V\) is a \(J_P\) complex
line.

The bundle \(\tilde{V} \to X\)
is actually pulled back from \(\Cstrucs{V}\), 
as the isotropy groups are
\[
\begin{tikzcd}[ampersand replacement=\&]
B \arrow{d}^{\Gamma_0} \arrow{r} \& G \arrow{d}^{\Gamma} \\
X \arrow{r} \& \Cstrucs{V}.
\end{tikzcd}
\]
Since \(\Gamma_0 \subset \Gamma \subset G\), and \(V\)
is a \(G\) module, so a \(\Gamma\) module,
we find that \(\tilde{V}\) is defined
as a bundle on \(\Cstrucs{V}\).
The bundle \(\Det\left(\tilde{V}\right)\) is
defined to be the complex determinant
line bundle, which is defined because 
\(\Gamma\) is a complex Lie group.
This \(\Det\left(\tilde{V}\right)\)
is soldered by
\[
A = \sqrt{-1} \left( \Omega^0_0 + \Omega^p_p \right).
\] 
Differentiating, we find the
curvature of \(\Det\left(\tilde{V}\right)\)
is
\[
F = dA = 
-\sqrt{-1} \Omega^0_{\bar{0}} \wedge \Omega^{\bar{0}}_0
-\sqrt{-1} \Omega^0_{\bar{p}} \wedge \Omega^{\bar{p}}_0 
-\sqrt{-1} \Omega^p_{\bar{0}} \wedge \Omega^{\bar{0}}_p
-\sqrt{-1} \Omega^0_{\bar{q}} \wedge \Omega^{\bar{q}}_p
\]
which is a \((1,1)\) form on \(\Cstrucs{V}\),
invariant under complex conjugation, i.e.
\(
\bar{F} = F.
\)

For each \(J \in \Cstrucs{V}\), let \(X_J \subset \Gro{2}{V}\)
be the base manifold of the associated
Hopf fibration, i.e. \(X_J\) is the
set of \(J\) complex lines, suitably oriented.
\begin{lemma}
Recall that there is a canonical 
identification \(\Gro{2}{V} \cong \Lin{P}{V/P}.\)
Suppose that \(X_J \subset \Gro{2}{V}\) 
is the base of a Hopf fibration.
Then \(T_P X_J \subset T_P \Gro{2}{V}\)
is identified with the subset of \(J\)
linear maps. The map \(X_J \to \Cstrucs{V}\)
constructed above is constant, mapping to \(J\).
\end{lemma}
\begin{proof}
By \(\SL{V}\)-invariance, we only have to
prove the result for the standard
complex structure \(J\) on \(V=\C{n+1}\).
The first result is a calculation, while the second
is immediate from the structure equations:
\[
\Omega^{p}_{\bar{0}}=\Omega^{0}_{\bar{p}}=\Omega^p_{\bar{q}}=0
\]
on \(X_J\), and the fact that \(X_J\) is connected.
\end{proof}

Take a great circle fibration \(S^{2n+1} \to X\) and a point
\(P \in X\), and the complex structure \(J_P\). 
The \emph{osculating} Hopf fibration to \(X\) at \(P\) is \(S^{2n+1} \to X_{J_P}\).
Note that 
\[
P \in X_{J_P} \cap X \subset \Gro{2}{V}.
\]

\begin{lemma} The invariant \(\left(s^p_q\right)\) vanishes
at a point \(P \in X\) on an elliptic
submanifold \(X \subset \Gro{2}{V}\)
precisely when the
osculating complex structure \(X_{J_P}\)
at \(P\) is tangent to \(X\) inside \(\Gro{2}{V}\).
\end{lemma}
\begin{proof} We showed already
that this invariant is determined
by the 1-jet of the immersion \(X \to \Gro{2}{V}.\)
\end{proof}

\begin{lemma} An elliptic
submanifold \(X \subset \Gro{2}{V}\)
is locally the base manifold of a Hopf fibration
\(X_J\) precisely if \(X\) all of the
\(s\) invariants vanish:
\[
s^p_q=s^0_r=s^p_{\bar{q}r}=s^p_{\bar{q}\bar{r}}=
s^0_{\bar{p}q}=s^0_{\bar{p},\bar{q}}=0
\]
which happens precisely if the map
\(X \to \Cstrucs{V}\) constructed above is
constant.
\end{lemma}
\begin{proof} The 1-forms
which are semibasic for the map \(G \to \Cstrucs{V}\)
are precisely the \(\Omega\) with no bars, i.e.
\(\Omega^P_Q\), so the \(\Omega^P_{\bar{Q}}\) are not
semibasic. This means that we pull them back
when we differentiate the map \(X \to \Cstrucs{V}\),
and nothing else. Hence the map \(X \to \Cstrucs{V}\)
has vanishing differential, and consequently
is locally constant, precisely when all of
these invariants vanish. Conversely, if they
all vanish, then our structure equations
become the same as for the Hopf fibration,
and the result follows by the Frobenius theorem.
\end{proof}

\begin{corollary} A submanifold \(X \subset \Gro{2}{V}\)
is the base of a Hopf fibration 
precisely when it is compact and connected
and all of the \(s\) invariants vanish.
\end{corollary}

We have an inclusion map
\[
\Cstrucs{V} \to \operatorname{GCF}(V)
\]
(where \(\operatorname{GCF}(V)\) means the 
space of great circle fibrations of the
sphere \((V \setminus 0)/\R{+}\), and
\(\Cstrucs{V}\) is the space of complex
structures) given
by mapping a complex structure to its
Hopf fibration. Conversely, if we pick
any nonzero vector \(v \in V\), we can 
associate to any great circle fibration
\(\pi \colon S^{2n+1} \to X\) the osculating complex
structure \(J_v^X\) to \(X\) at the point \(x=\pi\left(v\right)\).
The diagram
\[
\begin{tikzcd}[ampersand replacement=\&]
\& \operatorname{GCF}(V) \arrow{dr}{J_v} \& \\
\Cstrucs{V} \arrow{ur} \arrow{rr}{\text{id}} \& \& \Cstrucs{V} 
\end{tikzcd}
\]
which sheds some light on the algebraic
topology of the space of great circle
fibrations. 

\section{Recognizing the Hopf fibration}

Given a great circle fibration, if its Reznikov invariant \(s^p_q\) vanishes with its first derivative at a point, then the structure equations show that \(s^0_r\) also vanishes there.
Consider what if the Reznikov invariant vanishes everywhere.
The equations \(s^p_q=s^0_r=0\) for all \(p,q,r\) form an involutive system, with general solution depending on functions of \(n\) variables.
In particular, there are real analytic elliptic immersed submanifolds \(X \subset \Gro{2}{V}\) which have vanishing Reznikov invariant.
Nothing is known about them.
Previous versions of this paper wrongly stated that vanishing of the Reznikov invariant everywhere forces the great circle fibration to be a Hopf fibration.

%{
%\begin{theorem}
%Given a great circle fibration \(S^{2n+1} \to X\)
%with \(n>1\) (i.e. not \(S^3 \to X^2\)),
%the \(s^p_q\) invariant vanishes
%precisely when the fibration is a Hopf
%fibration.
%\end{theorem}
%\begin{proof}
%If the invariant \(s^p_q\) vanishes,
%then the structure equations determine
%that
%\[
%0=s^p_{qr} + \delta^p_q s^0_r 
%= s^p_{\bar{r}q} = s^0_{\bar{q}t} \delta^p_u - s^0_{\bar{q}u} \delta^p_t.
%\]
%So if \(n > 1\) (i.e. our sphere \(S^{2n+1}\)
%has dimension at least 5) this last equation says
%that
%\(
%s^0_{\bar{p}q} = 0.
%\)
%Differentiating our structure equations,
%we find that finally this forces
%all invariants to vanish, so that 
%by the Frobenius theorem the result follows.
%\end{proof}
%
%The \(s^p_q\) invariant for \(S^3 \to X^2\)
%great circle fibrations vanishes, although
%great circle fibrations of \(S^3\)
%are generically not Hopf fibrations.

\section{``Hyperplanes''}\label{sec:Hyperplanes}

Given a real valued 1-form \(\xi \in V^*\)
on \(V\) we can define a section of
the bundle \(\tilde{V}^* \to X\)
by setting 
\(\sigma_\xi(x)\) to be \(\xi\).
This determines a function \(f \colon B \to V^*\)
by
\[
f(g) = \xi g.
\]
This function satisfies
\[
0=d
\begin{pmatrix}
f_0 \\
f_{\bar{0}} \\
f_{p} \\
f_{\bar{p}}
\end{pmatrix}
-
\begin{pmatrix}
\Omega^0_0 & \Omega^{\bar{0}}_0 & \Omega^q_0 & \Omega^{\bar{q}}_0 \\
\Omega^0_{\bar{0}} & \Omega^{\bar{0}}_{\bar{0}} 
        & \Omega^q_{\bar{0}} & \Omega^{\bar{q}}_{\bar{0}} \\
\Omega^0_p & \Omega^{\bar{0}}_p & \Omega^q_p & \Omega^{\bar{q}}_p \\
\Omega^0_{\bar{p}} & \Omega^{\bar{0}}_{\bar{p}} & 
        \Omega^q_{\bar{p}} & \Omega^{\bar{q}}_{\bar{p}} \\
\end{pmatrix}
\begin{pmatrix}
f_0 \\
f_{\bar{0}} \\
f_q \\
f_{\bar{q}}
\end{pmatrix}.
\]
Because \(\xi\) is real,
\[
f_{\bar{0}} = \bar{f}_0 \text{ and }
f_{\bar{p}} = \bar{f}_p.
\]

The bundle \(\OO{1} \to X\) is a quotient bundle
of \(\tilde{V}\),
and these sections determine sections of the
quotient which are just the functions \(f_0\).
They satisfy
\[
df_0 - \Omega^0_0 f_0 = f_{\bar{0}} \Omega^{\bar{0}}_0
+ f_p \Omega^p_0 + f_{\bar{p}} \Omega^{\bar{p}}_0.
\]
In particular, they are not holomorphic 
unless certain \(s\) invariants vanish.

Suppose that \(\xi \ne 0\).
Consider the locus \((f_0=0)\) inside \(X\).
(Of course, we should write it as something
like \(\left(\sigma_0=0\right)\) since 
\(f_0\) is really a function on \(B\), not
on \(X\).) At each point of this set,
\[
df_0 = f_q \Omega^q + f_{\bar{q}} s^{\bar{q}}_{\bar{r}} \Omega^{\bar{r}}.
\]
This can't vanish, since \(f_0=0\) forces
some \(f_q \ne 0\). Therefore the 
``hyperplanes'' determined by 
vanishing of these sections of \(\OO{1}\)
are smooth submanifolds of \(X\).
They are the analogues of complex hyperplanes.
We will see that they are connected.

\section{The Sato map}

Sato \cite{Sato:1984} constructed a map from the base \(X\) of a great circle fibration into a complex projective space. Yang \cite{Yang:1990} pointed out that this map was not well defined.
We will now present a very similar (but obviously well defined) map, which we will call the \emph{Sato map}.

To motivate it, consider any complex structure \(J\) on our vector space \(V\). As we have remarked above, this \(J\) by definition is a linear map \(J \colon V \to V\) satisfying \(J^2 = -1\), and therefore \(J\) has two eigenspaces
\(V^{1,0},V^{0,1} \subset V_{\C{}}=V \otimes_{\R{}} \C{}\), where \(V^{1,0}\) is the eigenspace with eigenvalue \(\sqrt{-1}\), and \(V^{0,1}\) is the eigenspace with eigenvalue \(-\sqrt{-1}\).
We thereby obtain two complex projective subspaces 
\[
\CP{}\left(V^{1,0}\right),
\CP{}\left(V^{0,1}\right)
\subset \SatoSpace.
\]
We can identify \(V\) with \(V^{1,0}\)
by 
\[
v \in V \mapsto v-\sqrt{-1}Jv \in V^{1,0}
\]
and projectivize this map to identify
\[
\CP{}\left(V\right) = \CP{}\left(V^{1,0}\right)
\]
(where the left hand side means the
space of \(J\) complex lines in \(V\)).
The copies of \(\CP{n}\) inside \(\SatoSpace\)
that occur this way are precisely those
with no real points, not intersecting \(\RP{}(V)\).
So the generic \(\CP{n} \subset \SatoSpace\) 
occurs in this way.
We will now imitate this story in the
context of a general great circle
fibration.

Consider a great circle fibration \(S^{2n+1} \to X\).
Take any \(x \in X\). This \(x\) may
be identified with an oriented
2-plane \(x \subset V\) since
\(X \subset \Gro{2}{V}\). 
There is a complex structure \(J_x\) on \(V\) for which \(x\) is a complex line, defined in section \vref{sec:Osculating}.
Map \(x \to V_{\C{}}\) by
\[
\sigma \colon v \in x \mapsto v-\sqrt{-1}J_x v \in V_{\C{}}.
\]
Since the fiber of \(\OO{-1} \to X\) above \(x\) is just \(x \subset V\) itself, the Sato map \(\sigma\) is defined on the total space of \(\OO{-1}\).

\begin{lemma}
The Sato map satisfies
\(
\sigma\left(J_x v\right) = \sqrt{-1} \sigma(v)
\)
for any \(v \in x \in X\).
\end{lemma}
\begin{proof}
\begin{align*}
\sigma\left(J_x v\right) &= J_x v - \sqrt{-1} J_x J_x v \\
  &= \sqrt{-1}\left(v - \sqrt{-1} J_x v\right).
\end{align*}
\end{proof}
As a consequence, \(\sigma\) takes the \(J_x\) complex line \(x\) to a complex
line \(\sigma(x)\) in \(V_{\C{}}\). 
Define the Sato map
\(
\sigma \colon X \to \SatoSpace
\)
to be the map that assigns to \(x\) the
line \(\sigma(x)\). 
This gives a morphism of bundles
\[
\begin{tikzcd}
\OO{-1}_X \arrow{d} \arrow{r}{\sigma} & \OO{-1}_{\SatoSpace} \arrow{d} \\
X \arrow{r}{\sigma} & \SatoSpace
\end{tikzcd}
\]
which is complex linear on the fibers and \(\SL{V}\)-equivariant.
The action of \(\SL{V}\) on \(\SatoSpace\) is by biholomorphisms (indeed, by projective automorphisms).

\begin{lemma} The Sato map \(X \to \SatoSpace\) is injective.
\end{lemma}
\begin{proof}
Suppose that \(\sigma(x)=\sigma(y)\).
Then for any \(v \in x\) and \(w \in y\)
we must find \(\sigma(v)\) and 
\(\sigma(w)\) on the same complex
line in \(V_{\C{}}\). 
Therefore, 
\[
\sigma(v)=v-\sqrt{-1}J_x v = 
\left ( a + b \sqrt{-1} \right ) \left (w - \sqrt{-1} J_y w \right)
=
\left(a+bJ_y\right) w + \sqrt{-1} \left(b - a J_y \right ) w
\]
which implies that 
\[
v = \left(a+b J_y\right)w \in y.
\]
But the fiber above \(y\) is \(J_y\) invariant,
so \(v\) and \(w\) lie on the same great
circle fiber, and therefore \(x=y\).
\end{proof}

Let us look at 
the locus of points \(z \in V_{\C{}}\)
so that \(z \wedge \bar{z} = 0 \in \Lm{1,1}{V_{\C{}}}\).
This projects to a variety \(\Delta \subset \SatoSpace.\)
We can parameterize this variety as
follows: a vector \(z \in V_{\C{}}\) satisfies
\(z \wedge \bar{z}\) precisely
when it has the form \(z=v+\sqrt{-1}w\)
with \(v\) and \(w\) belonging to the
same real line in \(V\). Multiplying
by a suitable complex number, we 
can get \(z=v \in V\). So we can identify
each element of \(\Delta\) with a
vector in \(V\) up to real multiples:
\(
\Delta = \RP{}\left(V\right)
\)
i.e. \(\Delta\) consists in the real
points of \(\SatoSpace.\)

\begin{lemma}
The space \(\Sato = \SatoSpace \backslash \RP{}(V)\) is canonically
identified with the space of pairs
\((P,j)\) where \(P \subset V\) is a 2-plane,
and \(j \colon P \to P\) is a complex structure.
\end{lemma}
\begin{proof}
Membership in \(\RP{}(V)\) is precisely
given by the equation \(z \wedge \bar{z} = 0\)
and therefore the solutions in \(V_{\C{}}\)
consist precisely in vectors \(z=v+\sqrt{-1}w\)
so that \(v \wedge w=0\), i.e. \(v,w \in V\)
belong to the same real line. Therefore
away from \(\RP{}(V),\)
points \(\ray{v+\sqrt{-1}w} \in \SatoSpace \backslash \RP{}(V)\)
map to 2-planes \(P = \ray{v \wedge w}\).
To a point 
\(\ray{v + \sqrt{-1}w} \in \Sato = \SatoSpace \backslash \RP{}(V)\),
attach the complex structure \(j\) which
maps
\(
jv = -w, jw=v.
\)
On the other hand, given an oriented 2-plane \(P\) with 
complex structure \(j\)
and a point \(v \in P \backslash 0,\)
map it to \(\C{\times}\left(v-\sqrt{-1}jv\right)\).
This gives continuous maps
in each direction between \(\Sato\)
and the space of 2-planes with complex structures,
and these maps are easily seen to be
inverses of one another.
Since the space of 2-planes with
complex structures is a homogeneous space
of \(\SL{V}\), and the maps are
\(\SL{V}\)-invariant, both
sides are homogeneous spaces and
these maps are equivariant diffeomorphisms.
\end{proof}

Each fiber of \(\Sato \to \Gro{2}{V}\) is
thus the set of complex structures on a given
oriented 2-plane \(P\), forming a copy of the
hyperbolic plane. We can look at the
same 2-plane \(P\) with the opposite
orientation, and find another such
hyperbolic plane. These two planes
are glued together inside \(\SatoSpace\)
along the real locus of points of the
projectivization of \(P\).
In other
words, the points of these two hyperbolic
planes are the complex points of the
projectivization of \(P\). Thus the fibers
of \(\Sato \to \Gro{2}{V}\) are open
subsets (``halves'') of the complex
points of real projective lines in
\(\SatoSpace.\)
}

We can write out the structure equations for \(G \to \Sato\): the 1-forms
\(\Omega^p_0, \Omega^{\bar{p}}_0,\Omega^{\bar{0}}_0\)
are semibasic (because they are independent on \(G\) but vanish
on the structure group) and they satisfy
\begin{equation}\label{eqn:ComplexStruc}
d
\begin{pmatrix}
\Omega^p_0 \\
\Omega^{\bar{p}}_0 \\
\Omega^{\bar{0}}_0 
\end{pmatrix}
=
-
\begin{pmatrix}
\Omega^p_q - \delta^p_q \Omega^0_0 &
        \Omega^p_{\bar{q}} &
        0 \\
\Omega^{\bar{p}}_q & 
        \Omega^{\bar{p}}_{\bar{q}} - \delta^{\bar{p}}_{\bar{q}} \Omega^0_0 &
        0 \\
\Omega^{\bar{0}}_q &
        \Omega^{\bar{0}}_{\bar{q}} &
        \Omega^{\bar{0}}_{\bar{0}} - \Omega^0_0         
\end{pmatrix}
\wedge
\begin{pmatrix}
\Omega^q_0 \\
\Omega^{\bar{q}}_0 \\
\Omega^{\bar{0}}_0 
\end{pmatrix}
+
\underbrace{
\Omega^{\bar{0}}_0 \wedge
\begin{pmatrix}
\Omega^p_{\bar{0}} \\
\Omega^{\bar{p}}_{\bar{0}} \\
0 
\end{pmatrix}}_{\text{torsion}}.
\end{equation}
From this expression, we find that 
the 1-forms \(\Omega^p_0,\Omega^{\bar{p}}_0,
\Omega^{\bar{0}}_0\) vary in a complex 
linear representation under the action of
the structure group.
Declaring them to be  \((1,0)\)-forms
for an almost complex structure,
the torsion consists
of \((1,1)\) forms, so the almost complex structure is
a complex structure: \(\Sato\) is a complex manifold.

\begin{lemma} The complex structure on \(\Sato\)
given by the structure equations~\vref{eqn:ComplexStruc}
is the same as that given by the 
embedding \(\Sato \subset \SatoSpace\).
\end{lemma}
\begin{proof}
The generic \(\CP{n} \subset \SatoSpace\)
lies in \(\Sato\) and occurs as the
base of a Hopf fibration. From the structure
equations of a Hopf fibration we see that
these are complex submanifolds in that
complex structure. But they are obviously
complex submanifolds in the \(\SatoSpace\)
complex structure. The complex structures
on these \(\CP{n}\) submanifolds agree,
as we can see directly looking at the
construction of the Hopf fibration. Since
the tangent planes to these \(\CP{n}\)
form an open subset of \(n\) dimensional 
complex tangent planes, they
force both complex structures on \(\Sato\)
to agree.
\end{proof}

The two zeroes in these structure equations,
in the expressions for \(d \Omega^p_0\) and
\(d \Omega^{\bar{p}}_0\), show that this
manifold is fibered by complex curves,
defined by the equations
\(
\Omega^p_0 = \Omega^{\bar{p}}_0 = 0
\)
and these are of course the fibers
of 
\(
\Sato \to \Gro{2}{V}
\)
since the 1-forms \(\Omega^p_0, \Omega^{\bar{p}}_0\)
are semibasic for this map.
Note \(\Sato \to \Gro{2}{V}\) is not a holomorphic map for any complex structure on \(\Gro{2}{V}\), because the relevant representation of the structure group does not preserve a complex structure. 
The torsion terms up on \(\Sato\) reorganize into the soldering 1-forms on \(\Gro{2}{V}\).

\begin{corollary}
The Sato map \(\sigma \colon X \to \Sato \subset \SatoSpace\)
is an embedding. The complex 
line bundle \(\OO{-1}_X\) is the pullback
of the complex line bundle \(\OO{-1}_{\SatoSpace}\)
by \(\sigma\). The image of the Sato
map determines the great circle fibration
\(S^{2n+1} \to X\).
\end{corollary}
\begin{proof} The map
\(X \to \Gro{2}{V}\) is an embedding, and the
map \(\Sato \to \Gro{2}{V}\)
is an \(\SL{V}\)-equivariant fiber bundle.
\end{proof}

\begin{corollary} The ``hyperplanes''
in \(X\)
which were discussed in section~\vref{sec:Hyperplanes}
are the intersections of \(X\)
with complex hyperplanes in \(\SatoSpace\).
\end{corollary}
\begin{proof} This is just saying that
since \(\OO{-1}_X\) is the pullback
of \(\OO{-1}_{\SatoSpace}\), therefore
the dual of the pullback, i.e. \(\OO{1}_X,\)
is the pullback of the dual.
\end{proof}
Consider the requirement that
the 1-forms used to cut out the hyperplanes
are real valued. Take
a real valued 1-form \(\xi \in V^*\) and
split it into complex linear and conjugate
linear parts on \(V_{\C{}}\) by 
\[
\xi^{1,0}\left(v+\sqrt{-1}w\right) = \xi(v) + \sqrt{-1}\xi(w)
\]
and 
\[
\xi^{0,1}\left(v+\sqrt{-1}w\right) = \xi(v) - \sqrt{-1}\xi(w).
\]
Then the \(\xi^{1,0}\) part cuts out a hyperplane of \(v+\sqrt{-1}w\)
so that \(\xi(v)+ \sqrt{-1}\xi(w) = 0\), i.e. where
both \(\xi(v)\) and \(\xi(w)\) vanish.

The Sato map
is given by taking elements of 
\(V^*_{\C{}}\) and using them as
sections of \(\OO{-1}_X\) and then
using that to make an embedding
\[
P \mapsto \left[\sigma_0(P):\dots:\sigma_{2n}(P)\right].
\]

\section{Mapping to projective space}\label{section:hinge}

Two great circle fibrations \(X, X' \subset \Gro{2}{V}\) are \emph{parallel} at a point \(x \in X \cap X'\) if they induce the same complex structure on the 2-plane \(P=P_x \subset V\) corresponding to the point \(x \in \Gro{2}{V}\).
If a point \(x \in X\) is associated to a 2-plane \(P \subset V\), then the differential of the map \(X \to \Gro{2}{V}\) is
\[
T_x X \to T_P \Gro{2}{V} = P^* \otimes \pr{V/P},
\]
and becomes complex linear for unique complex structures on \(T_x X\), \(P\) and \(V/P\), determined algebraically from the linear map
\[
T_x X \to T_P \Gro{2}{V},
\]
i.e. from the tangent spaces of \(X\).
Any complex structure on \(P\) arises this way, for a suitable choice of tangent space by \(\SL{V}\)-equivariance. 
Hence if we can slightly perturb the tangent planes of two great circle fibrations at their intersection points, then they will be nowhere parallel.
Another way to phrase parallelism: points of parallelism are intersection points in the image of the Sato map inside \(\SatoSpace\).

Parallelism of great circle fibrations \(X_0, X_1\) at a point \(x\) clearly implies that the osculating complex structures \(J_{1,x}, J_{2,x}\) agree on the 2-plane \(P=P_x\), and so 
\[
0 = \det\pr{J_{1,x}-J_{2,x}}.
\]
In particular, if we can find some constant complex structure \(J \in V^* \otimes V\) whose Hopf fibration base \(X_J=\CP{n}_J \subset \Gro{2}{V}\) satisfies
\[
0 \ne \det\pr{J-J_x}
\]  
for all \(J_x \in K\), i.e. for all osculating complex structures to \(X\)
then \(X_J\) and \(X\) are nowhere parallel.

Take a great circle fibration with base space \(X \subset \Gro{2}{V}\).
Let \(\CP{n}=\CP{n}_J\) be the Hopf fibration base of some complex structure \(J\). 
By transversality, after a \(C^2\)-small perturbation of \(X \subset \Gro{2}{V}\), we can arrange that \(X\) has finitely many intersections with \(\CP{n} \subset \Gro{2}{V}\).
We can also arrange that these intersection points \(x \in X \cap \CP{n}\) are transverse, and have ``generic'' tangent space \(T_x X\).
In particular, achieving generic tangent spaces \(T_x X\) at intersection points \(x \in X \cap \CP{n}\) ensures that \(X\) and \(\CP{n}\) are nowhere parallel, i.e.  \(X\) and \(\CP{n}\) are disjoint inside \(\SatoSpace\).

Take two complex structures \(J\) and \(J'\) on \(V\), whose associated Hopf fibration bases are nowhere parallel to \(X\).
It is easy to check that small perturbation of \(J'\) among complex structures will ensure that \(J'\) and \(X\) are pairwise nowhere parallel and that \(J-J'\) is invertible.
Define a map
\[
f=f_{J,J',X} \colon v \in V\backslash 0 \mapsto \pr{J'-J}^{-1} \pr{J_x - J'}v \in V\backslash 0,
\]
where \(v \in S^{2n+1} \mapsto x \in X\) is the bundle map of our great circle fibration, with osculating complex structure \(J_x\).

The reader can easily check the following linear algebra fact: for any 3 complex structures \(J_0,J, J'\) on the same vector space \(V\), with \(J'-J\) invertible, and \(J_0\) and \(J'\) nowhere parallel, the linear equations
\begin{align*}
0 &=  v_0 + v + v', \\
0 &= J_0 v_0 + Jv + J'v', 
\end{align*}
determine any two of \(v,v'\) as real linear functions of \(v_0\).
From the form of these equations, if we replace \(v_0\) by \(J_0 v_0\), and so on, the solutions change to
\begin{align*}
0 &= J_0v_0 + Jv + J'v', \\
0 &= J_0J_0v_0 + JJv +J'J'v', 
\end{align*}
i.e. the solution for a given \(v_0\) gives a complex linear map
\[
v_0 \in \pr{V,J_0} \mapsto \pr{v,v'} \in \pr{V,J} \times \pr{V,J'}.
\]
Moreover, this solution is given by
\begin{align*}
v &= \pr{J'-J}^{-1}\pr{J_0 - J'} v_0, \\
v' &= -v_0 - v.
\end{align*}
In other words, if \(X=X_{J_0}\) is a Hopf fibration base, then the linear algebra tells us that the associated map \(v=f\pr{v_0}\) is complex linear on each fiber of the Hopf fibration.
Indeed in this simple case, \(f\) is a real linear automorphism of \(V\).

\begin{lemma}
Take any great circle fibration \(X \subset \Gro{2}{V}\).
Then there are complex structures \(J, J'\) for which \(J-J'\) is invertible and for which \(X\) is nowhere parallel to the Hopf fibration base \(X_{J'}=\CP{n}_{J'}\) and for which the associated smooth map \(f\) above takes \(V\backslash 0\) to \(V \backslash 0\), taking each 2-plane \(x=\left<v,J_x v\right> \subset V\) associated to a point \(x \in X\) to a \(J\)-linear 2-plane \(\left<w,Jw\right>\) where \(w=f(v)\).
The map \(f\) is linear on these planes.
If \(f\) is a local diffeomorphism of \(V \backslash 0\) then it descends to a diffeomorphism \(X \to X_{J'}=\CP{n}_{J'}\) to a Hopf fibration, linear on the fibers of the Hopf fibration.
\end{lemma}
\begin{proof}
Take a vector \(v \in V\backslash 0\) and let \(x=\left<v,J_x v\right> \in X\) be the associated oriented 2-plane.
If \(f(v)=0\) then one easily checks that \(J_x v = Jv\) and so \(X\) and \(X_J\) are parallel at \(x\).
Consequently \(f \colon V\backslash 0 \to V\backslash 0\) is a smooth map, equivariant under rescaling.

Applying the linear algebra of the osculating complex structures to \(X\) to the equations above, we find that \(f\) is complex linear on each 2-plane \(x \in X\).
Once we see that \(f\) is a local diffeomorphism, rescaling invariance allows us to treat \(f\) as a smooth local diffeomorphism of the sphere, so a diffeomorphism.
\end{proof}

In order to see whether \(f\) is a local diffeomorphism, we need to see whether, for a suitable choice of \(J'\), \(f'\pr{v_0}w_0 \ne 0\) for all \(w_0\).
For each vector \(v_0 \ne 0\), we let \(J_{v_0}\) be \(J_x\) for the associated \(x \in X\), and we let \(A_{v_0} \in V^* \otimes V\) be the linear map
\[
A_{v_0} w = J_{v_0}w + \left.\frac{d}{dt}\right|_{t=0} J_{v_0+tw} v_0.
\]
Calculate
\[
f'\pr{v_0}w
= 
\pr{A_{v_0}-J'}w.
\]
Consider the map
\[
v_0 \in S^{2n+1} \to A_{v_0} \in V^* \otimes V.
\]
If we can choose \(J'\) so that additionally 
\[
0 \ne \det{J'-A_{v_0}},
\]
for any \(v_0 \in V\backslash 0\), i.e. so that \(f'\pr{v_0}\) has no kernel then \(f\) is a diffeomorphism.

Consider again the Hopf fibration, and its base manifold \(X_J\). 
The base sits inside \(\SatoSpace\) as a linear subspace, cut out by the equation
\(
Jz=\sqrt{-1}z.
\)
So it is a linear projective subspace,
\[
X_J = \CP{n} \subset \CP{2n+1}=\SatoSpace.
\]
If we pick another linear \(\CP{n}\) subspace, call it \(\CP{n}_0\), which is not parallel to \(X_J\), we can look at the family of all linear \(\CP{n+1}\) subspaces containing \(\CP{n}_0\).
These \(\CP{n+1}\)'s will each intersect \(\CP{n}\) in a unique point transversely.
%\begin{figure}
%\begin{picture}(0,0)%
%\drawfilePSTEX{Crossing}%
%%\includegraphics{Crossing.pdf}%
%%\epsfig{file=Crossing.pstex}%
%\end{picture}%
%\setlength{\unitlength}{3947sp}%
%%
%\begingroup\makeatletter\ifx\SetFigFont\undefined%
%\gdef\SetFigFont#1#2#3#4#5{%
%  \reset@font\fontsize{#1}{#2pt}%
%  \fontfamily{#3}\fontseries{#4}\fontshape{#5}%
%  \selectfont}%
%\fi\endgroup%
%\begin{picture}(3237,5412)(514,-4786)
%\put(3451,-4786){\makebox(0,0)[lb]{\smash{\SetFigFont{12}{14.4}{\familydefault}{\mddefault}{\updefault}\(X_J=\CP{n}\)}}}
%\put(3751,-3511){\makebox(0,0)[lb]{\smash{\SetFigFont{12}{14.4}{\familydefault}{\mddefault}{\updefault}\(\CP{n}_0\)}}}
%\put(1726,-961){\makebox(0,0)[lb]{\smash{\SetFigFont{12}{14.4}{\familydefault}{\mddefault}{\updefault}\(\CP{n+1}\)}}}
%\end{picture}
%\caption{Slicing with \(\CP{n+1}\)'s
%containing a given hinge \(\CP{n}\)}\label{fig:SatoOne} 
%\end{figure}
We will call the \(\CP{n}_0\) a \emph{hinge}.
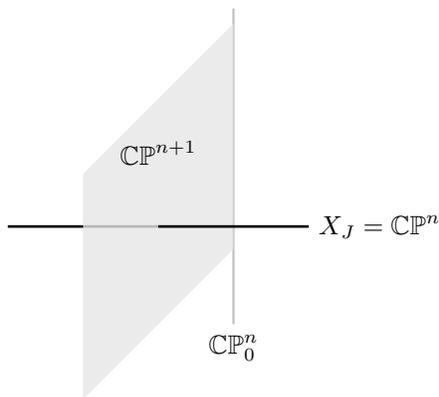
\begin{center}
\begin{tikzpicture}
\draw[gray!50,thick] (1,3.2) -- (1,-1) node[below,black] {\(\mathbb{CP}^n_0\)};
\fill[opacity=.8,gray!20] (0,-1) -- (1,0) -- (1,3) -- (0,2) -- cycle;
\draw[line width=1pt] (-2,.3) -- (2,.3) node[right] {\(X_J=\mathbb{CP}^n\)};
\fill[gray!20,opacity=.8] (-1,-2) -- (0,-1) -- (0,2) -- (-1,1) -- cycle;
\node[above] at (0,1) {\(\mathbb{CP}^{n+1}\)};
\end{tikzpicture}
\end{center}

On \(\Sato,\) the 1-forms 
\(\Omega^p_0,\Omega^{\bar{p}}_0,\Omega^{\bar{0}}_0\)
form a complex linear coframing,
and if this coframing comes from a point of
the adapted bundle \(B \to X_J\), then the
tangent space of \(X_J\) is given by
\(\Omega^{\bar{p}}_0=\Omega^{\bar{0}}_0=0\).
The \(\CP{n+1}\)'s striking this \(X_J=\CP{n}\) transversely at this point are given by a set of complex linear equations
\[
\Omega^p_0 = a^p_{\bar{q}} \Omega^{\bar{q}}_0 
+ a^p_{\bar{0}} \Omega^{\bar{0}}_0.
\]

The intersection of a \(\CP{n+1}\) with the submanifold \(X_J\) inside \(\SatoSpace\) is described in our bundle \(B\) by the equations
\[
0 = \Omega^p_0=\Omega^{\bar{p}}_0=\Omega^{\bar{0}}_0
\]
cutting out the \(\CP{n+1}\) together
with the equations
\[
0 = \Omega^p_{\bar{q}} =
\Omega^{\bar{p}}_0 =
\Omega^{\bar{0}}_p =
\Omega^{\bar{0}}_0
\]
cutting out the \(X_J=\CP{n}\).

Repeat this story using a general great circle fibration \(S^{2n+1} \to X\);
%as in figure~\vref{fig:SatoTwo},
a tangent space to \(X \subset \Sato\) is given in an adapted coframe by
\begin{align*}
\Omega^{\bar{p}}_0&=s^{\bar{p}}_{\bar{q}} \Omega^{\bar{q}}_{\bar{0}} \\
\Omega^{\bar{0}}_0&=s^{\bar{0}}_{\bar{q}} \Omega^{\bar{q}}_{\bar{0}}
\end{align*}
while the osculating Hopf fibration \(X_J\) (for \(J=J_P\), osculating at point \(P \in X\)) satisfies similar equations \(0=\Omega^{\bar{p}}_0=\Omega^{\bar{0}}_0\).
It is easy to calculate that the common tangent vectors in these tangent spaces are precisely those given by \(v \hook \Omega^p_0 = v^p\) where \(v^p\) is any eigenvector with eigenvalue 1 of the matrix
\[
\pr{
a^p_{\bar{t}}
s^{\bar{t}}_{\bar{r}}
+
a^p_{\bar{0}}
s^{\bar{0}}_{\bar{r}}
}
\pr{
a^{\bar{r}}_s
s^s_q
+
a^{\bar{r}}_0
s^0_q
}.
\]

If the spectrum of this matrix lies inside the unit disk at some point of the bundle \(B\), then the same is true for any point in the same fiber over \(X\), and we say that the \(\CP{n+1}\) is \emph{strongly transverse} to \(X\) at such a point.

%Therefore the \(\CP{n+1}\)'s transversal to \(X_J\) are precisely the same as those transversal to \(X\). 
%The intersections with \(X\) are positive.

Returning to the Hopf fibration \(X_J=\CP{n} \subset \CP{n+1}\), take any other \(\CP{n}_0 \subset \CP{n+1}\) so that \(\CP{n}_0\) is not parallel to \(X_J\), 
i.e. so that \(\CP{n}_0 \cup X_J\) is not contained in any \(\CP{n+1}\). 
Call such a \(\CP{n}_0\) a \emph{hinge}. 
Take all \(\CP{n+1}\) in \(\SatoSpace\) containing \(\CP{n}_0\). 
If each of these strikes \(X_J\) at a single point transversely, we can identify \(X_J\) with the projective quotient \(\SatoSpace/\CP{n}_0=\CP{n}\),
i.e. with the projectivized holomorphic normal bundle at any point of  \(\CP{n}_0\). 

The same procedure works for a general great circle fibration \(S^{2n+1} \to X\), under some additional hypotheses. 
A linear subspace  \(\CP{n} \subset \SatoSpace=\CP{2n+1}\) is \emph{parallel} to \(X\) at a point \(P \in X\) if it is parallel to the osculating Hopf fibration \(X_J\) for \(J=J_P\).
Conversely, a linear subspace \(\CP{n} \subset \SatoSpace\) is called a \emph{hinge} for \(X\) if first it is nowhere parallel to \(X\) and second it does not intersect \(X\) and third every \(\CP{n+1}\) containing that \(\CP{n}\) intersects \(X\) transversely and positively

%\begin{lemma}\label{lemma:NotParallel}
%For any great circle fibration \(S^{2n+1} \to X\), sitting \(X\)
%inside \(\SatoSpace\) by the Sato map, there is a hinge
%for \(X\). 
%The hinges for \(X\) form a dense open subset in the space of linear \(\CP{n}\) subspaces of \(\SatoSpace\), i.e. in
%\(
%\GrC{n+1}{2n+2}.
%\)
%\end{lemma}
%\begin{proof}
%To see that a hinge for \(X\)
%exists, we have only to examine
%dimensions. This \(X\) has \(2n\) real 
%dimensions, so generically it will have \(2n\)
%real dimensions of osculating 
%Hopf fibrations \(X_J\), for \(J=J_P\)
%for \(P \in X\). Let \(W\) be the
%set of triples \((P,A,B)\)
%where \(P \in X\), \(A\) is a linear
%\(\CP{n+1}\) subspace of \(\CP{2n+1}\)
%containing the osculating Hopf
%fibration \(X_{J_P}=\CP{n}\)
%to \(X\) at \(P\) and \(B\) is a linear 
%\(\CP{n}\) subspace contained in \(A\).
%We can see that at a generic point of 
%\(W\), \((P,B)\) determines \(A\) because
%\(A\) is the linear 
%subspace containing \(B\) and \(X_{J_P}\).
%The possible subspaces
%\(B\) are parameterized by the
%normal bundle to \(X_{J_P}\) at \(P\).
%Therefore the dimension
%of \(W\) as a real manifold is
%\(
%2n+2(n+1)=4n+2.
%\)
%So this is the dimension of the
%family of \(\CP{n}\) parallel to \(X\).
%
%
%On the other hand, the family
%of all \(\CP{n}\) inside \(\CP{2n+1}\)
%has real dimension \(2n^2+4n+2\)
%which is always bigger than \(4n+2\).
%Therefore by Sard's theorem,
%there is some \(\CP{n}_0\) which
%is nowhere parallel to \(X\),
%and we can pick it from an
%open set of \(\CP{n}\) linear
%subspaces, and moreover arrange
%that it doesn't intersect \(X\).
%\end{proof}

\begin{lemma}\label{lemma:Socks}
Take \(S^{2n+1} \to X\) a great circle fibration, and stick \(X\) into \(\SatoSpace\) via the Sato map. 
Take a hinge \(\CP{n}_0\) for \(X\).
The family of \(\CP{n+1}\) linear subspaces
containing that hinge is diffeomorphic
to \(\CP{n}\).
There is a diffeomorphism \(X \to \CP{n}\) given
by taking a point \(P \in X\) to the
unique \(\CP{n+1}\) linear space containing
\(\CP{n}_0\) which strikes \(X\) at \(P\).
\end{lemma}
\begin{proof}
Consider the manifold \(Y \cong \CP{n}\) 
of all \(\CP{n+1}\) containing the hinge \(\CP{n}_0\),
and the incidence correspondence \(Z\) of all
pairs \((P,\Pi)\) with \(\Pi \in Y\)
and \(P \in X\) and \(X_{J_P} \subset \Pi\).
By positivity of intersections, \(Z\) 
is a submanifold of \(X \times Y\)
of the same dimension (\(2n\)) as
both \(X\) and \(Y\), and \(X\)
and \(Y\) are compact, and so \(Z\) is too. Again by 
positivity of intersections, 
\(Z \to Y\) and \(W \to X\) are covering
maps, preserving orientation. 
But \(X\) and \(Y = \CP{n}\) are both simply connected,
so \(Z \to Y\) and \(Z \to X\) are diffeomorphisms.
\end{proof}

\begin{lemma}
If a great circle fibration admits a hinge then it admits one of the form \(\CP{n}_0=X_J\) for some complex structure \(J\) on \(V\).
\end{lemma}
\begin{proof}
By dimension count, we see that these \(X_J \subset \SatoSpace\) (which are in one to one correspondence with complex structures \(J\) on \(V\)) have
the same dimension as the space of all \(\CP{n}\) linear subspaces in \(\SatoSpace\), so they form an open subset. 
Indeed the \(X_J\) subspaces are precisely those \(\CP{n}\) linear subspaces with no real points on them.
\end{proof}

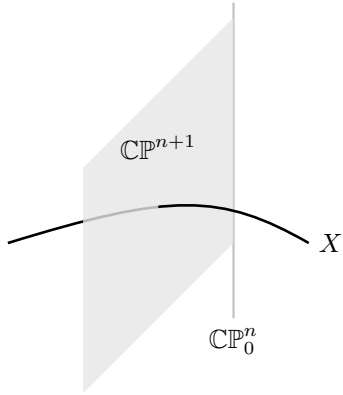
\begin{center}
\begin{tikzpicture}
\draw[gray!50,thick] (1,3.2) -- (1,-1) node[below,black] {\(\mathbb{CP}^n_0\)};
\fill[opacity=.8,gray!20] (0,-1) -- (1,0) -- (1,3) -- (0,2) -- cycle;
\draw[line width=1pt] (-2,0) .. controls (0.17,0.67) and (0.83,0.67) .. (2,0) node[right] {\(X\)};
\fill[gray!20,opacity=.8] (-1,-2) -- (0,-1) -- (0,2) -- (-1,1) -- cycle;
\node[above] at (0,1) {\(\mathbb{CP}^{n+1}\)};
\end{tikzpicture}
\end{center}

%\begin{figure}
%\begin{picture}(0,0)%
%\drawfilePSTEX{Crossing2}%
%%\includegraphics{Crossing2.pdf}%
%%\epsfig{file=Crossing2.pstex}%
%\end{picture}%
%\setlength{\unitlength}{3947sp}%
%%
%\begingroup\makeatletter\ifx\SetFigFont\undefined%
%\gdef\SetFigFont#1#2#3#4#5{%
%  \reset@font\fontsize{#1}{#2pt}%
%  \fontfamily{#3}\fontseries{#4}\fontshape{#5}%
%  \selectfont}%
%\fi\endgroup%
%\begin{picture}(3162,5412)(589,-4786)
%\put(3751,-3511){\makebox(0,0)[lb]{\smash{\SetFigFont{12}{14.4}{\familydefault}{\mddefault}{\updefault}\(\CP{n}_0\)}}}
%\put(1726,-961){\makebox(0,0)[lb]{\smash{\SetFigFont{12}{14.4}{\familydefault}{\mddefault}{\updefault}\(\CP{n+1}\)}}}
%\put(3451,-4786){\makebox(0,0)[lb]{\smash{\SetFigFont{12}{14.4}{\familydefault}{\mddefault}{\updefault}\(X\)}}}
%\end{picture}
%\caption{Slicing with \(\CP{n+1}\)'s
%containing a given \(\CP{n}\)}\label{fig:SatoTwo}
%\end{figure}
%

Now that we can map the base of our great circle fibration to a complex projective space, we have to map the great circles to the
fibers of the Hopf fibration. 
Consider first how to do this for a pair Hopf fibrations, given by complex structures \(J_1\) and \(J_2\).
We take any other Hopf fibration, given by a complex structure \(J_0\), so that
if the associated fibrations are \(S^{2n+1} \to X_k\) for \(k=0,1,2\), then \(X_0 \subset \SatoSpace\) will not be parallel to \(X_1\) or to \(X_2,\)
i.e. \(X_0\) is a hinge for \(X_1\) and also for \(X_2\). 
Of course, these \(X_k\) are all \(\CP{n}\) linear subspaces inside \(\CP{2n+1} = \SatoSpace\) given by linear equations
\[
X_k = \left ( J_k = \sqrt{-1} \right ).
\]

Our map \(X_1 \to X_2\) is constructed by taking all linear \(\CP{n+1}\) subspaces containing \(X_0\) and matching the point \(\CP{n+1} \cap X_1\) to the point \(\CP{n+1} \cap X_2\).

First, we need to lift the entire picture up
to a linear picture in \(V_{\C{}}\). Then
we will see what is happening in \(V\) itself.
Up in \(V_{\C{}}\) we have 3 complex linear
subspaces \(V_0,V_1,V_2\), all isomorphic to \(\C{n+1}\),
which are just the preimages of the \(X_0,X_1,X_2\)
linear subspaces in \(\SatoSpace.\)
Another way to say this: 
\[
V_k = \left \{ v \, | \, J_k v = \sqrt{-1} v \right \}.
\]
Since the \(V_0\) has no vectors in common
with \(V_1\) or with \(V_2\), 
every vector in \(V_{\C{}}\) can be written
as a combination \(w=w_0+w_1\) of vectors
from \(V_0\) and \(V_1\). The vectors in \(V_k\) 
are of the form 
\[
v - \sqrt{-1} J_k v
\]
for vectors \(v \in V\).
The equation
\[
v_2 - \sqrt{-1} J_2 v_2 = v_0 - \sqrt{-1}J_0 v_0 + v_1 - \sqrt{-1} J_1 v_1
\]
for \(v_0,v_1,v_2 \in V\)
breaks into real and imaginary parts:
\begin{align*}
v_2 &= v_0 + v_1 \\
J_2 v_2 &= J_0 v_0 + J_1 v_1
\end{align*}
and has the solution 
\begin{align*}
v_0 &= \left(J_2 - J_0\right)^{-1}\left(J_1 - J_2\right)v_1 \\
v_2 &= v_0 + v_1 \\
\end{align*}
so taking \(v_1 \in V \to v_2 \in V\) by a linear
isomorphism, taking the \(\sqrt{-1}\) eigenspace
of \(J_1\) to that of \(J_2\). 
The equation
\(
J_2 v_2 = J_0 v_0 + J_1 v_1
\)
ensures us that this map takes \(J_1 v_1\) to \(J_2 v_2\).
So it is a complex linear map
\[
\left(V,J_1\right) \to \left(V,J_2\right).
\]
By construction it matches the 2-planes \(P_1\) and \(P_2\) in \(V\), just looking back at the construction of the Sato map.

As before, we can apply the same idea to a great circle fibration as follows:
we take \(S^{2n+1} \to X\) our great circle fibration, and \(S^{2n+1} \to X_2\)
a Hopf fibration, given by a complex structure \(J_2\). 
Now we pick another Hopf fibration \(S^{2n+1} \to X_0\) given by a complex structure \(J_0\), so that \(X_0\) is nowhere parallel to \(X\) and to \(X_2\).
For each 2-plane \(P \in X\) we take the osculating complex structure \(J_P\)
and use the above process to produce a complex linear map \(\left(V,J_P\right) \to \left(V,J_2\right).\)
This will map \(P\) to a 2-plane \(P_2\)
complex linearly:
\[
\left(V,P,J_P\right) \to \left(V,P_2,J_2\right).
\]
Hence it identifies the great circle fibrations.
We have proven:
\begin{theorem}
If a great circle fibration \(S^{2n+1} \to X\) admits a hinge, then the base manifold \(X\) is diffeomorphic to \(\CP{n}\) by a diffeomorphism which lifts to a diffeomorphism \(S^{2n+1} \to S^{2n+1}\) linear on the fibers. 
For each choice of hinge \(X_{J_0} \subset \SatoSpace\) for
\(X\) we obtain an isomorphism
\[
\begin{tikzcd}
S^{2n+1} \arrow{r} \arrow{d} & S^{2n+1} \arrow{d} \\
X \arrow{r} & \CP{n} 
\end{tikzcd}
\]
identifying the great circle fibration with a Hopf fibration. 
This isomorphism depends smoothly on the choice of hinge.
\end{theorem}

\begin{corollary}
Given any two great circle fibrations \(S^{2n+1} \to X_0\) and \(S^{2n+1} \to X_1\), both admitting a hinge,  the embedded submanifolds \(X_0, X_1 \subset \Gro{2}{V}\) are homotopic inside \(\Gro{2}{V}\).
\end{corollary}
\begin{proof} 
We pick a hinge for each \(X_j\) and an orientation preserving linear isomorphism identifying the hinges.
The linear isomorphism lies on a path inside the general linear group.
So we can assume that the two hinges are the same.
Deform \(X_0\) to \(X_1\) along the straight lines contained in each \(\CP{n+1} \backslash \CP{n}_0 = \C{n+1}\). 
Then take the image in \(\Gro{2}{V}\).
\end{proof}

\section{\texorpdfstring{Nonlinear $J$: twisted complex structures}{Nonlinear J: twisted complex structures}}
{\label{sec:NonlinearJ}

A great circle fibration \(S^{2n+1} \to X\) determines a 2-plane \(P\) through each point \(v \in V \backslash 0\), and a complex structure \(J_P \colon V \to V\)
for which \(P\) is a complex line.
Define the map \(J_X \colon V \backslash 0 \to V \backslash 0\) by \(J_X v = J_P v\). 
Clearly \(J_X\) is a smooth map, and satisfies \(J_X^2 = -1\). 
The map \(J_X\) is linear precisely if the great circle fibration is a Hopf
fibration.

A \emph{twisted complex structure} is a map \(J \colon V \backslash 0 \to V \backslash 0\) 
which 
\begin{enumerate}
\item
is continuous,
\item
satisfies \(J^2 = -1\)
\item
leaves invariant each 2-plane \(\text{span}\left<v,Jv\right>\)
and 
\item
is linear on those 2-planes.
\end{enumerate}
A real vector space \(V\) equipped with
a twisted complex structure is a \emph{twisted complex vector space}. 
(Twisted complex structures were discovered simultaneously by 
the author \cite{McKay:2001b} and by Jean-Claude Sikorav
\cite{Sikorav:2000}.)

Extend our twisted complex structures from \(V \backslash 0\) to \(V\)
by defining \(J0=0\), so that \(J \colon V \to V\) is a homeomorphism, but not differentiable at \(0\) unless \(J\) is linear, i.e. a complex structure.
A twisted complex structure is \emph{smooth} if \(J\) is smooth on \(V \backslash 0\).
Note that a twisted complex structure extends to \(0\) to be  smooth at \(0\) just when it is linear.

Given a smooth or topological twisted complex structure \(J\), define a smooth or topological great circle fibration
by taking \(X(J) \subset \Gro{2}{V}\) to be the set of oriented 2-planes
of the form \(\left<v,J_Xv\right>\) for \(v \in V\), oriented by setting \(v \wedge J_X v\) to be positive, taking the great circles in \(S^{2n+1} = (V \backslash 0)/\R{+}\) to be the quotients by \(\R{+}\) of the \(J\)-invariant 2-planes.

\begin{lemma}
For any smooth great circle fibration \(S^{2n+1} \to X\),
\(
X\left(J_X\right) = X.
\)
\end{lemma}
\begin{proof}
The 2-planes which are \(J_X\)-invariant are precisely those belonging to \(X\),
with the required orientation. 
So it is the same submanifold of the Grassmannian.
\end{proof}

A twisted complex structure \(J\) is a \emph{pseudocomplex structure}
if it satisfies \(J=J_X\) for \(S^{2n+1} \to X\) a great circle fibration.

\begin{lemma}
The generic smooth twisted complex structure \(J\) is \emph{not} pseudocomplex, i.e.
\(
J_{X(J)} \ne J.
\)
\end{lemma}
\begin{proof}
Let \(X=X(J)\) and \(J'=J_{X(J)}\).
Clearly \(J'\) leaves a 2-plane invariant precisely when \(J\) does.
These 2-planes have the same orientation, so they determine the same great circle fibration, \(X\). 
Each of \(J\) and \(J'\) also determine sections of \(\Sato \to \Gro{2}{V}\).
Recall that \(\Sato\) is the space of pairs \((P,j)\) so that \(P \subset V\) is a 2-plane and \(j \colon P \to P\) is a complex structure
on \(P\). 
Define
\(
\sigma_J \colon X \to \Sato
\)
by
\(
\sigma_J(P) = \left(P,J|_P\right)
\)
restricting \(J\) to \(P\).
But we can vary \(J\) to an arbitrary section of \(\Sato \to \Gro{2}{V}\)
over the same \(X\), while \(J'\) is fixed by the choice of \(X\).
\end{proof}

\begin{corollary}
The space of smooth twisted complex structures retracts to the space of great circle fibrations, i.e. to the space of pseudocomplex structures.
\end{corollary}
\begin{proof}
Indeed, using the notation of the preceding lemma, \(J=J'\) precisely when the sections \(X \to \Sato|_X\) agree.
Note that these are disk bundles, so the space of \(J\) with fixed \(X(J)\) is contractible. 
Or, more canonically, just take \(J \to X(J)\).
\end{proof}

If vector spaces \(V_0, V_1\) have twisted complex structures \(J_0, J_1\), then the \emph{sum} of these is \(V = V_0 \oplus V_1\)
with twisted complex structure 
\[
J\left(v_0,v_1\right) = \left(J_0 v_0, J_1 v_1\right).
\]

An obvious result:
\begin{proposition}
The sum of twisted complex vector spaces is twisted complex.
The sum is smooth just when one of the summands is linear, i.e. a complex structure.
\end{proposition}

It appears that Yang's proof of the topological Blaschke conjecture for Blaschke manifolds with the cohomology of complex projective spaces involves proving that every such Blaschke manifold is the base manifold of a sum of two great circle fibrations.

\section{The space $M$}

Consider the homogeneous \(\SL{V}\) space \(\M\) which consists of choices of pairs \((P,J)\) of 2-plane \(P \subset V\) and complex structure \(J:V \to V\)
so that \(P\) is a complex \(J\)-line.
We have described a map \(X \to \M\) from the base manifold \(X\) of any great circle fibration \(S^{2n+1}\to X\).
Indeed the structure group of our bundle \(B \to X\) is precisely the isotropy group of a point of \(\M\). 
There are obvious maps
\[
\begin{tikzcd}[ampersand replacement=\&]
\& \M \arrow{dr} \arrow{dl} \& \\
\Gro{2}{V} \& \& \Cstrucs{V}.
\end{tikzcd}
\]
The fibers of the map \(\M \to \Cstrucs{V}\) are copies of \(\CP{n}\). 
From the structure equations, we see that \(\M\) is a complex manifold.
In fact it is possible to see this from a different point of view: since the complex structures on \(V\) are identified with the complex linear subspaces \(W \subset V_{\C{}}\) with no real points, i.e. \(W \cap V = 0\), this is an
open subset of the Grassmannian \(\GrC{n+1}{2n+2}\). 
Above \(\GrC{n+1}{2n+2}\) we have the universal bundle \(\mathcal{U} \to \GrC{n+1}{2n+2}\) whose fiber above \(W \in \GrC{n+1}{2n+2}\)
is \(W\) itself. 
Projectivizing this bundle, we have the bundle 
\[
\mathbb{CP}(\mathcal{U}) \to \GrC{n+1}{2n+2}
\]
of complex projective spaces. 
A choice of \((P,J) \in \M\) is precisely a choice of complex subspace \(W \in \GrC{n+1}{2n+2}\) not containing any real vectors: \(W \cap V = 0\), and a choice of complex line inside \(W\): the line consisting of the vectors
\[
v - \sqrt{-1} J v
\]
for \(v \in P\). 
Consequently, \(\M\) is just the pullback
\[
\begin{tikzcd}[ampersand replacement=\&]
\M \arrow{d} \arrow{r} \& \mathbb{CP}(\mathcal{U}) \arrow{d} \\
\Cstrucs{V} \arrow{r} \& \GrC{n+1}{2n+2}.
\end{tikzcd}
\]
We can also identify \(\M\) with the homogeneous space
\[
\M = \SL{V}/\Gamma_0.
\]
Since a great circle fibration 
\[
\begin{tikzcd}[ampersand replacement=\&]
S^1 \arrow{r} \& S^{2n+1} \arrow{d} \\
           \& X^{2n}
\end{tikzcd}
\]
gives rise to a right principle \(\Gamma_0\) bundle
\[
\begin{tikzcd}[ampersand replacement=\&]
\Gamma_0 \arrow{r} \& B \arrow{d} \\
                \& X
\end{tikzcd}
\]
with \(B \subset \SL{V}\), we have
\[
\begin{tikzcd}[ampersand replacement=\&]
\Gamma_0 \arrow{r} \& B \arrow{d} \arrow{r} \& \SL{V} \arrow{d} \& \Gamma_0 \arrow{l} \\
                \& X \arrow{r}        \& \M.
\end{tikzcd}
\]
Consider also the map
\[
X \to \M \to \Cstrucs{V} \subset \GrC{n+1}{2n+2}.
\]
Picking a hinge, i.e. a choice of complex subspace \(W_0 \in \GrC{n+1}{2n+2} \cap \Cstrucs{V}\) which is transverse to every complex subspace \(W \in \GrC{n+1}{2n+2} \cap \Cstrucs{V}\) arising as osculating complex structure 
to \(X\), we can then trivialize the universal bundle over the part of \(\GrC{n+1}{2n+2}\) consisting of subspaces transverse to \(W_0\).
We do that by taking complex coordinates \(z,w\) on \(V_{\C{}}\) so that the complex \(n+1\) plane \(z=0\) is the hinge. 
Then all of the other complex \(n+1\) planes have the form \(w=pz\) so that \(p = \left(p^P_Q\right)\) is our Pl{\"u}cker coordinate system. 
Given two \(n+1\) planes, say with Pl{\"u}cker coordinates \(p_0\) and \(p_1\),
use the map
\[
\left(z,p_0z\right) \mapsto \left(z,p_1z\right)
\]
to identify them. 
This is a complex linear map, so it identifies complex lines with complex lines. 

On \(X \subset M\), map
\[
\left(z,pz\right) \in X \mapsto \left(z,0\right) \in X_0
\] 
where \(X_0\) is the subspace associated to some Hopf fibration. 
This map takes complex lines to complex lines, and therefore identifies the complex structures on the 2-planes belonging to \(X\) with those on \(X_0\).
\begin{center}
\begin{tikzpicture}
\fill[gray!5] (0,0) rectangle (1,2);
\foreach \i in {1,...,6}
{
	\draw[gray!50] ({.2*(\i-1)},0) -- ({.2*(\i-1)},2);
}
\draw[thick] (.5,0) .. controls (.8,0.25) and (.2,1.75) .. (.6,2) node[left] {\(X\)};
\draw[gray!50] (0,-.5) -- (1,-.5) node[right,black] {\(\operatorname{Gr}_{\mathbb{C}}\left(n+1,2n+2\right)\)};
\node[right,black] at (1,1) {\(\mathbb{CP}(\mathcal{U})\)};
\end{tikzpicture}
\end{center}

\newcommand{\linedone}{\\[20pt]}

\section{Further remarks}

 For applications to pseudoholomorphic
curves, one would like to define a concept
of totally real subspace \(R \subset V\),
\(\dim_{\R{}} R = n+1\). This should be
precisely a subspace of \(V\) which,
thought of as a great \(n\) sphere in 
\(S^{2n+1}\), has no great circles in
it from our fibration. Generic \(R\) should
have this property.

The general story of great sphere fibrations (the topological Blaschke theory) can probably be studied as follows: 
each Hopf fibration 
\[
\begin{tikzcd}[ampersand replacement=\&]
S^k \arrow{r} \& S^N \arrow{d} \\
           \& X_{\text{Hopf}}
\end{tikzcd}
\]
has a symmetry group \(\Gamma \subset \SL{N+1,\R{}}\), and a point of \(X_{\text{Hopf}}\) has an isotropy subgroup \(\Gamma_0 \subset \Gamma\).
This gives an embedding of homogeneous \(\SL{N+1,\R{}}\) spaces
\[
X_{\text{Hopf}} = \Gamma/\Gamma_0 \subset \SL{N+1,\R{}}/\Gamma_0.
\]
Then \(\SL{N+1,\R{}}\) acts on this picture to move the embedded submanifold
around in a fibration, with base \(\SL{N+1,\R{}}/\Gamma\).
\[
\begin{tikzcd}[ampersand replacement=\&]
\Gamma/\Gamma_0 \arrow{r} \& \SL{N+1,\R{}}/\Gamma_0 \arrow{d} \\
                       \& \SL{N+1,\R{}}/\Gamma.
\end{tikzcd}
\]
Now any great sphere fibration
\[
\begin{tikzcd}[ampersand replacement=\&]
S^k \arrow{r} \& S^N \arrow{d} \\
           \& X
\end{tikzcd}
\]
probably gives, via moving frame calculations similar to those above, an embedding
\[
X \subset \SL{N+1,\R{}}/\Gamma_0
\]
which is ``close'' to vertical.
Locally trivializing the fiber bundle
\[
\begin{tikzcd}[ampersand replacement=\&]
\Gamma/\Gamma_0 \arrow{r} \& \SL{N+1,\R{}}/\Gamma_0 \arrow{d} \\
                       \& \SL{N+1,\R{}}/\Gamma
\end{tikzcd}
\]
in some ``nice'' way, one will probably find that this gives an explicit diffeomorphism \(X \to X_{\text{Hopf}}\).

Taking \(M\) a manifold and \(SM=\left(TM \backslash 0\right)/\R{+}\), the tangent sphere bundle, one could consider a great sphere fibration of each sphere in \(SM\). 
Call this a \emph{Blaschke system}. 
We can interpret such a system as a first order system of partial differential equations, so that a solution to such an equation is an immersed submanifold of \(M\) whose tangent spaces project via \(TM \backslash 0 \to SM\) to be great spheres belonging to our
fibration. 
If there are involutive differential equations constructible in this way, with a suitable notion of taming, then this will provide a theory of pseudoquaternionic curves and pseudo-octave curves, and perhaps a theory of Gromov--Witten invariants for hyper-K{\"a}hler and octavic spaces.

For applications to elliptic partial differential equations, the most important result one would like to prove is probably the existence of a taming symplectic structure. 
Here that means an element \(\omega \in \Lm{2}{V^*}\) so that \(\omega > 0\) on each 2-plane in \(X\). 
In local complex coordinates \(z,w\) on \(V\), we can write \(2\)-planes as
\[
dw^i = p^i dz + q^i d\bar{z}
\]
and then the symplectic form
\[
\omega = \frac{\sqrt{-1}}{2} \left ( dz \wedge d \bar{z}
+ dw^i \wedge dw^{\bar{i}} \right )
\]
becomes on that 2-plane
\[
\omega = \frac{\sqrt{-1}}{2} \left ( 1 + |p|^2 - |q|^2 \right )
dz \wedge d\bar{z}.
\]
Therefore the null 2-planes for \(\omega\) form a real hyperquadric in \(\Gro{2}{V}.\)
When we look up in \(\Sato\) we find two more dimensions to the space, and to the
subspace of null 2-planes. 
It is still a real hyperquadric. 
The problem is then to show that this hypersurface does not intersect the base manifold of a great circle fibration, after suitable linear transformation.

In studying families of great circle fibrations, it would be  helpful to have a retraction from the space of great circle fibrations to the space of complex structures, i.e. Hopf fibrations. 
Such a retraction is most likely to be found using a parabolic heat flow.

\begingroup
\begin{table}
\[
\begin{array}{@{}lll@{}}
\toprule
\text{Right principal bundle} & 
\text{Semibasic 1-forms} & 
\text{Structure group}
\\ \midrule
B \to S^{2n+1} & \Omega^p_0, \Omega^0_{\bar{0}} + \Omega^{\bar{0}}_0 
& 
\Gamma_0 = 
\left \{
\begin{pmatrix}
g^0_0 & g^0_q \\
0 & g^p_q
\end{pmatrix}
; g^0_0 \in \R{}\right \}
\linedone
B \to X & \Omega^p_0 & 
\Gamma_1 = 
\begin{pmatrix}
g^0_0 & g^0_q \\
0 & g^p_q 
\end{pmatrix}
\linedone
G \to S^{2n+1} &
\Omega^p_0,\Omega^0_{\bar{0}} &
G_0 = \begin{pmatrix}
g^0_0 & g^0_1 & g^0_j \\
0 & g^1_1 & g^1_j \\
0 & g^i_1 & g^i_j
\end{pmatrix} 
\linedone
G \to \Gro{2}{V} &
\Omega^p_0, \Omega^p_{\bar{0}} &
\Gcircle = 
\begin{pmatrix}
g^0_0 & g^0_1 & g^0_j \\
g^1_0 & g^1_1 & g^1_j \\
0 & 0 & g^i_j
\end{pmatrix} 
\linedone
G \to \Cstrucs{V} &
\Omega^0_{\bar{0}}, 
\Omega^p_{\bar{0}},\Omega^0_{\bar{p}},\Omega^p_{\bar{q}} &
\GL{V,J_0}
\linedone
G \to \Sato & 
\Omega^0_{\bar{0}}, \Omega^p_0, \Omega^p_{\bar{0}} &
\begin{pmatrix}
g^0_0 & -g^1_0 & g^0_j \\
g^1_0 & g^0_0 & g^1_j \\
0 & 0 & g^i_j
\end{pmatrix} 
\linedone
G \to \M &
\Omega^p_0, \Omega^{\bar{p}}_0,
\Omega^{\bar{0}}_0, \Omega^{\bar{0}}_p,
\Omega^{\bar{p}}_q
&
\Gamma_0
\\ \bottomrule
\end{array}
\]
\caption{Actual structure groups are intersections of the above ones with \(G=\SL{V}\), \(V=\R{2n+2}\).}
\end{table}
\endgroup

\begin{table}
\begin{align*}
\mu,\nu,\sigma&=1,\dots,2n+1 \\
i,j,k&=2,\dots,2n+1 \\
p,q,r&=1,\dots,n \\
P,Q,R&=0,\dots,n.
\end{align*}
\caption{Index conventions}
\end{table}

%\nocite{*}
\bibliographystyle{amsplain}
\bibliography{greatCircle}

\providecommand{\bysame}{\leavevmode\hbox to3em{\hrulefill}\thinspace}
\providecommand{\MR}{\relax\ifhmode\unskip\space\fi MR }
% \MRhref is called by the amsart/book/proc definition of \MR.
\providecommand{\MRhref}[2]{%
  \href{http://www.ams.org/mathscinet-getitem?mr=#1}{#2}
}
\providecommand{\href}[2]{#2}
\begin{thebibliography}{10}

\bibitem{Besse:1978}
Arthur~L. Besse, \emph{Manifolds all of whose geodesics are closed},
  Springer-Verlag, Berlin, 1978, With appendices by D. B. A. Epstein, J.-P.
  Bourguignon, L. B\'erard~Bergery, M. Berger and J. L. Kazdan. \MR{80c:53044}

\bibitem{Bryant:2014}
Robert~L. Bryant, \emph{Notes on exterior differential systems}, ArXiv e-prints
  (2014), Notes from the December 2013 workshop Exterior Differential Systems
  and Lie Theory at the Fields Institute in Toronto.

\bibitem{Clelland:2016}
Jeanne Clelland, \emph{From {F}renet to {C}artan: The method of moving frames},
  Amer. Math. Soc., 2016.

\bibitem{DNF:1985}
B.~A. Dubrovin, A.~T. Fomenko, and S.~P. Novikov, \emph{Modern
  geometry---methods and applications. {P}art {I}{I}}, Springer-Verlag, New
  York, 1985, The geometry and topology of manifolds, Translated from the
  Russian by Robert G. Burns. \MR{86m:53001}

\bibitem{GluckWarnerYang:1983}
Herman Gluck, Frank Warner, and C.~T. Yang, \emph{Division algebras, fibrations
  of spheres by great spheres and the topological determination of space by the
  gross behavior of its geodesics}, Duke Math. J. \textbf{50} (1983), no.~4,
  1041--1076. \MR{85i:53047}

\bibitem{GluckWarner:1983}
Herman Gluck and Frank~W. Warner, \emph{Great circle fibrations of the
  three-sphere}, Duke Math. J. \textbf{50} (1983), no.~1, 107--132.
  \MR{84g:53056}

\bibitem{Gromov:1986}
Mikhael Gromov, \emph{Partial differential relations}, Springer-Verlag, Berlin,
  1986. \MR{90a:58201}

\bibitem{Hartshorne:1967}
Robin Hartshorne, \emph{Foundations of projective geometry}, Lecture Notes,
  Harvard University, vol. 1966/67, W. A. Benjamin, Inc., New York, 1967.
  \MR{0222751}

\bibitem{Kobayashi/Nagano:1964}
Shoshichi Kobayashi and Tadashi Nagano, \emph{On projective connections}, J.
  Math. Mech. \textbf{13} (1964), 215--235. \MR{0159284}

\bibitem{McCleary:1985}
John McCleary, \emph{User's guide to spectral sequences}, Publish or Perish
  Inc., Wilmington, DE, 1985. \MR{87f:55014}

\bibitem{McKay:2001}
Benjamin McKay, \emph{Analogues of complex geometry}, eprint, 2001.

\bibitem{McKay:2001b}
\bysame, \emph{Dual curves and pseudoholomorphic curves}, eprint, 2001.

\bibitem{Reznikov:1985}
A.~G. Reznikov, \emph{Blaschke manifolds of the projective plane type},
  Funktsional. Anal. i Prilozhen. \textbf{19} (1985), no.~2, 88--89.
  \MR{800932}

\bibitem{Sato:1984}
Hajime Sato, \emph{On topological {B}laschke conjecture. {I}. {C}ohomological
  complex projective spaces}, Geometry of geodesics and related topics (Tokyo,
  1982), North-Holland, Amsterdam, 1984, pp.~231--238. \MR{86c:53025a}

\bibitem{Sikorav:2000}
Jean-Claude Sikorav, \emph{Dual elliptic structures on {$\CP{2}$}}, eprint,
  2000.

\bibitem{Yang:1990}
C.~T. Yang, \emph{Smooth great circle fibrations and an application to the
  topological {B}laschke conjecture}, Trans. Amer. Math. Soc. \textbf{320}
  (1990), no.~2, 507--524. \MR{91e:55025}

\bibitem{Yang:1993}
\bysame, \emph{On smooth great circle fibrations of a round sphere},
  Differential geometry (Shanghai, 1991), World Sci. Publishing, River Edge,
  NJ, 1993, pp.~301--309. \MR{96f:53052}

\end{thebibliography}

\end{document}